\documentclass[10pt, reqno]{amsart}
\usepackage{amsmath, amssymb, amscd, amsthm, amsfonts, amsbsy, bbm, mathrsfs, bm, linearb, upgreek, mathtools, extarrows}
\usepackage[pdftex, final]{graphicx}
\usepackage{indentfirst}
\usepackage{setspace}
\usepackage{lipsum}
\usepackage{xcolor}
\usepackage{float}
\usepackage{tikz-cd}
\usepackage{enumitem}
\usepackage[all,cmtip]{xy}
\usepackage[mathscr]{euscript}
\usepackage{wrapfig}
\graphicspath{{./pics/}}
\usepackage[parfill]{parskip}
\usepackage{caption}
\usepackage{subcaption}
\usepackage{dcolumn}
\usepackage{color}
\usepackage[pdfauthor={},%
            pdftitle={},%
            pdftex]{hyperref}
\hypersetup{
	colorlinks,
	citecolor=red,
	filecolor=black,
	linkcolor=blue,
	urlcolor=blue
}
\usepackage[utf8]{inputenc}
\usepackage[english]{babel}

\DeclareSymbolFont{symbolsC}{U}{txsyc}{m}{n}
\SetSymbolFont{symbolsC}{bold}{U}{txsyc}{bx}{n}
\DeclareMathSymbol{\multimapdotboth}{\mathrel}{symbolsC}{22}

\usepackage[left=3.00cm, right=3.00cm, top=2.00cm, bottom=3.00cm]{geometry} 
\newcommand{\Z}{\mathbb{Z}}
\newcommand{\bgd}{\begin{displaymath}}
\newcommand{\edd}{\end{displaymath}}
\newcommand{\hf}{\hspace*{0.5cm}}

\theoremstyle{plain}

\newtheorem{theorem}{Theorem}[section]
\newtheorem{prop}[theorem]{Proposition}
\newtheorem{cor}[theorem]{Corollary}
\newtheorem{lemma}[theorem]{Lemma}
\theoremstyle{remark}

\def\D{\mathcal D}
\def\B{\mathcal B}
\theoremstyle{definition}
\newtheorem{defn}[theorem]{Definition}
\newtheorem{rem}[theorem]{Remark}
\newtheorem{eg}[theorem]{Example}
\usepackage[most]{tcolorbox}

\newtheorem{innercustomgeneric}{\customgenericname}
\providecommand{\customgenericname}{}
\newcommand{\newcustomtheorem}[2]{%
  \newenvironment{#1}[1]
  {%
   \renewcommand\customgenericname{#2}%
   \renewcommand\theinnercustomgeneric{##1}%
   \innercustomgeneric
  }
  {\endinnercustomgeneric}
}
\newcustomtheorem{customthm}{Theorem}
\numberwithin{equation}{section}

\usepackage{biblatex} 
\begin{document}

\vspace*{-0.3cm}
\title{Graph polynomial for vertex-colored embedded graphs}

\author{Somnath Basu}
\address{Department of Mathematics \& Statistics, Indian Institute of Science Education \& Research, Mohanpur 741246, West Bengal, India}
\email{somnath.basu@iiserkol.ac.in}

\author{Dhruv Bhasin}
\address{Department of Mathematical Sciences, Indian Institute of Science Education and Research Pune, Maharashtra - 411008, India}
\email{bhasin.dhruv@students.iiserpune.ac.in}

\author{Siddhartha Lal}
\address{Department of Physical Sciences, Indian Institute of Science Education \& Research, Mohanpur 741246, West Bengal, India}
\email{slal@iiserkol.ac.in}

\author{Siddhartha Patra}
\address{Centre for Quantum Engineering, Research and Education, TCG CREST, Kolkata 700091, India}
\email{siddhartha.patra@tcgcrest.org}

\begin{abstract}
Topological invariants arising from planar graphs have been used to study topological entanglement entropy in physics. Motivated by a prior work by some of us, we study finite (vertex) colored graphs embedded in surfaces by associating a polynomial to it, called colored island polynomial. The tools used in developing a theory of such graph polynomials are algebraic topological while the polynomial itself is inspired from ideas arising in physics. We also analyze a special case of these polynomials for uncolored embedded graphs, called island polynomials. We describe the changes in these polynomials under several (standard as well as non-standard) graph theoretic operations. The collection of colored island polynomials is defined to be the structure set and has interesting properties and applications. We conclude with several applications of this polynomial including detection of trees and cycles, and important connections of this polynomial with topological entanglement entropy. 
\end{abstract}

\maketitle

\renewcommand{\baselinestretch}{0.6}\normalsize
\tableofcontents
\renewcommand{\baselinestretch}{1.0}\normalsize


\section{Introduction}

\hf In algebraic graph theory, invariants of graphs which are polynomials form natural objects of interest. The importance stems from their applicability and the underlying power of algebra to efficiently package information. Among the commonly known polynomial invariants, one of the most famous is the Tutte polynomial. It is the most general graph invariant defined by deletion-contraction recursion. However, it does not naturally generalize to embedded graphs in surfaces. There are stronger and more modern invariants like Bollob\'{a}s-Riordan polynomial for ribbon graphs \cite{BolRio01, BolRio02}, Las Vergnas polynomial \cite{Ver81}, \cite{EMMof15} and Krushkal polynomial \cite{Kru11}. It is known (see Theorem 27.10 of \cite{EmMof22}) that the Krushkal polynomial is powerful enough to encode these other polynomials (see figure 27.1 in chapter 27 of \cite{EmMof22}). The Krushkal polynomial and its modern generalized form due to Butler \cite{But18} is one of the most powerful polynomial invariant of graphs embedded in surfaces. It also enjoys a duality relation for cellularly embedded graphs. \\
\hf Given a graph $\Gamma$ with a vertex-coloring $\mathfrak{c}$ (not necessarily proper) and an embedding (not necessarily cellular) $\sigma:\Gamma\hookrightarrow \Sigma$ inside a surface $\Sigma$, we define a polynomial $\beta_{(\Gamma,\sigma,\mathfrak{c})}$ (definition \ref{defn:beta-c}), in one indeterminate $x$, called the {\it colored island polynomial}. It is built out of signed counts of the number of faces of all color-induced subgraphs. There is another special case of this polynomial $\beta_{(\Gamma,\sigma)}$ (definition \ref{defn:beta}), called the {\it island polynomial}, where the vertex-coloring is injective. The origin of these ideas and the subsequent developments of it, as presented in this article, are heavily inspired by an earlier work \cite{SiddharthaTEE2021} by some of us. In that joint work, we had shown that the topological entanglement entropy (TEE) of a topologically ordered phase, a quantity of importance in physics, arises from a polynomial assigned to the underlying planar graph. In particular, the TEE was related to $\beta_{(\Gamma,\sigma)}(-1)$, where $\sigma$ was a planar embedding. The TEE was also shown to be robust under several types of deformations. \\
\hf In this article, we prove several results about the colored island polynomial $\beta_{(\Gamma,\sigma,\mathfrak{c})}$ which indicate the non-trivial nature and the potential utility of this invariant. A key ingredient is the {\it structure set} of $\Gamma$, defined as the collection of all colored island polynomials across all possible embeddings and colorings (see definition \ref{defn:str-set}); it is denoted by $\mathscr{S}(\Gamma)$. The structure set has important graph-detection properties (see \S \ref{subsec:str-set} and \S \ref{subsec:cycle-tree}). The first main result, which is actually a combination of Theorem \ref{thm:Dnm}, Theorem \ref{thm:tree-detection} and Theorem \ref{thm:cycle-detect}, is the following. 
\begin{customthm}{A}
{\it The island polynomial has the following properties:\\
\textup{(a)} If $\Gamma$ is a graph such that the island polynomial $\beta_{(\Gamma,\sigma)}(x)=(1+x)^{n-1}+(n-1)(1+x)^{n-2}$ for any embedding, then $\Gamma$ is a tree on $n$ vertices. Moreover, if $\Gamma$ is a connected planar graph such that the island polynomial is given by
\bgd
{\beta}_{\Gamma}(x)=a(1+x)^{n-1}+b(1+x)^{n-2}
\edd
for $a\in \mathbb{N}$, then $\Gamma$ is obtained from a tree on $n$ vertices by adding $(a+b-n)$ loops and $(n-b-1)$ parallel edges. \\
\textup{(b)} The island polynomial for $(C_n,\sigma)$ is given by
\begin{equation*}
{\beta}_{(C_n,\sigma)}(x)=\left\{\begin{array}{rl}
n(1+x)^{n-2}+2x^{n-1} & \textup{if $\Sigma-\sigma(C_n)$ has $2$ components}\\
n(1+x)^{n-2}+x^{n-1} & \textup{otherwise}.
\end{array}\right.
\end{equation*}
If $\Gamma$ is a graph such that the island polynomials satisfy
$$\left\{\beta_{(\Gamma,\sigma)}(x)\,|\,\sigma:\Gamma\hookrightarrow \Sigma\right\}\subseteq\{x^{n-1}+n(1+x)^{n-2},2x^{n-1}+n(1+x)^{n-2}\}=\mathscr{S}(C_n)$$
then $\Gamma$ is a cycle on $n$ vertices. }
\end{customthm}
\vspace*{-0.2cm}
The second main result (see Corollary \ref{cor:collapse-tree} and Theorem \ref{thm:tree-cycle-color-detect}), analogous to Theorem A above but for colored island polynomials, involves proper vertex-coloring. 
\begin{customthm}{B}
{\it The colored island polynomial for any graph $(\Gamma,\sigma,\mathfrak{c})$ is the same as that for $(\Gamma/T,\overline{\sigma},\mathfrak{c})$, where $T$ is a uni-colored tree. In particular, any cycle or tree can be iteratively reduced (via tree collapses) to obtain a cycle or tree (respectively) which is proper vertex-colored. \\
\textup{(a)} Let $(T_n,\mathfrak{c})$ be proper vertex-colored tree on $n$ vertices and $N$ colors. Then the colored island polynomial is given by
\begin{equation*}
{\beta}_{(T_{n},\sigma,\mathfrak{c})}(x)= (1+x)^{N-1}+(n-1)(1+x)^{N-2}.
\end{equation*}
\textup{(b)} Let $(C_n,\mathfrak{c})$ be proper vertex-colored cycle on $n$ vertices and $N$ colors. Then the colored island polynomial is given by
\begin{equation*}
{\beta}_{(C_{n},\sigma,\mathfrak{c})}(x)=\alpha x^{N-1}+n(1+x)^{N-2}.
\end{equation*}
where $\alpha$ is the number of components of $\Sigma-\sigma(C_n)$.}
\end{customthm}
A graph can have several inequivalent vertex-colored embeddings inside the same surface, leading to different colored island polynomials. On the one hand, the theory yields results requiring very few hypotheses when we are dealing with planar graphs. On the other hand, since any finite graph can be embedded in surfaces with sufficiently high genus, no finite graphs are excluded from the ambit of our approach. \\
\hf The colored island polynomial $\beta_{(\Gamma,\sigma,\mathfrak{c})}$ is essential in studying how the island polynomial $\beta_{(\Gamma,\sigma)}$ changes when we do an edge contraction or an edge subdivision. We have analyzed how $\beta_{(\Gamma,\sigma,\mathfrak{c})}$ transforms under basic operations like disjoint unions (Proposition \ref{prop:disjoint-union}), adding a pendant edge (Propositon \ref{prop:Gamma-app}), edge subdivision (Proposition \ref{prop:subdivide}), adding a loop (Proposition \ref{prop:self-loop}), and deleting a bridge (Proposition \ref{prop:bridge-del}). For instance, the effect of edge contraction for a trivial coloring is given below (see Theorem \ref{thm:collapse}).
\begin{customthm}{C}
{\it Let $(\Gamma,\sigma,\mathfrak{c})$ be trivially colored. Then the polynomial for $\Gamma/e$ is given by
\begin{equation*}
{\beta}_{(\Gamma,\sigma)}(x)=x{\beta}_{(\Gamma/e,\overline{\sigma})}(x)+{\beta}_{(\Gamma-v,\sigma)}(x)+{\beta}_{(\Gamma-w,\sigma)}(x)-(1+x){\beta}_{(\Gamma-\{v,w\},\sigma)}(x),
\end{equation*}
where $e$ joins $v$ and $w$.}
\end{customthm}
We analyze other transformations like adding a parallel edge or adding a clean shortcut (see definition \ref{defn:csc}). Let $\Gamma\cup e$ be obtained from a connected graph $\Gamma$ via an edge $e$ joining distinct non-adjacent vertices $v$ and $w$. Then $e$ is called a {\it clean shortcut} if there exists a path $P_k$ (with $k\geq 3$) in $\Gamma$ joining $v$ and $w$ such that every non-pendant vertex of $P_k$ has valency two; such a $P_k$ is called a {\it clean path} with respect to $e$.
\begin{customthm}{D}
{\it \textup{(a)} Let $(\tilde{\Gamma},\tilde{\sigma},\mathfrak{c})$ be a graph on $n>2$ vertices, obtained from $(\Gamma,\sigma,\mathfrak{c})$ by adding a parallel edge $e'$, joining vertices $v_1$ and $v_2$. Let $N$ be the number of colours and $S=\{\mathfrak{c}(v_1),\mathfrak{c}(v_2)\}$. \\
\textup{(i)} If $f_{\tilde{\sigma}}(\mathscr{F}_S(\tilde{\Gamma},\mathfrak{c}))=f_{\sigma}(\mathscr{F}_S(\Gamma,\mathfrak{c}))+1$ then
\begin{equation*}
\beta_{(\tilde{\Gamma},\tilde{\sigma},\mathfrak{c})}(x)=\beta_{(\Gamma,\sigma,\mathfrak{c})}(x)+x^{|S|-1}(1+x)^{N-|S|}.
\end{equation*}
\textup{(ii)} If $\tilde{\sigma}(e')$ does not disconnect $(\Sigma-\sigma(\Gamma))\cup\{v_1,v_2\}$, then ${\beta}_{(\tilde{\Gamma},\tilde{\sigma},\mathfrak{c})}(x)={\beta}_{(\Gamma,\sigma,\mathfrak{c})}(x)$. \\
\textup{(b)} Let $(\Gamma,\sigma,\mathfrak{c})$ be a connected graph with a clean shortcut $e$ joining $v_1$ to $v_2$ with $P_k$ a clean path. Let $\Gamma_1:=(\Gamma-P_k)\cup e$ and assume that $|V(\Gamma_1)|\geq 3$. Let $\mathfrak{c}$ be a vertex-coloring of $\Gamma$ such that all the non-pendant vertices of $P_k$ receive distinct colors different from the $N$ colors in $\mathfrak{c}(\Gamma_1)$. Let $\Gamma'$ be the graph induced on the vertex set $V(P_k)$. If $f_{\sigma}(\Gamma')=f_{\sigma}(\Gamma'-P_k)+1$ or $\sigma(e)$ does not disconnect $(\Sigma-\sigma(\Gamma-P_k))\cup\{v_1,v_2\}$, then $\beta_{(\Gamma,\sigma,\mathfrak{c})}$ can be computed explicitly (see \eqref{eqn:clean-sc1}, \eqref{eqn:clean-sc2} in Proposition \ref{prop:clean-sc}) in terms of $\beta_{(\Gamma_1,\sigma,\mathfrak{c})}(x)$, the size of $S=\{\mathfrak{c}(v_1),\mathfrak{c}(v_2)\}$ and $N$.}
\end{customthm}
\hf We have provided a comparison (see \S \ref{subsec:Krushkal}) with other known state-of-the-art polynomials and how our polynomials differ from those. To our understanding, the colored island polynomial does not seem to be obtainable from these known polynomials. Moreover, the colored island polynomial $\beta_{(\Gamma,\sigma,\mathfrak{c})}$, should ideally be interpreted as a polynomial with three input parameters: one indeterminate $x$, a vertex-coloring $\mathfrak{c}$ and an embedding $\sigma$. This makes it structurally different from other graph polynomials. \\
\hf Graph polynomials have been traditionally used to study properties of graphs, including adjacency (characteristic polynomial), coloring (Birkhoff's chromatic polynomial \cite{Big93}), Euler tours, rank polynomial \cite{GoSo01} and more. The Ihara zeta function \cite{Iha66} is used in the study of free groups, spectral graph theory and symbolic dynamics. The Tutte polynomial \cite{Tut54, Tut67, Bol98} encodes information about connectivity of induced subgraphs. Historically, there have been connections between work in physics and graph polynomials. As as example, after Potts' work \cite{Pot52} on partition function of certain models in statistical mechanics in 1952, Fortuin and Kasteleyn \cite{FoKa72} found connections between the Tutte polynomial and their work on random cluster model, a generalisation of the Potts model. In their work in 2010, Chang and Shrock \cite{ChSh10} had defined a polynomial $Ph(G, q, w)$, which generalizes the chromatic polynomial $P(G,q)$. A recent study \cite{SiddharthaTEE2021} of the entanglement entropy of a topologically ordered state of quantum matter sheds light on the connection between multipartite quantum information measure and questions in graph theory.\\
\hf As indicated earlier, the value of the island polynomial at $x=-1$, i.e., $\beta_{(\Gamma,\sigma)}(-1)$, plays an important role and is called the {\it island count} of the graph $(\Gamma,\sigma)$. If $\Gamma$ has $N$ vertices, then from the perspective of physics, the island count is related to the $N$-partite information among $N$ subsystems in a topologically ordered ground state~\cite{wen2013,wen2017colloquium,wen2019review}. We map a collection of $N$ subsystems to a graph $\Gamma$ by representing each subsystem with a vertex and the connectivity between two subsystems by an edge connecting the two corresponding vertices. Thus, holes in the subsystem are represented by cycles in the graph.\\
\hf Apart from the main results, we have several applications (see \eqref{eqn:disjoint-cor}, Proposition \ref{prop:tree-cycle-gr} and discussions in \S \ref{subsec:TEE} and \S \ref{subsec:Euler}) which have been collected in the following. 
\begin{customthm}{E}
{\it \textup{(a)} The island count vanishes for a disjoint union of graphs as well as for a wedge sum of graphs. The island count also vanish for planar graphs built out of a tree with at least $3$ vertices by applying a finite sequence of one of two operations: we add a loop (or a parallel edge) and subdivide this loop (or this parallel edge) finitely many times. \\
\textup{(b)} The vanishing of the multipartite information measure for a planar collection of subsystems is equivalent to the vanishing of the island count for the associated graph.\\
\textup{(c)} For a cellularly embedded graph $(\Gamma,\sigma)$ in a compact surface $\Sigma$, the alternating sum of the island count over all subgraphs, on at least $3$ vertices, is $\chi-2f$, where $f$ is the number of components of $\Sigma-\sigma(\Gamma)$ and $\chi$ is the Euler characteristic of $\Sigma$, i.e.,
$$\chi-2f=(-1)^n\beta_{(\Gamma,\sigma)}(-1) + \sum_{\Gamma'\in \mathscr{F}_{n-1}(\Gamma)}(-1)^{n-1}\beta_{(\Gamma',\sigma)}(-1)+\cdots+\sum_{\Gamma'\in \mathscr{F}_{3}(\Gamma)}(-1)^{3}\beta_{(\Gamma',\sigma)}(-1)
$$
where $\mathscr{F}_j(\Gamma)$ denotes the collection of all color-induced subgraphs on $j$ colors.}
\end{customthm}
\vspace*{-0.3cm}
We may interpret Theorem E (c) as the statement that $\chi-2f$ is an emergent feature of the island counts associated to any embedded graph.\\
\hf Graph theory has gradually become an essential part of computer science \cite{deo2017graph,riaz2011applications,majeed2020graph} as well as network analysis \cite{network_graph_barnes1983,network_graph_bernhardt2015network,network_graph_derrible2009network,network_graph_scott1988social} (see also \cite{battiston2021} for a recent review on applications in physics). We firmly believe that the theory and results presented here will be useful to the community working on network analysis and applications of quantum information theory to quantum condensed matter physics (such as TEE), apart from its use within the graph theory community.\\[-0.2cm]

\noindent \textsc{Organization of the paper}. In \S \ref{sec:invariant}, we define the colored island polynomial and explore some of its basic properties in \S \ref{subsec:edge-loop}-\S \ref{subsec:split_edge}. Two other properties are analyzed in \S \ref{subsec:arrange_nearestneighbor}-\S \ref{subsec:arrange_shortcircuit}, the latter being the notion of {\it shortcuts} which we introduce in this article. The structure set, which plays an important role in some of the subsequent results, is introduced in \S \ref{subsec:str-set}. A comparison of the colored island polynomial (in terms of properties and identities it satisfy) with the Bollob\'{a}s-Riordan polynomial and the modern version of Krushkal polynomial is presented in \S \ref{subsec:Krushkal}. In \S \ref{sec:appl}, we give several applications of the theory presented here. In \S \ref{subsec:cycle-tree} we prove detection results (see Theorem A and B earlier). In \S \ref{subsec:tree-cycle-graphs}, we construct a class of planar graphs (called tree-cycle graphs) for which the island count vanishes. In \S \ref{subsec:TEE} we recall topological entanglement entropy (TEE) and show that the island count is a multiplicative factor in the TEE of a multipartite planar system. We show in \S \ref{subsec:Euler} that for cellularly embedded graphs, the alternating sum of the island counts over induced subgraphs on $3$ or more vertices is $\chi-2f$. In \S \ref{subsec:pair-of-pants} we discuss how the island polynomial changes under an important transformation which changes the combinatorial topology of the embedded graph.\\[-0.2cm]
	
\noindent \textsc{Acknowledgments.} The authors would like to thank Kaneenika Sinha, Moumanti Podder, Soumya Bhattacharya and Niranjan Balachandran for initial discussions on this topic. S. Basu would like to thank SERB for support through MATRICS grant MTR/2017/000807. D. Bhasin would like to acknowledge NBHM grant 0203/2/2021/RD-II/3033. S. Lal thanks the SERB, Govt.
of India for funding through MATRICS grant MTR/2021/000141 and Core Research Grant CRG/2021/000852. S. Patra would like to thank CSIR and IISER Kolkata for funding through a research fellowship.\\



\section{Invariants for vertex-colored finite graphs}
\label{sec:invariant}

\hf A graph is denoted by $\Gamma=(V, E)$, where $V:=V(\Gamma)$ is the set of vertices and $E:=E(\Gamma)$ is the set of edges. We shall consider graphs where multiple edges between vertices are allowed. In particular, we allow for loops as valid edges. Usually, such graphs are called \textit{multigraphs}, but we will refer to these as graphs in what follows. Thus, $E$ is a multiset and not just a subset of $V\times V$. If $\Gamma$ is a finite graph, i.e., both $V$ and $E$ are finite sets, then $|V|$ and $|E|$ will denote the number of vertices and edges, respectively. A subgraph $\Gamma'$ with a vertex set $S$ is called \textit{induced} if any edge in $\Gamma$ joining two vertices in $S$ is also in the subgraph. We will typically be dealing with non-trivial induced subgraphs, i.e., an induced subgraph where the vertex set is neither $\varnothing$ nor $V(\Gamma)$. Note that a graph can be given a natural (quotient) topology by identifying edges with $[0,1]$ and subsequently identifying endpoints of intervals according to the adjacency relations in $\Gamma$. \\
\hf We also need the notion of vertex-colored graphs. 
\begin{defn}
A \textit{vertex-coloring} of a graph $\Gamma$ is a map $\mathfrak{c}:V \to \mathcal{C}$ to a set of colors $\mathcal{C}$. A \textit{proper vertex-coloring} is a vertex-coloring in which no two adjacent vertices receive the same color.
\end{defn}
We shall assume, from now on, that any vertex-coloring is surjective, unless otherwise mentioned. As we deal with graphs with loops, the more classical notion of proper vertex-coloring, where no two adjacent vertices are assigned the same color, is not relevant in our context. When the vertex-coloring $\mathfrak{c}$ is injective, then we say that the graph is trivially colored. When the vertex-coloring is constant, then we say that the graph is uni-colored.
\begin{figure}[!ht]
\centering
\includegraphics[trim={0 40 0 0}, clip, scale=0.25]{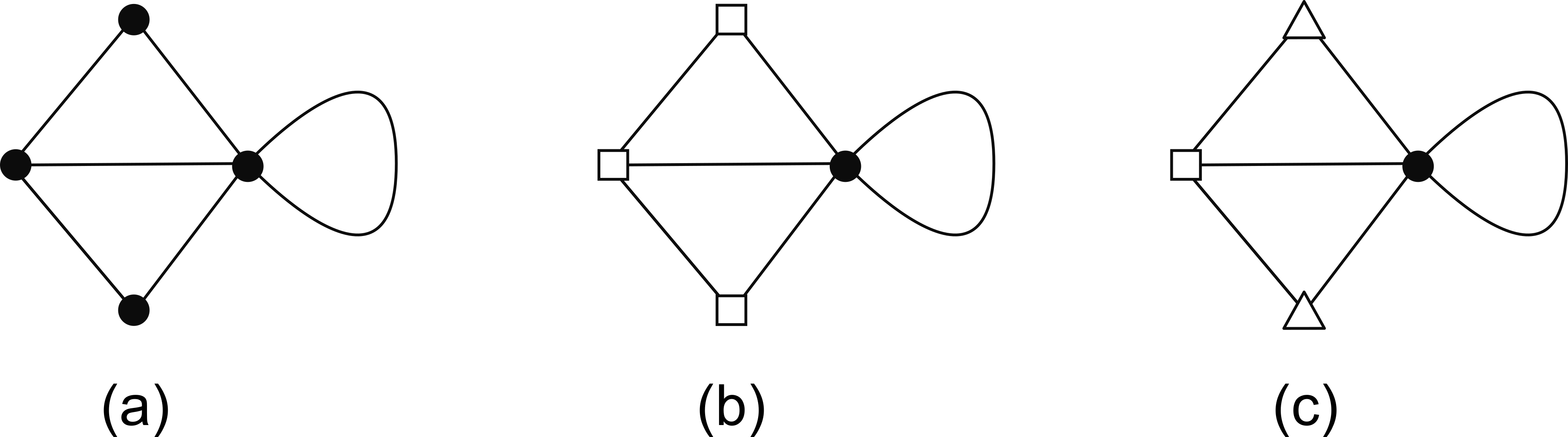}
\caption{Vertex-coloring (colors indicated by different shapes) of the same graph}
\label{fig:ver-col}
\end{figure}
\begin{defn}
Given a vertex-coloring $\mathfrak{c}:V(\Gamma)\to\mathcal{C}$, let $N$ denote the size of $\mathfrak{c}(V)$ (or the size of $\mathcal{C}$, as our vertex-coloring is assumed to be surjective). For $i=1,\ldots,N$, we define $\mathscr{F}_i(\Gamma,\mathfrak{c})$ to be the collection of induced subgraphs of $\Gamma'$ that satisfy\\
\hf (a) $\Gamma'$ requires $i$ colors;\\
\hf (b) if $c_1,\ldots,c_i$ are the colors for $\Gamma'$, then all vertices of color $c_1,\ldots,c_i$ are in $\Gamma'$.\\    
We call $\mathscr{F}_i(\Gamma,\mathfrak{c})$ the set of {\it colored island subgraphs} on $i$ colors. If $S\subseteq \mathcal{C}$ is of size $i$, then $\mathscr{F}_S(\Gamma,\mathfrak{c})$ denotes the unique element of $\mathscr{F}_i(\Gamma,\mathfrak{c})$ that is vertex-colored by $S$. By our notation, $\mathscr{F}_N(\Gamma,\mathfrak{c})=\mathscr{F}_\mathcal{C}(\Gamma,\mathfrak{c})=\{\Gamma\}$.
\end{defn}
\vspace*{-0.2cm}
\begin{figure}[!ht]
\centering
\includegraphics[scale=0.25]{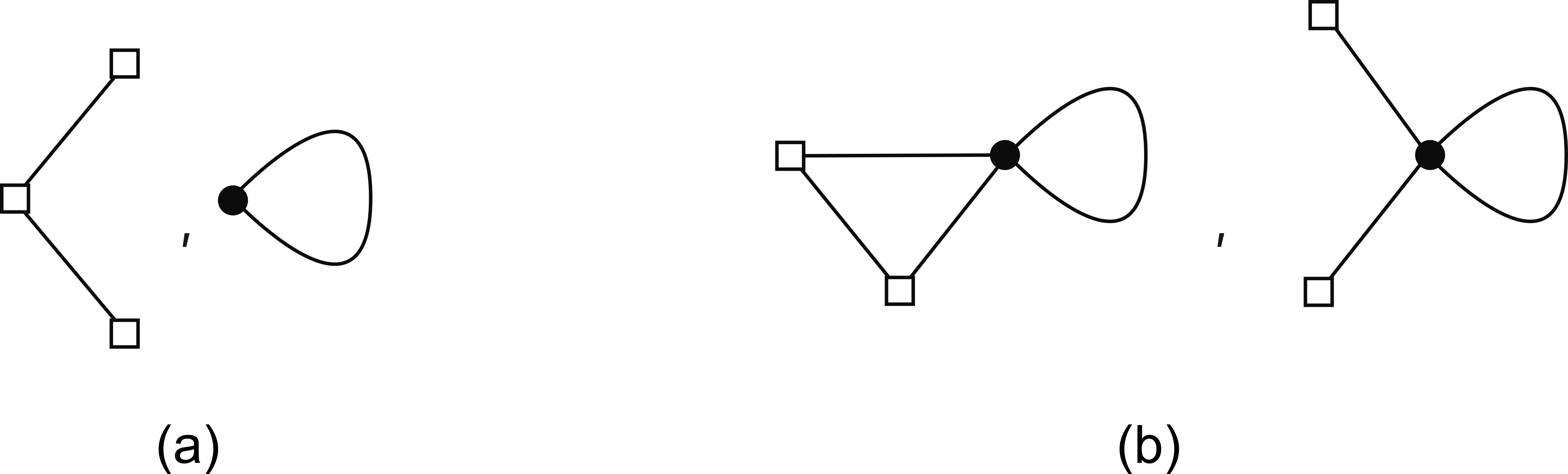}
\caption{(a) Colored island subgraphs $\mathscr{F}_1$ for fig \ref{fig:ver-col}(b); (b) graphs that are not in $\mathscr{F}_2$ for fig \ref{fig:ver-col}(b)}
\label{fig:color-is-sub}
\end{figure}
\begin{defn}[Colored island boundary count] \label{defn:beta-c}
Let $\sigma:\Gamma\hookrightarrow \Sigma$ be an embedding of $\Gamma$ in a surface $\Sigma$. Let $\mathfrak{c}:V \to \mathcal{C}$ be a vertex-coloring. For any connected subgraph $\Gamma'$ of $\Gamma$, let $\Sigma'\subset \Sigma$ be the path component containing $\sigma(\Gamma')$, Define $f_\sigma(\Gamma')$ to be the number of path components of $\Sigma' - \sigma(\Gamma')$. For a general subgraph $\Gamma'$ with components (or connected induced subgraphs) $\Gamma'_1,\ldots, \Gamma'_k$, we define
\bgd
f_\sigma(\Gamma'):=\sum_{i=1}^k f_\sigma(\Gamma'_i).
\edd
Define the $i^\textup{th}$ \textit{colored island boundary count} of $\Gamma$ with respect to the data $(\sigma,\mathfrak{c})$ to be
\begin{equation*}\label{eq:Di_col2}
\mathscr{D}_i(\Gamma,\sigma,\mathfrak{c})=\sum_{\mathscr{F}_i(\Gamma,\mathfrak{c})}  f_\sigma(\Gamma').
\end{equation*}
The \textit{colored island polynomial} for $(\Gamma,\sigma,\mathfrak{c})$ is defined to be 
\begin{equation}\label{eqn:col-isl-pol}
{\beta}_{(\Gamma,\sigma,\mathfrak{c})}(x):=\sum_{i=1}^{N}  \mathscr{D}_i(\Gamma,\sigma,\mathfrak{c})x^{i-1}
\end{equation}
where $N$ is the number of colors needed in the coloring $\mathfrak{c}$, i.e., $N=|\mathfrak{c}(V)|$.
\end{defn}
\begin{rem}
1. We do not require the embedding $\sigma:\Gamma\hookrightarrow \Sigma$ to be cellular, i.e., the faces of $\Sigma - \sigma(\Gamma)$ need not be homeomorphic to $2$-dimensional disks. \\
2. The colored island polynomial is the same for any trivial coloring. \\
3. If $\mathfrak{c}$ is a constant function, i.e., $(\Gamma,\sigma)$ is uni-colored, then 
\bgd
{\beta}_{(\Gamma,\sigma,\mathfrak{c})}(x)=\mathscr{D}_{1}(\Gamma,\sigma,\mathfrak{c})=\sum_{j=1}^k f_\sigma(\Gamma_j'),
\edd
where $\Gamma=\Gamma_1'\sqcup \cdots\sqcup \Gamma_k'$ is a decomposition of $\Gamma$ into its connected components.
\end{rem}
\begin{rem}\label{infinite}
The definition of colored island polynomial can be extended to a colored embedded infinite graph $\Gamma$ which is finitely colored if it satisfies the condition that $f_\sigma(\Gamma')$ is finite for any $\Gamma'\in \mathscr{F}_i(\Gamma,\mathfrak{c})$ and any positive integer $i$. In particular, this forces $\Gamma'$ to have finitely many components. 
\end{rem}
\begin{defn}\label{defn:beta}
The \textit{island polynomial} for $(\Gamma,\sigma)$ is defined to be the colored island polynomial for any trivial coloring, i.e., $\beta_{(\Gamma,\sigma)}:=\beta_{(\Gamma,\sigma,\mathfrak{c})}$ for any injective coloring $\mathfrak{c}:V\to\mathcal{C}$. The \textit{island count} for $(\Gamma,\sigma)$ is defined to be $\beta_{(\Gamma,\sigma)}(-1)$. 
\end{defn}
We shall use the notation
$$\mathscr{D}_i(\Gamma,\sigma):=\mathscr{D}_i(\Gamma,\sigma,\mathfrak{c}),\,\,\,\mathscr{F}_i(\Gamma):=\mathscr{F}_i(\Gamma,\mathfrak{c})$$
for any trivial coloring $\mathfrak{c}$ of $\Gamma$. These are the {\it $i^\textup{th}$ island boundary count} and the set of {\it island subgraphs on $i$ colors} respectively.
\begin{rem}
(1) The notation $f$ is used to remind us of the fact that the number of path components of $\Sigma- \sigma(\Gamma')$ is the number of faces, if we consider a triangulation of the surface by $\Gamma'$. For a planar connected graph $\Gamma$, i.e., an embedded graph $\sigma:\Gamma\hookrightarrow S^2$, the number of faces $f_\sigma(\Gamma)$ is the number of boundary components of a thickening of $\sigma(\Gamma)$ (figure \ref{fig:bdry-count}).  \\
(2) We are thinking of connected induced embedded subgraphs as \textit{islands} inside the surface. This is the reason behind our terminology. There is no real justification behind the nomenclature of $\beta_{(\Gamma,\sigma)}(-1)$ being the island count. It is a short term that connects with $\beta_{(\Gamma,\sigma)}(x)$ being called the island polynomial. We often drop $\sigma$ from the subscript and write $\beta_\Gamma(-1)$, when the embedding $\sigma$ is planar. \\
(3) The graph embeddings we are considering in Definition \ref{defn:beta-c} are not necessarily cellular. See figure \ref{fig:bouquet} for examples of non-cellular embeddings with island count $3,2$ and $1$ respectively.
\end{rem}
\begin{figure}[!h]
\centering
\includegraphics[trim={0 29 0 0}, clip, scale=0.7]{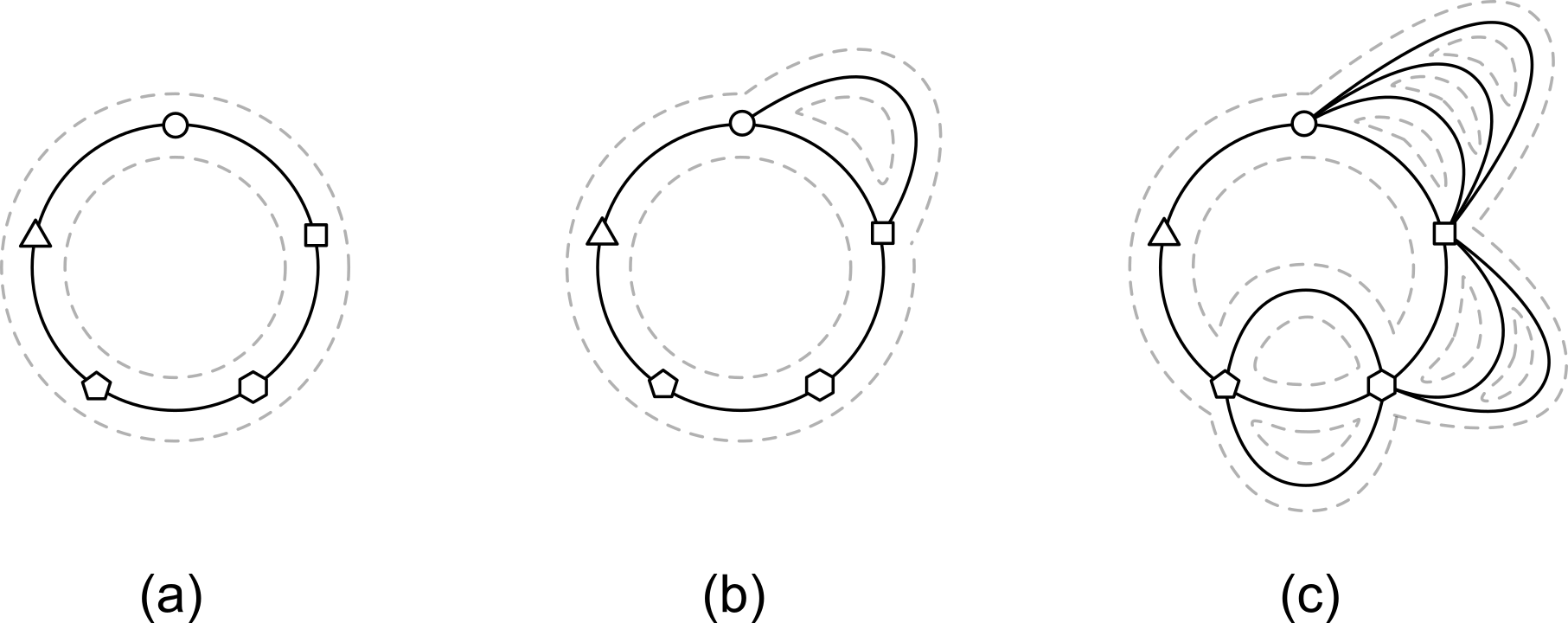}
\caption{For planar graphs with trivial coloring, the term $\mathcal{D}_{|V|}(\Gamma)$ counts the number of boundary components, indicated by dashed gray curves.}
\label{fig:bdry-count}
\end{figure}
\begin{eg}[Planar graphs]
For any connected graph $\Gamma$ which can be embedded in the plane (or, equivalently the $2$-sphere) via $\sigma:\Gamma\hookrightarrow \Sigma$, we have 
\begin{equation*}\label{eqn:face-planar}
f_\sigma(\Gamma)=2-v(\Gamma)+e(\Gamma).
\end{equation*}
This is due to the famous Euler's formula $v-e+f=2$ for connected planar graphs. Thus, $f_\sigma(\Gamma)$ is independent of the planar embedding, provided $\Gamma$ is connected.
\end{eg}
\begin{lemma}\label{planar-H0}
Let $\Gamma$ be a finite graph on $N$ vertices such that any proper induced subgraph is a disjoint union of trees. Then for any embedding $\sigma:\Gamma\hookrightarrow \Sigma$ inside a surface and $i<N$
\bgd
\mathscr{D}_i(\Gamma,\sigma)=\sum_{\Gamma'\in \mathscr{F}_i(\Gamma)} \textup{no. of connected components of $\Gamma'$}.
\edd
\end{lemma}
It can be verified that $\Gamma$, as above, is either the cycle graph $C_N$ or a disjoint union of trees. In this case, $\mathscr{D}_i(\Gamma,\sigma)$ counts the total number of connected induced subgraphs formed out of all possible induced subgraphs of $\Gamma$ with $i$ number of vertices. This count was first used to measure the multi-partite information of a topologically ordered system \cite{wen2013,Kitaev_2003,pollmann_2010} in theoretical condensed matter physics, where the physical attributes of a subgraph $\Gamma'\subseteq \Gamma$ were dependent on the number of connected induced subgraphs of $\Gamma'$. 
\begin{proof}
Let $\Gamma'$ be an induced subgraph of $\Gamma$; it will have components $\Gamma_1,\ldots,\Gamma_k$, all of which are trees. For an embedding into a surface, note that $f_\sigma(\Gamma')=\sum_j f_\sigma(\Gamma_j)=k$ as the complement of any embedded tree is connected. 
\end{proof}
\begin{defn}\label{defn:str-set}
    For a graph $\Gamma$, let 
$$\mathscr{S}_d(\Gamma):=\{\beta_{(\Gamma,\sigma,\mathfrak{c})}(x)\,|\,\sigma:\Gamma\hookrightarrow \Sigma, \mathfrak{c}:V(\Gamma)\twoheadrightarrow \mathcal{C}, |\mathcal{C}|=d+1\}$$
denote the colored island polynomials obtained from vertex-colorings of $\Gamma$ using $d+1$ colors, varying over embeddings. We define the {\it structure set of polynomials} of $\Gamma$ to be
$$\mathscr{S}(\Gamma):=\sqcup_{d\geq 0} \,\mathscr{S}_d(\Gamma).\label{eqn:str-set}$$
\end{defn}
Let $T_n$ denote any tree on $n$ vertices. It is evident that $\mathscr{S}(T_2)=\{1,x+2\}$. More interestingly, 
$$\mathscr{S}(T_3)=\{1,x+2,x+3,x^2+4x+3\}$$
$$\mathscr{S}(T_4)=\{1,x+2,x+3,x+4, x^2+4x+3, x^2+5x+4, x^3+6x^2+9x+4\}$$
implying that the structure set is independent of the combinatorial structure of the tree on $n\leq 4$ vertices. In fact, this is shown to be true (Corollary \ref{cor:str-tree}) for all trees. 


\subsection{Basic transformations: adding pendant edges and loops}\label{subsec:edge-loop}
\begin{figure}[!h]
\centering
\includegraphics[trim={0 40 0 0}, clip, scale=0.25]{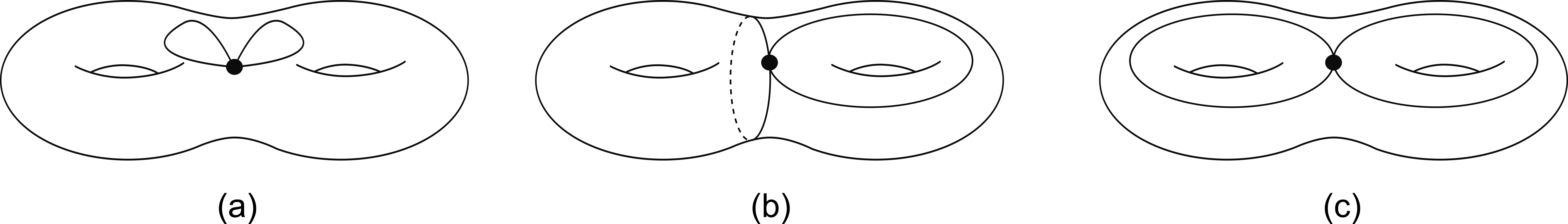}
\caption{Various non-isotopic embeddings of $B_2$ in a genus $2$ surface.}
\label{fig:bouquet}
\end{figure}

\hf The intricate nature of vertex-coloring does not have a major impact for trivial colorings. This lets us deduce nice structural properties of the island polynomial. We shall start with an analysis of two elementary families of graphs: (i) edgeless graphs and (ii) bouquet graphs. As any two embeddings of $E_n$, the edgeless graph on $n$ vertices, in a connected surface are isotopic, the $i^\textup{th}$ island boundary count is given by ${n \choose i}i$. It follows that for the trivial coloring, we obtain ${\beta}_{E_n}(x)=n(1+x)^{n-1}$. In particular, $\beta_{E_n}(0)=n$ captures the fact that the graph has $n$ components. \\
\hf Let $B_n$ denote the bouquet graph, i.e., a graph on one vertex and with $n$ loops, equipped with an embedding $\sigma:B_n\hookrightarrow \Sigma$. Note that the colored island polynomial is the island polynomial as the graph has only one vertex. If $\Sigma=S^2$ or $\Sigma=\mathbb{R}^2$, then $\Sigma - \sigma(B_n)$ has $n+1$ components (due to Euler's Theorem, which is related to the fact that any simple loop separates the plane). If $\Sigma$ is any higher genus surface, then $\Sigma - \sigma(B_n)$ has at most $n+1$ components. Thus, $\beta_{(B_n,\sigma)}(x)\leq n+1$, as loops may not always separate $\Sigma$ (refer to figure \ref{fig:bouquet}). 
\begin{prop}[Disjoint union]\label{prop:disjoint-union}
For vertex-colored embedded graphs $(\Gamma_i,\sigma_i,\mathfrak{c}_i)$, both embedded in $\Sigma$, $i=1,2$ with the images disjoint, i.e., $\sigma_1(\Gamma_1)\cap\sigma_2(\Gamma_2)=\varnothing$, and no color in common, i.e., $\mathfrak{c}_1(V(\Gamma_1))\cap \mathfrak{c}_2(V(\Gamma_2))=\varnothing$, we have
\begin{equation}\label{eqn:disjoint-color}
{\beta}_{(\Gamma_1\sqcup\Gamma_2, \sigma_1\sqcup\sigma_2,\mathfrak{c}_1\sqcup\mathfrak{c}_2)}(x)=(1+x)^{c_2}{\beta}_{(\Gamma_1,\sigma_1,\mathfrak{c}_1)}(x)+(1+x)^{c_1}{\beta}_{(\Gamma_2,\sigma_2, \mathfrak{c}_2)}(x),
\end{equation}
where $c_i$ is the number of colors used in the vertex-coloring of $\Gamma_i$. \\
\end{prop}
\begin{proof}
Let $\Gamma'=\Gamma_1'\sqcup \Gamma_2'\in \mathscr{F}_i(\Gamma_1\sqcup\Gamma_2,\mathfrak{c}_1\sqcup \mathfrak{c}_2)$. Note that $\Gamma_j'\in \mathscr{F}_{i_j}(\Gamma_j, \mathfrak{c}_j)$, for $j=1,2$, where $i_1+i_2=i$, and 
$$f_\sigma(\Gamma')=f_{\sigma_1}(\Gamma_1')+f_{\sigma_2}(\Gamma_2')$$
This implies that 
$$\mathscr{D}_i(\Gamma_1\sqcup \Gamma_2,\sigma_1\sqcup\sigma_2,\mathfrak{c}_1\sqcup\mathfrak{c}_2) = \sum_{\Gamma'\in \mathscr{F}_i(\Gamma_1\sqcup\Gamma_2,\mathfrak{c}_1\sqcup \mathfrak{c}_2)}f_{\sigma_1}(\Gamma_1')+f_{\sigma_2}(\Gamma_2')$$
where each $f_{\sigma_1}(\Gamma_1')$ and $f_{\sigma_2}(\Gamma_2')$ appear ${c_2\choose i_2}$ and ${c_1\choose i_1}$ times respectively in the right hand side. Thus, 
$$\mathscr{D}_i(\Gamma_1\sqcup \Gamma_2,\sigma_1\sqcup\sigma_2,\mathfrak{c}_1\sqcup\mathfrak{c}_2)= \sum_{i_1} {c_2\choose i-i_1}\mathscr{D}_{i_1}(\Gamma_1,\sigma_1,\mathfrak{c}_1)+\sum_{i_2} {c_1\choose i-i_2}\mathscr{D}_{i_2}(\Gamma_2,\sigma_2,\mathfrak{c}_2)$$
This implies the formula.
\end{proof}
\begin{rem}\label{rem:disjoint}
The same formula in Proposition \ref{prop:disjoint-union} holds if $\sigma_i:\Gamma_i\hookrightarrow \Sigma_i$ with $\Sigma_1$ and $\Sigma_2$ being disjoint surfaces, and the colorings are disjoint. 
\end{rem}
\begin{cor}\label{cor:disjoint-color}
Let $\Gamma=\Gamma_1\sqcup\cdots\sqcup \Gamma_k$ is a disjoint union of graphs, where $\Gamma_i$ need not be connected, and $\sigma:\Gamma\to\Sigma$ be an embedding. Let $\mathfrak{c}:V(\Gamma)\to\mathcal{C}$ be a vertex-coloring such that $\mathfrak{c}(V(\Gamma_i))\cap \mathfrak{c}(V(\Gamma_j))=\varnothing$ if $i\neq j$. Then, if $\mathfrak{c}_i=\mathfrak{c}|_{\Gamma_i}$ and $\sigma_i=\sigma|_{\Gamma_i}$, we have
\begin{equation}\label{eqn:disjoint-cor}
{\beta}_{(\Gamma,\sigma,\mathfrak{c})}(x)=(1+x)^{c-c_1}{\beta}_{(\Gamma_1,\sigma_1,\mathfrak{c}_1)}(x)+\cdots+(1+x)^{c-c_k}{\beta}_{(\Gamma_k,\sigma_k,\mathfrak{c}_k)}(x)
\end{equation}
where $c=|\mathfrak{c}(V(\Gamma))|$ and $c_i=|\mathfrak{c}(V(\Gamma_i))|$. In particular, if $\mathfrak{c}$ is a trivial coloring, then
\begin{equation}\label{eqn:disjoint-pol}
{\beta}_{(\Gamma, \sigma)}(x)=(1+x)^{n-n_1}{\beta}_{(\Gamma_1,\sigma_1)}(x)+(1+x)^{n-n_2}{\beta}_{(\Gamma_2,\sigma_2)}(x)+\cdots+(1+x)^{n-n_k}{\beta}_{(\Gamma_k,\sigma_k)}(x)
\end{equation}
where $n=|V(\Gamma)|$ and $n_i=|V(\Gamma_i)|$.
\end{cor}
We may deduce another useful result.
\begin{cor}\label{cor:vertex-deletion}
Let $(\Gamma,\sigma)$ be an embedded graph on $n$ vertices. If $v\in V(\Gamma)$, then we have the identity
\begin{equation*}\label{eqn:vertex-deletion}
\beta_{(\Gamma,\sigma)}(x)-(1+x)\beta_{(\Gamma-v,\sigma)}(x)=(1+x)^{n-n'}\big(\beta_{(\Gamma',\sigma)}(x)-(1+x)\beta_{(\Gamma'-v,\sigma)}(x)\big)
\end{equation*}
where $\Gamma'$ is the component of $\Gamma$ containing $v$ and has $n'$ vertices.
\end{cor}

\begin{figure}[h!]
\centering
\includegraphics[scale=0.8]{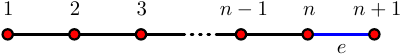}
\caption{$P_{n+1}$ is $P_n$ with a pendant edge}
\label{fig:appendix}
\end{figure}

\hf If $T$ is a tree on $n$ vertices with the trivial coloring, then it follows from Lemma \ref{planar-H0} that ${\beta}_{(T,\sigma)}(x)$ is independent of $\sigma$. Recall that a tree is built in iteration, starting from a vertex, by adding an edge to an existing vertex. We shall prove a general result regarding such graphs, called pendant graphs.
\begin{defn}[Adding a pendant edge]
Let $\Gamma$ be a graph. A graph $\Gamma_\textup{app}$, formed by adding an edge at an existing vertex of $\Gamma$, will be called {\it $\Gamma$ with a pendant edge}. 
\end{defn}
The new vertex in $\Gamma_\textup{app}$ is more commonly called a {\it pendant vertex}. In figure \ref{fig:appendix}, the vertices $1$ and $n+1$ are both pendant vertices as they have valency $1$.
Note that if $(\Gamma,\sigma)$ is an embedded graph, then there may be more than one extension of $\sigma$ to an embedded appended graph $(\Gamma_\textup{app},\tilde{\sigma})$. 
\begin{figure}[h!]
\centering
\includegraphics[scale=0.3]{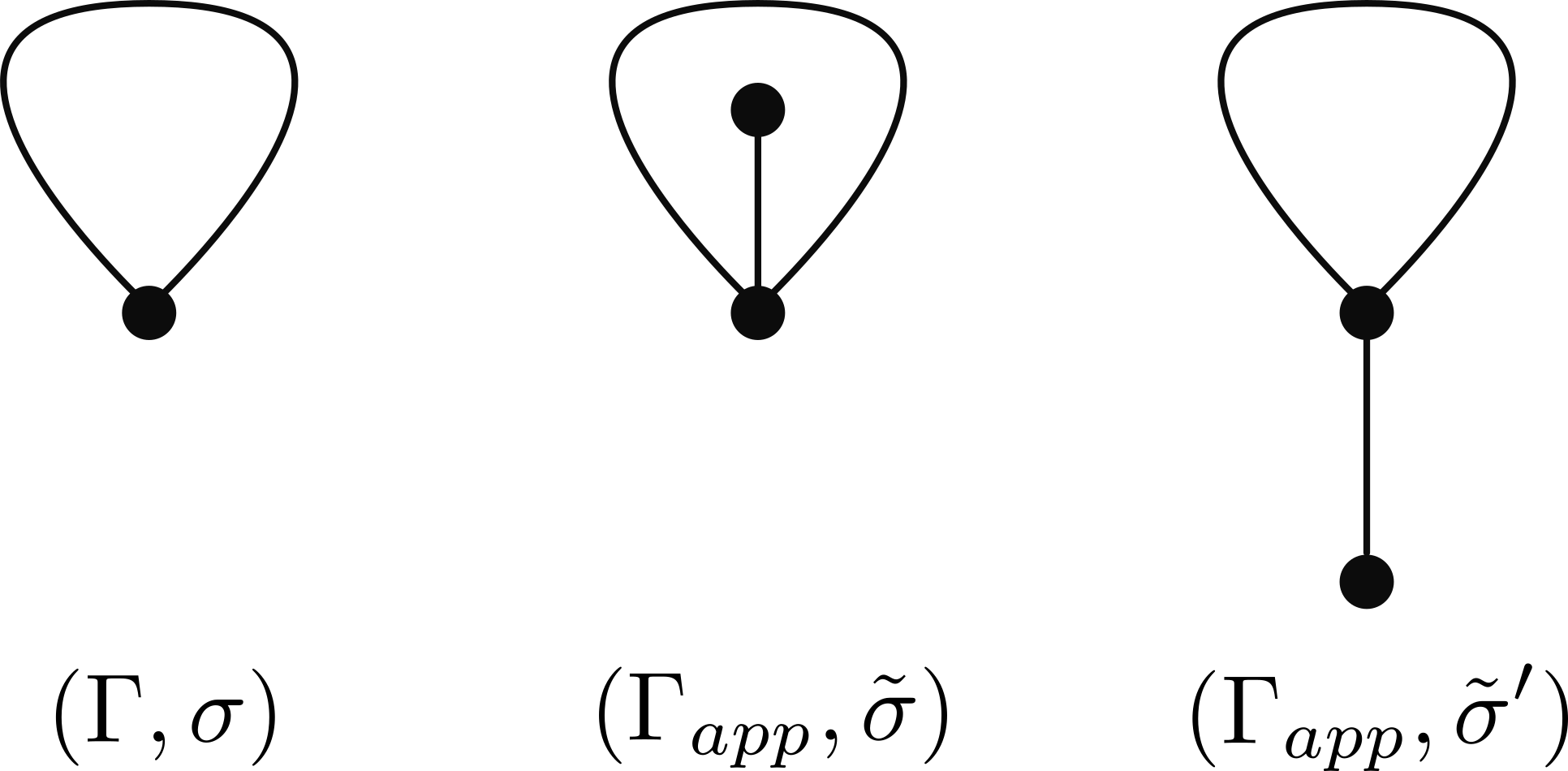}
\caption{Non-isotopic planar embeddings (graphs in (b) and (c) are isomorphic \& obtained from (a))}
\label{fig:pendant}
\end{figure}
\begin{prop}[Adding a pendant edge]\label{prop:Gamma-app}
Let $(\Gamma,\sigma,\mathfrak{c})$ and let $(\Gamma_\textup{app},\tilde{\sigma}, \tilde{\mathfrak{c}})$ be given such that the vertex-coloring of $\Gamma_\textup{app}$ is an extension of the vertex-coloring of $\Gamma$, i.e., $\tilde{\mathfrak{c}}|_{V(\Gamma)}=\mathfrak{c}$. Let $e_p$ denote the edge containing the pendant vertex $v_p$ and $v_0\in V(\Gamma)$. There are three possibilities for $\tilde{\mathfrak{c}}$:\\
\hf (i) $\tilde{\mathfrak{c}}(v_p)\not\in \mathfrak{c}(V(\Gamma))$;\\
\hf (ii) $\tilde{\mathfrak{c}}(v_p)=\mathfrak{c}(v_0)$;\\
\hf (iii) $\tilde{\mathfrak{c}}(v_p)\in \mathfrak{c}(V(\Gamma))$ but $\tilde{\mathfrak{c}}(v_p)\neq \mathfrak{c}(v_0)$.\\
Then, with $N=|\mathfrak{c}(V(\Gamma))|$, the following holds
\begin{eqnarray}
  {\beta}_{(\Gamma_\textup{app},\tilde{\sigma},\tilde{\mathfrak{c}})}(x) & = & (1+x)\big({\beta}_{(\Gamma,\sigma,\mathfrak{c})}(x)+(1+x)^{N-2}\big) \label{eqn:pen-new-color}\\
  {\beta}_{(\Gamma_\textup{app},\tilde{\sigma},\tilde{\mathfrak{c}})}(x) & = & {\beta}_{(\Gamma,\sigma,\mathfrak{c})}(x) \label{eqn:pen-same-color}\\
  {\beta}_{(\Gamma_\textup{app},\tilde{\sigma},\tilde{\mathfrak{c}})}(x) & = & {\beta}_{(\Gamma,\sigma,\mathfrak{c})}(x)+(1+x)^{N-2}. \label{eqn:pen-same-non-adj-color}
\end{eqnarray}
for the cases (i), (ii) and (iii) respectively.
\end{prop}
Note that case (iii) does not appear when $\Gamma$ has one vertex. We may now deduce the following results.
\begin{cor}\label{cor:tree-island}
The island polynomial for any tree $T_n$ on $n$ vertices is given by
\begin{equation}\label{eqn:tree-island}
\beta_{T_n}(x)=(1+x)^{n-1}+(n-1)(1+x)^{n-2}.
\end{equation}
Let $\Gamma=T_1\sqcup\cdots\sqcup T_r$ be a disjoint union of trees $T_i$'s on $n_i$ vertices respectively. Then 
\begin{equation}\label{eqn:r-comp-tree}
\beta_{T_1\sqcup\cdots\sqcup T_r}(x)=r(1+x)^{n-1}+(n-r)(1+x)^{n-2} 
\end{equation}
where $n=n_1+\cdots+n_r$.
\end{cor}
\begin{proof}
The first identity \eqref{eqn:tree-island} can be proven using induction and \eqref{eqn:pen-new-color}. The second identity \eqref{eqn:r-comp-tree} can be deduced from \eqref{eqn:tree-island} and \eqref{eqn:disjoint-pol}. We have omitted the embedding as for graphs which are disjoint union of trees, the island polynomial does not depend on the embedding.
\end{proof}
\begin{cor}\label{cor:it-pendant}
Let $\Gamma'$ be obtained from a graph $\Gamma$ by iteratively adding pendant edges $\{e_1,\ldots,e_k\}$. Let $\sigma:\Gamma'\hookrightarrow \Sigma$ be an embedding and let $\mathfrak{c}$ be a vertex-coloring of $\Gamma'$ such that the colors of the pendant vertices $v_j$ of $e_j$'s are distinct and $\mathfrak{c}(\Gamma)\cap \{\mathfrak{c}(v_1),\ldots,\mathfrak{c}(v_k)\}=\varnothing$. Then 
\begin{equation}\label{eqn:it-pendant}
    \beta_{(\Gamma',\sigma,\mathfrak{c})}(x)=(1+x)^k \beta_{(\Gamma,\sigma,\mathfrak{c})}(x)+k(1+x)^{N+k-2}
\end{equation}
where $\Gamma$ requires $N$ colors. 
\end{cor}
The proof involves iterating \eqref{eqn:pen-new-color}. 
\begin{proof}[Proof of Proposition \ref{prop:Gamma-app}]
We first derive the result for case (ii). Let $c_0=\tilde{\mathfrak{c}}(v_p)=\mathfrak{c}(v_0)$. Fix a set $S$ of $i$ colors $\mathfrak{c}(V(\Gamma))$. If $c_0\not\in S$, then there is the identity bijection between $\mathscr{F}_S(\Gamma_\textup{app},\tilde{\mathfrak{c}})$ and $\mathscr{F}_S(\Gamma,\mathfrak{c})$ The corresponding contributions to $\mathscr{D}_i$'s are the same. If $c_0\in S$, then there is a natural bijection between $\mathscr{F}_S(\Gamma,\mathfrak{c})$ and $\mathscr{F}_S(\Gamma_\textup{app},\tilde{\mathfrak{c}})$; it maps $\Gamma'$ to $\Gamma'\cup e_p$. As adding an edge (via a pendant vertex) to a graph does not change the connectivity of the complement, it follows that $\mathscr{D}_i$'s are the same. This implies \eqref{eqn:pen-same-color}.\\
\hf For case (i), note that for any pair of spaces $(Y,B)$, the connected components of $Y-B$ are in bijection with the connected components of $Y/B -[B]$, i.e., the number of connected components of $Y-B$ is the rank of $\tilde{H}_0(Y/B;\mathbb{Z})$. We also note that $(\Sigma,\tilde{\sigma}(\Gamma_\textup{app}))$ and $(\Sigma,\sigma(\Gamma))$ are good pairs. For a good pair $(X,A)$, there are isomorphisms $H_j(X,A;\mathbb{Z})\cong \tilde{H}_j(X/A;\mathbb{Z})$. As $\Gamma_\textup{app}$ deformation retracts to $\Gamma$, the long exact sequence (in homology) for the triple $\sigma(\Gamma)\subset \tilde{\sigma}(\Gamma_\textup{app})\subset \Sigma$, implies that the inclusion of pairs $\iota:(\Sigma,\sigma(\Gamma))\to (\Sigma,\tilde{\sigma}(\Gamma_\textup{app}))$ induce isomorphisms
$$\iota_\ast: H_j((\Sigma,\sigma(\Gamma);\mathbb{Z})\to H_j(\Sigma,\tilde{\sigma}(\Gamma_\textup{app});\mathbb{Z})$$
In particular, when $j=0$ we obtain 
\begin{equation}\label{eqn:Dapp}
\mathscr{D}_{N+1}(\Gamma_\textup{app},\tilde{\sigma},\tilde{\mathfrak{c}})=\mathscr{D}_N(\Gamma,\sigma,\mathfrak{c})
\end{equation}
The count $\mathscr{D}_S(\Gamma_\textup{app},\tilde{\sigma},\tilde{\mathfrak{c}})$ for $1\leq |S|=j\leq N$ comes from three mutually exclusive and exhaustive cases\\
\hf (a) $S\subset \mathfrak{c}(V(\Gamma))$, and\\
\hf (b) $c_0:=\tilde{\mathfrak{c}}(v_p)\in S$, $\mathfrak{c}(v_0)\not\in S$, and\\
\hf (c) $c_0\in S$, $\mathfrak{c}(v_0)\in S$.\\
The count for (a) is $\mathscr{D}_S(\Gamma,\sigma,\mathfrak{c})$. The count for (c) is $\mathscr{D}_{S-\{c_0\}}(\Gamma,\sigma,\mathfrak{c})$, where $\mathfrak{c}(v_0)\in S$. The count for (b) is $\mathscr{D}_{S-\{c_0\}}(\Gamma,\sigma,\mathfrak{c})+1$, where $\mathfrak{c}(v_0)\not\in S$ and the $+1$ is for the isolated vertex $v_p$ in $\mathscr{F}_S(\Gamma_\textup{app},\tilde{\sigma},\widetilde{\mathfrak{c}})$. Summing over all subsets of $\mathfrak{c}(V(\Gamma))\cup\{c_0\}$ of size $j$, we obtain
\begin{equation}\label{eqn:Djapp}
\mathscr{D}_j(\Gamma_\textup{app},\tilde{\sigma},\tilde{\mathfrak{c}})=\mathscr{D}_{j}(\Gamma,\sigma,\mathfrak{c})+\mathscr{D}_{j-1}(\Gamma,\sigma,\mathfrak{c})+{N-1\choose j-1}.
\end{equation}
Combining \eqref{eqn:Dapp}, \eqref{eqn:Djapp} with the definition of the total connected induced subgraph boundary polynomial ${\beta}$, we obtain \eqref{eqn:pen-new-color}
\begin{equation*}
{\beta}_{(\Gamma_\textup{app},\tilde{\sigma},\tilde{\mathfrak{c}})}(x)=(1+x){\beta}_{(\Gamma,\sigma,\mathfrak{c})}(x)+(1+x)^{N-1}. 
\end{equation*}
\hf For case (iii), \eqref{eqn:Dapp} still holds, i.e.,
\begin{equation*}
\mathscr{D}_{N}(\Gamma_\textup{app},\tilde{\sigma},\tilde{\mathfrak{c}})=\mathscr{D}_N(\Gamma,\sigma,\mathfrak{c})
\end{equation*}
The count $\mathscr{D}_S(\Gamma_\textup{app},\tilde{\sigma},\tilde{\mathfrak{c}})$ for $1\leq |S|=j\leq N-1$ comes from four mutually exclusive and exhaustive cases\\
\hf (a) $c_0:=\tilde{\mathfrak{c}}(v_p),\mathfrak{c}(v_0)\in S$, and\\
\hf (b) $c_0\in S$, $\mathfrak{c}(v_0)\not\in S$, and\\
\hf (c) $c_0\not\in S$, $\mathfrak{c}(v_0)\in S$, and\\
\hf (d) $c_0\not\in S$, $\mathfrak{c}(v_0)\not\in S$.\\
In all three cases (a), (c) and (d), the count is $\mathscr{D}_S(\Gamma,\sigma,\mathfrak{c})$ while for (b) it is $\mathscr{D}_S(\Gamma,\sigma,\mathfrak{c})+1$, where the $+1$ comes from the isolated vertex $v_p$ in $\mathscr{F}_S(\Gamma_\textup{app},\tilde{\sigma},\widetilde{\mathfrak{c}})$. As there are ${N-2 \choose j-1}$ such graphs satisfying (b), on summing over $S$ of size $j$, we obtain \eqref{eqn:pen-same-non-adj-color}
\begin{equation*}
{\beta}_{(\Gamma_\textup{app},\tilde{\sigma},\tilde{\mathfrak{c}})}(x)={\beta}_{(\Gamma,\sigma,\mathfrak{c})}(x)+(1+x)^{N-2}. 
\end{equation*}
\end{proof}
\hf Given an embedded graph $(\Gamma,\sigma)$, fix a vertex $v$ and add a loop $e$ at $v$ to create a graph $\tilde{\Gamma}$. An embedding of $\tilde{\Gamma}$, extending $\sigma$, is a choice of an embedded loop in $\big(\Sigma-\sigma(\Gamma)\big)\cup\{\sigma(v)\}$. 
Typically these are not isotopic choices (figure \ref{fig:self-loop-ext}). Let $\tilde{\sigma}$ be such an extension of $\sigma$. 
\begin{figure}[!ht]
\centering
\includegraphics[scale=0.38]{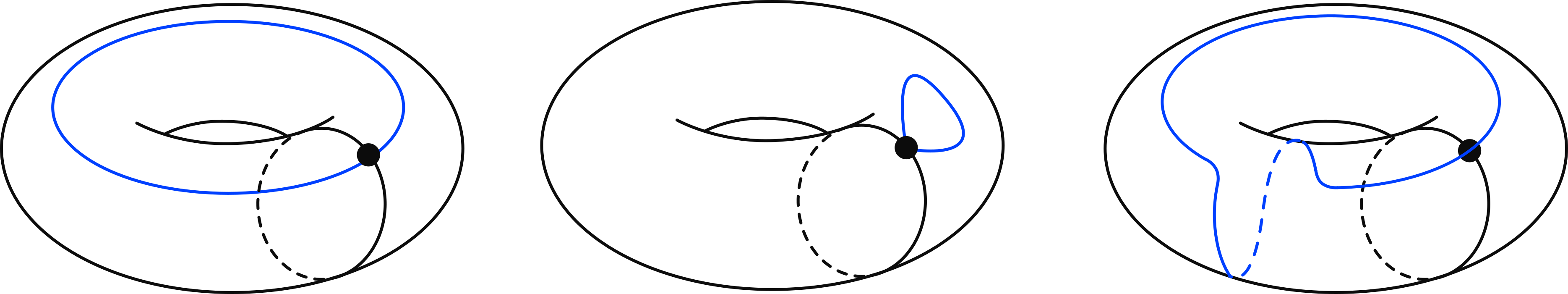}
\caption{Non-isotopic extensions of an embedded graph}
\label{fig:self-loop-ext}
\end{figure}

\begin{figure}[!ht]
\centering
\includegraphics[scale=0.32]{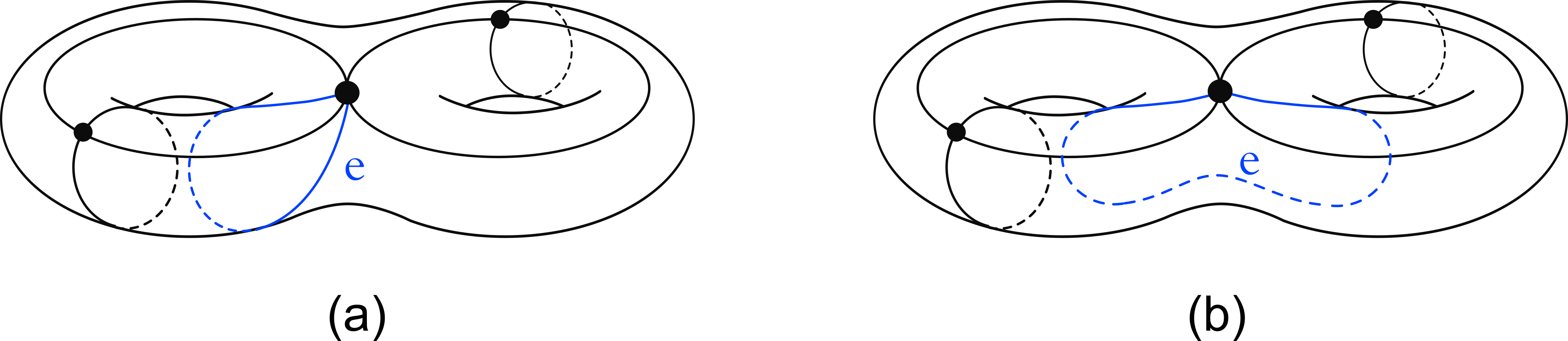}
\caption{Two homologically non-trivial loops that disconnect $\Sigma-\Gamma$}
\label{fig:genus_2}
\end{figure}

\begin{prop}[Adding a loop]\label{prop:self-loop}
Let $(\tilde{\Gamma},\tilde{\sigma},\mathfrak{c})$ be a vertex-colored, embedded graph on $n>1$ vertices, obtained from $(\Gamma,\sigma,\mathfrak{c})$ by adding a loop $e$ at a vertex $v_0$. If $\Sigma$ is oriented, then there are three mutually exclusive and exhaustive possibilities for $\tilde{\sigma}(e)$:\\
\hf (i) $\tilde{\sigma}(e)$ is homologically trivial, i.e., $[\tilde{\sigma}(e)]=0\in H_1(\Sigma;\mathbb{Z})$;\\
\hf (ii) $\tilde{\sigma}(e)$ does not disconnect $(\Sigma-\sigma(\Gamma))\cup\{\sigma(v_0)\}$; \\
\hf (iii) $\tilde{\sigma}(e)$ is homologically non-trivial and disconnects $(\Sigma-\sigma(\Gamma))\cup\{\sigma(v_0)\}$. \\
Then, with $N=|\mathfrak{c}(V(\Gamma))|$, the following holds
\begin{eqnarray}
  {\beta}_{(\tilde{\Gamma},\tilde{\sigma},\tilde{\mathfrak{c}})}(x) & = & {\beta}_{(\Gamma,\sigma,\mathfrak{c})}(x)+(1+x)^{N-1} \label{eqn:loop-trivial}\\
  {\beta}_{(\tilde{\Gamma},\tilde{\sigma},\tilde{\mathfrak{c}})}(x) & = & {\beta}_{(\Gamma,\sigma,\mathfrak{c})}(x) \label{eqn:loop-non-trivial}
\end{eqnarray}
for the cases (i) and (ii) respectively.
\end{prop}
For planar embeddings, (ii) and (iii) are vacuous conditions as every loop is homologically trivial and separates the plane. However, for other surfaces, this is subtle as there could be loops that are homologically non-trivial and disconnect $(\Sigma-\sigma(\Gamma))\cup\{v_0\}$. Figure \ref{fig:genus_2} exhibit two such loops with island polynomials $2x^2+5x+3$ and $2x^2+4x+3$ respectively. 
The analysis for (iii) in Proposition \ref{prop:self-loop} is intricate and there are no clean, nice formula as in (i) or (ii). For example, in figure \ref{fig:genus_2}, the two graphs are obtained from a common graph whose island polynomial is $x^2+4x+3$. The change in this polynomial when adding a loop (see figure \ref{fig:genus_2}) is by $x^2+x$ in (a) and by $x^2$ in (b). 
\begin{proof}
We prove the two cases, in order.\\
\textbf{Case i}: $[\gamma]=0\in H_1(\Sigma;\Z)$\\
Let $S\subseteq \mathfrak{c}(V(\Gamma))$ be if size $j$. If $\mathfrak{c}(v_0)\not\in S$, then $\mathscr{D}_S(\tilde{\Gamma},\tilde{\sigma}, \mathfrak{c})=\mathscr{D}_S(\Gamma,\sigma, \mathfrak{c})$. If $\mathfrak{c}(v_0)\in S$, then $\mathscr{D}_S(\tilde{\Gamma},\tilde{\sigma}, \mathfrak{c})=\mathscr{D}_S(\Gamma,\sigma, \mathfrak{c})+1$, where the $+1$ for the disconnection due to $\tilde{\sigma}(e)$. Therefore,
\bgd
\mathscr{D}_j(\tilde{\Gamma},\tilde{\sigma},\mathfrak{c})=\mathscr{D}_j(\Gamma,\sigma,\mathfrak{c})+\sum 1
\edd
where the sum is over all subgraphs $\mathscr{F}_S(\tilde{\Gamma},\tilde{\sigma},\mathfrak{c})$ where $|S|=j$ and $\mathfrak{c}(v_0)\in S$. As there are ${N-1\choose j-1}$ such subgraphs, we obtain \eqref{eqn:loop-trivial}.\\
\textbf{Case ii}: $[\gamma]$ does not disconnect $(\Sigma-\sigma(\Gamma))\cup\{v\}$\\
The condition means that $[\gamma]\neq 0\in H_1(\Sigma;\Z)$. In this case the counts $\mathscr{D}_j$ for $\Gamma$ and $\tilde{\Gamma}$ agree, implying an equality of polynomials ${\beta}_{\tilde{\Gamma},\tilde{\sigma},\mathfrak{c}}={\beta}_{\Gamma,\sigma,\mathfrak{c}}$.
\end{proof}

\subsection{Basic transformations: coloring over vertices}\label{subsec:c=c'}

\begin{figure}[!h]
\centering
\includegraphics[scale=0.62]{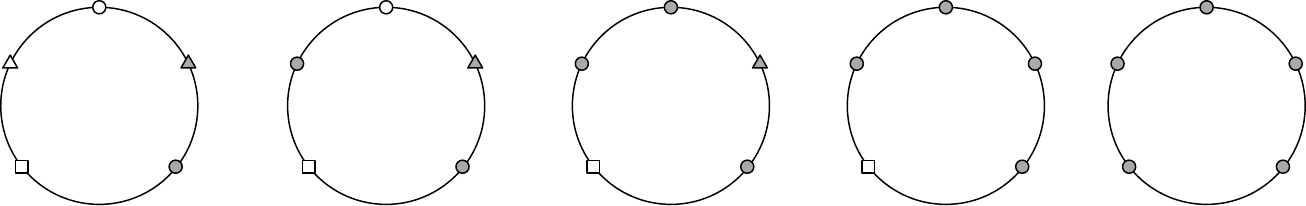}
\caption{Various vertex-colorings of $C_5$}
\label{fig:C5-col}
\end{figure}

\hf Let $(\Gamma,\sigma, \mathfrak{c})$ be an embedded colored graph. Let $c$ and $c'$ be distinct colors used in the vertex-coloring of $\Gamma$. We denote by $\Gamma_{c=c'}$ the underlying embedded graph of $\Gamma$ with the erstwhile $c$-colored vertices now colored by $c'$. Note that the colored island polynomials for $\Gamma_{c=c'}$ and $\Gamma_{c'=c}$ are identical. We will relate the polynomial of $\Gamma$ to that of $\Gamma_{c=c'}$. We shall denote by $\Gamma_{\neq c},\Gamma_{\neq c'}$ and $\Gamma_{\neq c,c'}$ the largest subgraphs of $\Gamma$ on vertices not colored by $c$, $c'$ and $\{c,c'\}$ respectively. These three graphs come naturally equipped with embeddings and vertex-colorings, denoted by $\sigma$ and $\mathfrak{c}$.
\begin{theorem}\label{thm:color}
Let $\mathfrak{c}$ be a vertex-coloring of $\Gamma$ which is not uni-colored. The colored island polynomial for $\Gamma$ and $\Gamma_{c=c'}$ are related by the identity
\begin{equation}\label{eqn:ver-col}
{\beta}_{(\Gamma,\sigma,\mathfrak{c})}(x)=x{\beta}_{(\Gamma_{c=c'},\sigma,\mathfrak{c})}(x)+{\beta}_{(\Gamma_{\neq c},\sigma,\mathfrak{c})}(x)+{\beta}_{(\Gamma_{\neq c'},\sigma,\mathfrak{c})}(x)-(1+x){\beta}_{(\Gamma_{\neq c,c'},\sigma,\mathfrak{c})}(x).
\end{equation}
\end{theorem}
\begin{proof}
The count $\mathscr{D}_j(\Gamma)$ is a sum of $\mathscr{D}_S(\Gamma)$'s where $S$ runs over $j$-element subsets of the color set. Any such set $S$ is one of the following four types: \\
\hf (i) $c,c'\in S$ (ii) $c\in S$, $c'\not\in S$ (iii) $c'\in S$, $c\not\in S$ (iv) $c\not\in S$, $c'\not\in S$. \\
The corresponding subgraph $\mathscr{F}_S(\Gamma,\sigma,\mathfrak{c})$ then is one of the following graphs:\\
\hf (i) $\mathscr{F}_{S-c}(\Gamma_{c=c'},\sigma,\mathfrak{c})$ and $c'\in S-c$; \\
\hf (ii) $\mathscr{F}_S(\Gamma_{\neq c'},\sigma,\mathfrak{c})$ and $c\in S$; \\
\hf (iii) $\mathscr{F}_S(\Gamma_{\neq c},\sigma,\mathfrak{c})$ and $c'\in S$; \\
\hf (iv) $\mathscr{F}_S(\Gamma_{\neq c,c'},\sigma,\mathfrak{c})$ (which implies $c\not\in S, c'\not\in S$). \\
Observe that $\mathscr{D}_j(\Gamma,\sigma,\mathfrak{c})$, denoted by $\mathscr{D}_j(\Gamma)$ in short, can be broken into a signed sum 
\bgd
\mathscr{D}_j(\Gamma)=\mathscr{D}_j^{c,c'}(\Gamma)+\mathscr{D}_j^{\neq c}(\Gamma)+\mathscr{D}_j^{\neq c'}(\Gamma)-\mathscr{D}_j^{\neq c,c'}(\Gamma)
\edd
of four terms:\\
\hf $\mathscr{D}_j^{c,c'}$ - contribution from colored island subgraphs that contain colors $c$ and $c'$;\\
\hf $\mathscr{D}_j^{\neq c}$ - contribution from colored island subgraphs that do not contain color $c$;\\
\hf $\mathscr{D}_j^{\neq c'}$ - contribution from colored island subgraphs that do not contain color $c'$;\\
\hf $\mathscr{D}_j^{\neq c,c'}$ - contribution from colored island subgraphs that do not contain colors $c$ or $c'$.\\
We have the following identities
\begin{eqnarray*}
\mathscr{D}_j^{c,c'}(\Gamma) & = & \mathscr{D}_{j-1}^{c'}(\Gamma_{c=c'})\\
\mathscr{D}_j^{\neq c}(\Gamma) & = & \mathscr{D}_{j}(\Gamma_{\neq c})\\
\mathscr{D}_j^{\neq c'}(\Gamma) & = & \mathscr{D}_{j}(\Gamma_{\neq c'})\\
\mathscr{D}_j^{\neq c,c'}(\Gamma) & = & \mathscr{D}_j(\Gamma_{\neq c,c'})\, =\, \mathscr{D}_{j}^{\neq c'}(\Gamma_{c=c'})
\end{eqnarray*}
where $\mathscr{D}_{k}^{c'}(\Gamma_{c=c'})$ and $\mathscr{D}_{j-1}^{\neq c'}(\Gamma_{c=c'})$ denote the $k^\textup{th}$ colored island boundary count of $\Gamma_{c=c'}$ using the color $c'$ and without using the color $c'$ respectively. The above identities will help us rewrite $\mathscr{D}_j(\Gamma)$ as 
\begin{eqnarray*}
\mathscr{D}_j(\Gamma) & = &   \mathscr{D}^{c'}_{j-1}(\Gamma_{c=c'})+\mathscr{D}_j(\Gamma_{\neq c})+\mathscr{D}_j(\Gamma_{\neq c'})-\mathscr{D}_j(\Gamma_{\neq c,c'})\\
& = & \mathscr{D}^{c'}_{j-1}(\Gamma_{c=c'})+\mathscr{D}^{\neq c'}_{j-1}(\Gamma_{c=c'})+\mathscr{D}_j(\Gamma_{\neq c})+\mathscr{D}_j(\Gamma_{\neq c'})-\big(\mathscr{D}_j(\Gamma_{\neq c,c'})+\mathscr{D}^{\neq c'}_{j-1}(\Gamma_{c=c'})\big)\\
& = & \mathscr{D}_{j-1}(\Gamma_{c=c'})+\mathscr{D}_j(\Gamma_{\neq c})+\mathscr{D}_j(\Gamma_{\neq c'})-\big(\mathscr{D}_j(\Gamma_{\neq c,c'})+\mathscr{D}_{j-1}(\Gamma_{\neq c,c'})\big)
\end{eqnarray*}
It follows that we have the identity \eqref{eqn:ver-col} of polynomials.
\end{proof}
\hf Start with a graph $(\Gamma,\sigma)$ with all vertices colored by distinct colors. Given a subset $S$ of $V(\Gamma)$, we may color all of $S$ by one color $c$ and obtain a colored graph $\Gamma_{S=c}$. The colored island polynomial of $\Gamma_{S=c}$ can be computed from that of $\Gamma$ by iteratively applying Theorem \ref{thm:color}, where we color all of $S$ by $c$, changing the color of one vertex to $c$ in $S$ at a time. A (vertex-colored) graph $\Gamma_{S=c}$ may be called a graph with \textit{multiple copies of the same node} (look at figure \ref{fig:C5-col}). 
\begin{rem}
We should clarify that two nodes of the same color need not have the same valency or the same number of loops. By ``multiple copies of the same node", we merely indicate the imagery that the same colored node is present in many positions. 
\end{rem}


\subsection{Basic transformations: contracting edges}
\label{subsec:collapse_edge}

\hf The graph obtained from $\Gamma$ by collapsing an edge $e$ is denoted by $\Gamma/e$ as the topology on it is the quotient topology obtained by identifying all points of $e$ to one point. This operation shall be called {\it edge contraction}. Given an embedding $\sigma:\Gamma\hookrightarrow \Sigma$, we want to provide an embedding of $\Gamma/e$ inside $\Sigma$. Recall that $\sigma(e)$ is homeomorphic to $[0,1]$. Thus, there exists a neighbourhood $U$ of $\sigma(e)$ in $\Sigma$ which is homeomorphic to a disk, 
\begin{figure}[h!]
\centering
\includegraphics[scale=0.55]{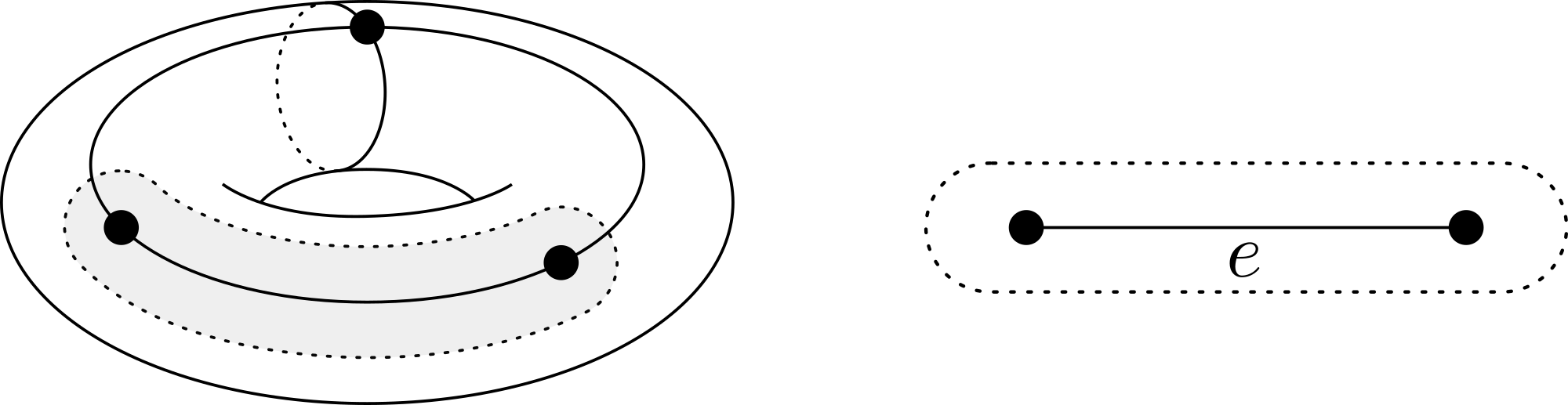}
\caption{A tubular neighbourhood of an edge}
\label{fig:band-aid}
\end{figure}
but thought of as a \textit{band-aid} (refer figure \ref{fig:band-aid}). Let $B_v$ and $B_w$ be open balls with center $v$ and $w$ respectively such that $B_v\cup B_w\subset U$.  
\begin{figure}[!ht]
\centering
\includegraphics[scale=0.7]{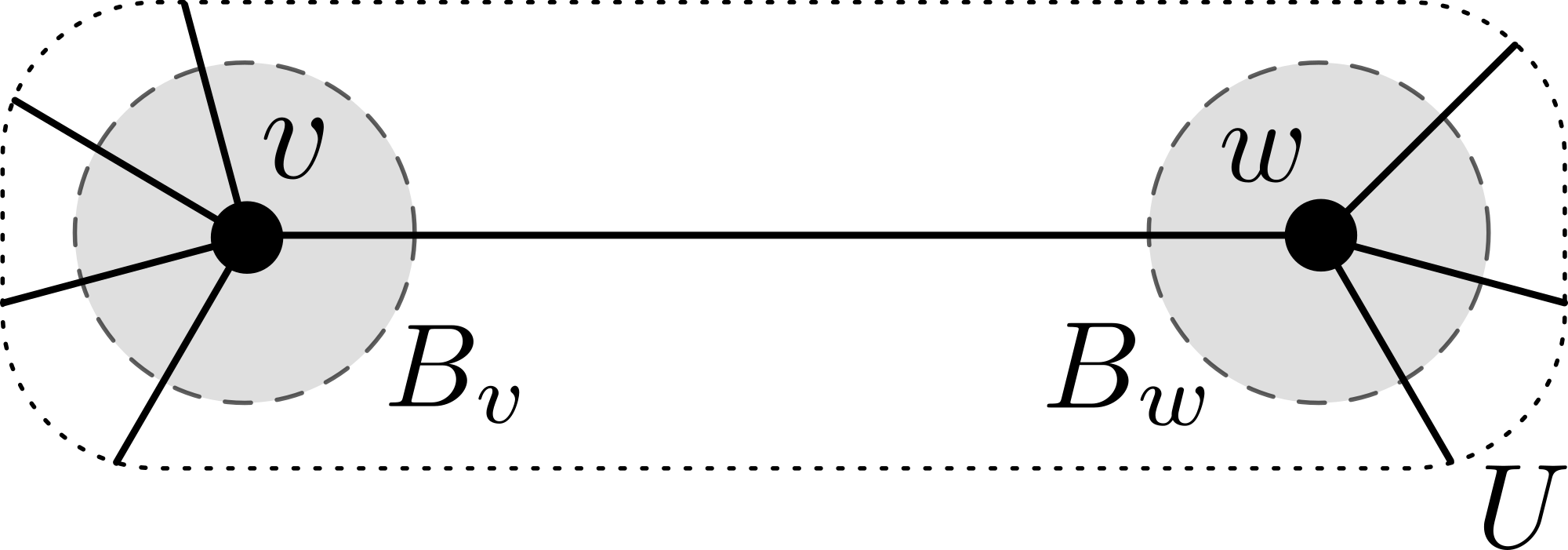}
\caption{A typical tubular neighbourhood of an edge}
\label{fig:band-aid2(a)}
\end{figure}
Note that there is an ordering of the edges emanating from each vertex.
\begin{defn}[Edge contraction]
Given an edge joining $v$ to $w$ in $(\Gamma,\sigma)$, we choose an open neighbourhood $U$ of $\sigma(e)$ as shown in figure \ref{fig:band-aid2(a)}. Introduce a new vertex $v_0$ at the middle of the unit interval joining $v$ to $w$. The edges terminating at $v$ are modified, preserving the ordering at $v$, as shown in 
\begin{figure}[!h]
\centering
\includegraphics[scale=0.7]{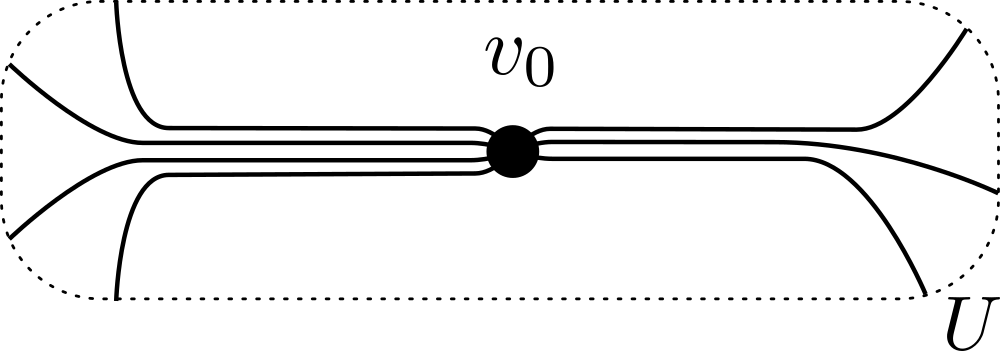}
\caption{Contracting an edge}
\label{fig:band-aid2(b)}
\end{figure}
figure \ref{fig:band-aid2(b)}, so as to join $v_0$. Similar modifications are done to the edges terminating at $w$. This new embedded graph will be denoted by $(\Gamma/e,\bar{\sigma})$.
\end{defn}
It can be shown that the isotopy class of the embedding $(\Gamma/e,\bar{\sigma})$ is independent of the choices made in the definition above. Moreover, the modification done to obtain $\Gamma/e$ is localized around $\sigma(e)$ and the new embedding $\bar{\sigma}$ agrees with $\sigma$ on the complement of $U$.\\
\hf Let $\mathfrak{c}$ be a vertex-coloring of $\Gamma$. Let $e$ be an edge between $v$ and $w$ with $v\neq w$. If $\mathfrak{c}(v)=\mathfrak{c}(w)$, then there is a natural coloring on $\Gamma/e$, where the vertex $v_0=[v]=[w]$ is colored by $\mathfrak{c}(v)$; we shall denote this vertex-colored graph by $(\Gamma/e,\overline{\sigma},\mathfrak{c})$. If $\mathfrak{c}(v)\neq \mathfrak{c}(w)$, then we need the following convention/definition.
\begin{defn}[Induced colorings on $\Gamma/e$]
Given an edge joining $v$ to $w$ in $(\Gamma,\sigma,\mathfrak{c})$, we consider the graph $(\Gamma/e, \overline{\sigma})$. Let $\mathfrak{c}(v)\neq \mathfrak{c}(w)$ be the colors of $v$ and $w$ respectively. Let $\mathfrak{c}_{v=w}$ denote the graph $(\Gamma,\sigma)$ where both vertices $v$ and $w$ are colored by $\mathfrak{c}(w)$. Let $\mathfrak{c}_{w=v}$ denote the graph $(\Gamma,\sigma)$ where both vertices $v$ and $w$ are colored by $\mathfrak{c}(v)$. \\
\hf We may color $\Gamma/e$ by leaving unchanged the colors of vertices other than $[v]=[w]$ and coloring the vertex $[v]$ with either $\mathfrak{c}(v)$ or $\mathfrak{c}(w)$. The first graph will be denoted by $(\Gamma/e,\overline{\sigma}, \mathfrak{c}_{w=v})$ and the latter will be denoted by $(\Gamma/e,\overline{\sigma}, \mathfrak{c}_{v=w})$. 
\end{defn}
Note that the number of colors may decrease by $1$ in the process of edge contraction (cf. figure \ref{fig:contract-edge}).
\begin{figure}[!h]
\centering
\includegraphics[scale=0.2]{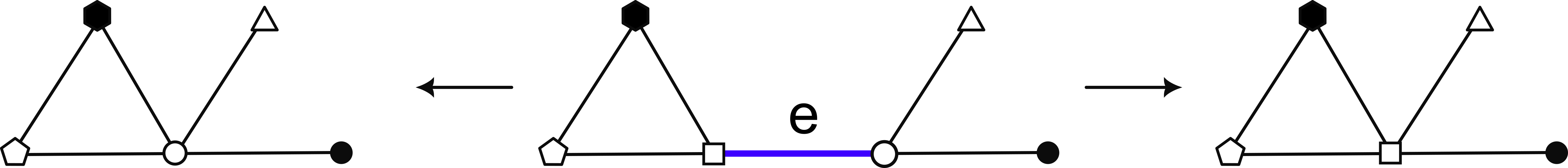}
\caption{Possible edge-contractions of a colored graph}
\label{fig:contract-edge}
\end{figure}
\begin{lemma}\label{lmm:collapse}
Let $e$ be an edge joining $v$ to $w$ in $(\Gamma, \sigma, \mathfrak{c})$. The colored island polynomials satisfy
\begin{eqnarray}
    \beta_{(\Gamma/e,\overline{\sigma}, \mathfrak{c}_{v=w})}(x) & = & \beta_{(\Gamma,\sigma,\mathfrak{c}_{v=w})}(x)\label{eqn:beta-v=w}\\
    \beta_{(\Gamma/e,\overline{\sigma}, \mathfrak{c}_{w=v})}(x) & = & \beta_{(\Gamma,\sigma,\mathfrak{c}_{w=v})}(x)\label{eqn:beta-w=v}
\end{eqnarray}
In particular, if $\mathfrak{c}(v)=\mathfrak{c}(w)$, then $\beta_{(\Gamma/e,\overline{\sigma}, \mathfrak{c})}(x)=\beta_{(\Gamma,\sigma,\mathfrak{c})}$.
\end{lemma}
\begin{proof}
We shall prove the first identity \eqref{eqn:beta-v=w}; the second one \eqref{eqn:beta-w=v} is similar. There are two scenarios: (a) the color $\mathfrak{c}(v)$ occurs in the vertex-coloring of both graphs, or (b) the color $\mathfrak{c}(v)$ does not occur in the vertex-coloring of both graphs. Note that (a) holds if and only if there exists $v'\not\in V(\Gamma)$ such that $\mathfrak{c}(v')=\mathfrak{c}(v)$ and (b) holds if $v$ is the only vertex of $\Gamma$ of its color. In (b), we require one less color than that of $\Gamma$ while in (a) we the color set in unchanged. \\
\hf We compare the counts $\mathscr{D}_S(\Gamma/e,\overline{\sigma}, \mathfrak{c}_{v=w})$ and $\mathscr{D}_S(\Gamma,\sigma,\mathfrak{c}_{v=w})$ for any subset $S$ of the vertex color set in both scenarios together. \\
Case (i): If $\mathfrak{c}(w)\in S$, $\mathfrak{c}(v)\not\in S$, then $\mathscr{F}_S(\Gamma,\sigma,\mathfrak{c}_{v=w})$ is obtained from $\mathscr{F}_S(\Gamma/e,\overline{\sigma}, \mathfrak{c}_{v=w})$ by adding a pendant edge to the vertex $w$. This does not change the colored island count.\\
Case (ii): If $\mathfrak{c}(w)\in S$, $\mathfrak{c}(v)\in S$, then $\mathscr{F}_S(\Gamma/e,\overline{\sigma}, \mathfrak{c}_{v=w})$ is obtained from $\mathscr{F}_S(\Gamma,\sigma,\mathfrak{c}_{v=w})$ by collapsing $e$. This does not change the colored island count. \\
Case (iii): If $\mathfrak{c}(w)\not\in S$, then both counts equal $\mathscr{D}_S(\Gamma-\{v,w\},\sigma, \mathfrak{c})$.\\
Thus, by summing over subsets we obtain \eqref{eqn:beta-v=w}.
\end{proof}
\begin{cor}\label{cor:collapse-tree}
Let $T$ be a tree in $(\Gamma,\sigma,\mathfrak{c})$ such that $T$ is uni-colored. Then $(\Gamma/T,\overline{\sigma}, \mathfrak{c})$ may be defined iteratively by collapsing one pendant edge of $T$ at a time since a tree is built out of a vertex by adding pendant edges. We have the following equality
\begin{equation}\label{eqn:collapse-tree}
\beta_{(\Gamma/T,\overline{\sigma}, \mathfrak{c})}(x)=\beta_{(\Gamma,\sigma,\mathfrak{c})}(x).
\end{equation}
\end{cor}
\begin{figure}[!h]
\centering
\includegraphics[scale=0.2]{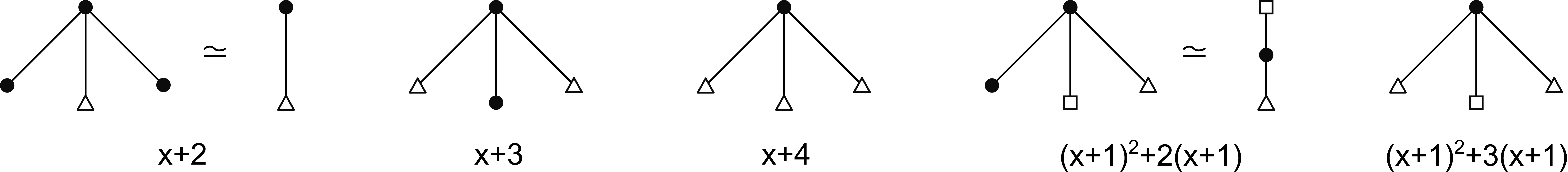}
\caption{All colored island polynomials (with at $2$ or $3$ colors) for a tree on $4$ vertices}
\label{fig:contract-tree}
\end{figure}
\begin{figure}[!h]
\centering
\includegraphics[scale=0.2]{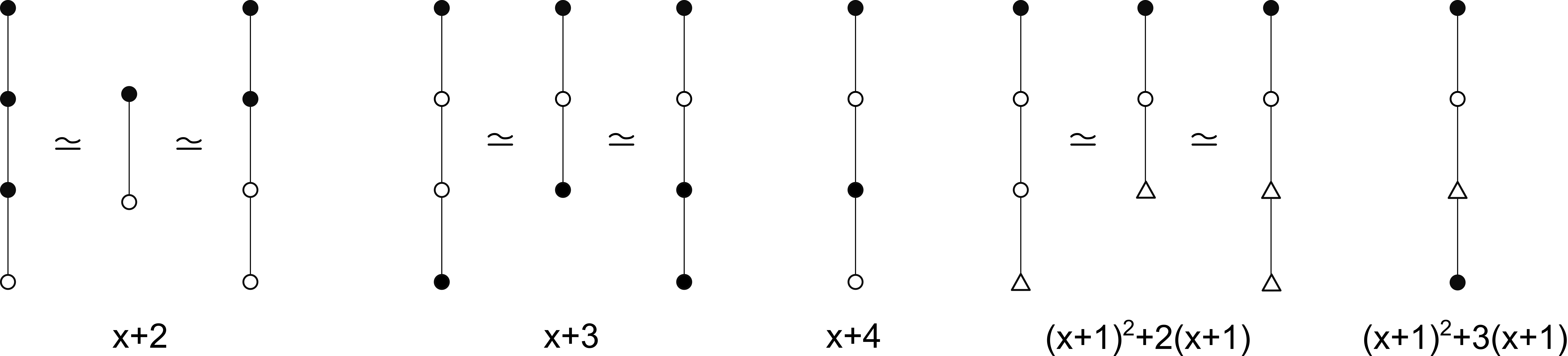}
\caption{All colored island polynomials (with at $2$ or $3$ colors) for $P_4$}
\label{fig:contract-path4}
\end{figure}
Aided with Theorem \ref{thm:color} and Lemma \ref{lmm:collapse}, we deduce the following result.
\begin{theorem}\label{thm:collapse}
Let $(\Gamma,\sigma,\mathfrak{c})$ be trivially colored. Then the polynomial for $\Gamma/e$ is given by
\begin{equation*}\label{eqn:collapse}
{\beta}_{(\Gamma,\sigma)}(x)=x{\beta}_{(\Gamma/e,\overline{\sigma})}(x)+{\beta}_{(\Gamma-v,\sigma)}(x)+{\beta}_{(\Gamma-w,\sigma)}(x)-(1+x){\beta}_{(\Gamma-\{v,w\},\sigma)}(x),
\end{equation*}
where $e$ joins $v$ and $w$.
\end{theorem}
There is no simpler formula for ${\beta}_{(\Gamma/e,\overline{\sigma})}$. Any reasonable formula would involve intrinsic features of the embedded graph $\Gamma$, like how many $k$ cycles are there involving either both $v$ and $w$ or neither. The formula for colored island polynomial is harder due to the fact that in Theorem \ref{thm:color}, we deal with graphs where all vertices with color $c'$ are now colored by $c$ or where we delete all vertices of color $c$. This is a drastic change in $\Gamma$ as opposed to deleting a vertex or coloring two adjacent vertices by the same color, as it is in Theorem \ref{thm:collapse}. \\
\hf We analyze the third coloring option for $(\Gamma/e,\overline{\sigma})$, where the vertex $v_0=[v]=[w]$ is colored by $c_0\not\in \mathfrak{c}(V)$; we denote this coloring by $\hat{\mathfrak{c}}$. We may assume that there are vertices in $\Gamma$ other than $v$ (respectively $w$) colored by $\mathfrak{c}(v)$ (respectively $\mathfrak{c}(w)$). We also assume that $\mathfrak{c}$ is not uni-colored. Applying Theorem \ref{thm:color} to $(\Gamma/e,\overline{\sigma},\hat{\mathfrak{c}})$ we get
$${\beta}_{(\Gamma/e,\overline{\sigma},\hat{\mathfrak{c}})}(x)=x{\beta}_{((\Gamma/e)_{c_0=\mathfrak{c}(v)},\overline{\sigma},\hat{\mathfrak{c}})}(x)+{\beta}_{((\Gamma/e)_{\neq c_0},\overline{\sigma},\hat{\mathfrak{c}})}(x)+{\beta}_{((\Gamma/e)_{\neq \mathfrak{c}(v)},\overline{\sigma},\hat{\mathfrak{c}})}(x)-(1+x){\beta}_{((\Gamma/e)_{\neq c_0,\mathfrak{c}(v)},\overline{\sigma},\hat{\mathfrak{c}})}(x).$$
Combining $((\Gamma/e)_{c_0=\mathfrak{c}(v)},\overline{\sigma},\hat{\mathfrak{c}})=(\Gamma/e, \overline{\sigma},\mathfrak{c}_{w=v})$ with \eqref{eqn:beta-w=v} and the isomorphism $((\Gamma/e)_{\neq c_0},\overline{\sigma},\hat{\mathfrak{c}})\cong (\Gamma-\{v,w\},\sigma,\mathfrak{c})$ we may write
$${\beta}_{(\Gamma/e,\overline{\sigma},\hat{\mathfrak{c}})}(x)=x{\beta}_{(\Gamma, \sigma,\mathfrak{c}_{w=v})}(x)+{\beta}_{(\Gamma-\{v,w\},\sigma,\mathfrak{c})}(x)+{\beta}_{((\Gamma/e)_{\neq \mathfrak{c}(v)},\overline{\sigma},\hat{\mathfrak{c}})}(x)-(1+x){\beta}_{((\Gamma/e)_{\neq c_0,\mathfrak{c}(v)},\overline{\sigma},\hat{\mathfrak{c}})}(x).$$


\subsection{Basic transformations: bridge deletion and wedge sum}
\label{subsec:delete_edge}

\begin{figure}[!ht]
\centering
\includegraphics[scale=0.45]{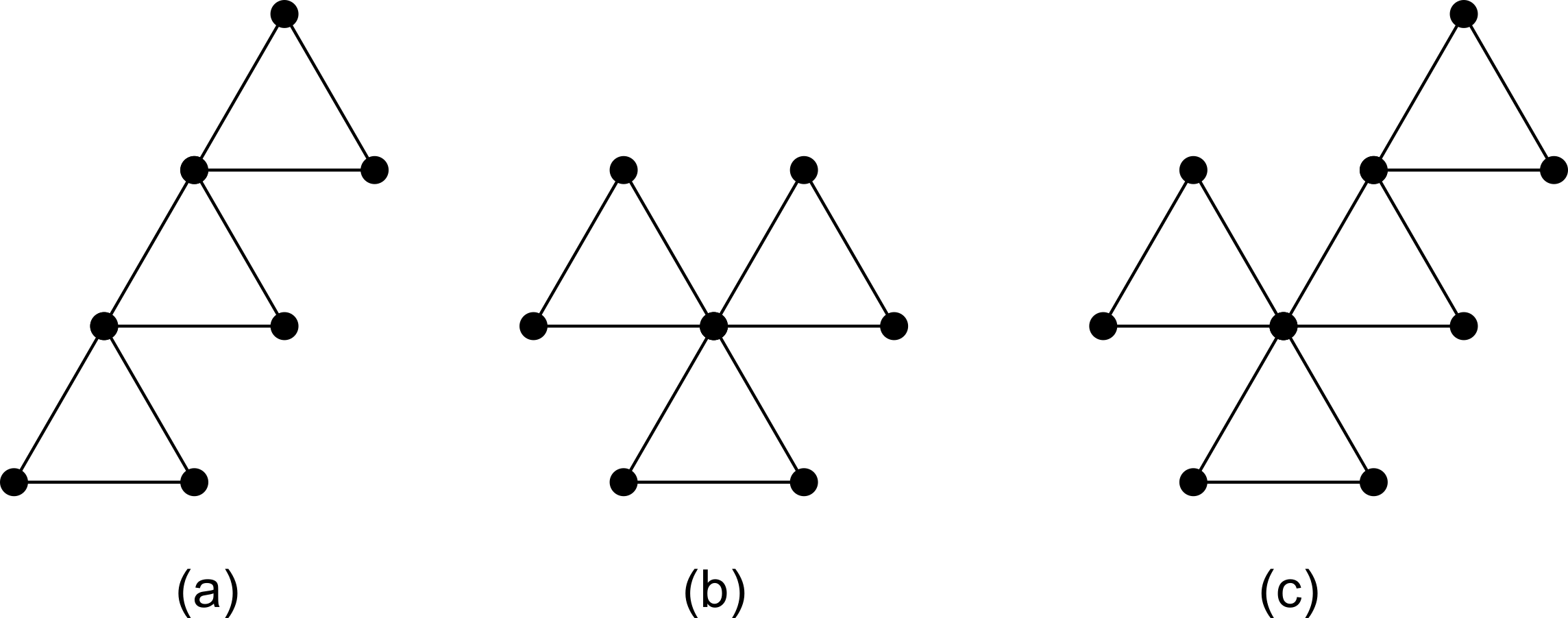}
\caption{(a) An iterated wedge sum of three $C_3$'s (b) A wedge sum of three $C_3$'s (c) An iterated wedge sum of four $C_3$'s}
\label{fig:wedge}
\end{figure}

\hf Recall that a bridge (or cut-edge) in a graph $\Gamma$ is an edge $e$ such that its removal increases the connectivity of the graph (by one). Thus, if $\Gamma$ is connected and $e$ is a bridge joining $v$ and $w$, then $\Gamma-e=\Gamma_1\sqcup\Gamma_2$, where $\Gamma_i$'s are connected subgraphs. We can recover $\Gamma$ from $\Gamma_1\sqcup\Gamma_2$ by adding an edge joining $v$ and $w$. This motivates the following definition.
\begin{defn}[Bridge graph]\label{defn:bridge}
Given two graphs $\Gamma_1$ and $\Gamma_2$, we may create a new graph $\Gamma_1\!\!\multimapdotboth\!\Gamma_2$ by introducing a new edge $e$ joining a pair of chosen vertices $v_i\in V(\Gamma_i)$, $i=1,2$. We shall call this edge the \textit{bridge} between $\Gamma_1$ and $\Gamma_2$. \\
\hf For embedded colored graphs $(\Gamma_i,\sigma_i,\mathfrak{c}_i), i=1,2$ in $\Sigma$, we assume that the embeddings are disjoint and a path $\gamma_e$ in $\Sigma$, corresponding to $e$, is chosen such that $\textup{image}(\gamma_e)\cap \sigma_i(\Gamma_i)=\sigma_i(v_i)$. The graph with vertex set $V(\Gamma_1)\sqcup V(\Gamma_2)$, edge set $E(\Gamma_1)\sqcup E(\Gamma_2)\sqcup\{e\}$, coloring $\mathfrak{c}_1\sqcup \mathfrak{c}_2$ and embedding $\sigma_\gamma:=\sigma_1\cup \gamma_e\cup \sigma_2$ will be called a {\it bridge graph} and denoted by $\Gamma_1\!\!\multimapdotboth\!\Gamma_2$.
\end{defn}
Note that the definition is dependent on the choice of $\gamma_e$ as different choices typically lead to non-isotopic embeddings. However, $\beta_{(\Gamma, \sigma_\gamma,\mathfrak{c}_1\sqcup\mathfrak{c}_2)}(x)$ is independent of this choice.
\begin{prop}\label{prop:bridge-del}
    Let $\Gamma:=\Gamma_1\!\!\multimapdotboth\!\Gamma_2$ be as defined in Definition \ref{defn:bridge}. Let $S=\{\mathfrak{c}_1(v_1),\mathfrak{c}_2(v_2)\}$ and $e$ be the bridge joining $v_1$ and $v_2$. Then, with $N$ being the number of colors used, we have 
    \begin{equation}\label{eqn:bridge-del}
    \beta_{(\Gamma_1\sqcup \Gamma_2,\sigma_1\sqcup\sigma_2, \mathfrak{c}_1\sqcup\mathfrak{c}_2)}(x)=\left\{\begin{array}{rl}
    \beta_{(\Gamma, \sigma_\gamma,\mathfrak{c}_1\sqcup\mathfrak{c}_2)}(x) + (1+x)^{N-1} & \textup{if $|S|=1$}\\
    \beta_{(\Gamma, \sigma_\gamma,\mathfrak{c}_1\sqcup\mathfrak{c}_2)}(x) + x(1+x)^{N-2} & \textup{if $|S|=2$}
    \end{array}\right.
    \end{equation}
\end{prop}
\begin{proof}
We compute $f_{\sigma_\gamma}(\mathscr{F}_T(\Gamma,\mathfrak{c}_1\sqcup\mathfrak{c}_2))$ in three mutually exclusive and exhaustive cases:\\
\hf (i) $S\cap T=\varnothing$: As $\mathscr{F}_T(\Gamma,\mathfrak{c}_1\sqcup\mathfrak{c}_2)=\mathscr{F}_T(\Gamma_1\sqcup\Gamma_2,\mathfrak{c}_1\sqcup\mathfrak{c}_2)  $, we get 
$$f_{\sigma_\gamma}(\mathscr{F}_T(\Gamma,\mathfrak{c}_1\sqcup\mathfrak{c}_2))=f_{\sigma}(\mathscr{F}_T(\Gamma_1\sqcup\Gamma_2,\mathfrak{c}_1\sqcup\mathfrak{c}_2))$$
\hf (ii) $S\subseteq T$: In this case, $S\cap T =S$ and we get
$$f_{\sigma_\gamma}(\mathscr{F}_T(\Gamma,\mathfrak{c}_1\sqcup\mathfrak{c}_2))+1=f_{\sigma}(\mathscr{F}_T(\Gamma_1\sqcup\Gamma_2,\mathfrak{c}_1\sqcup\mathfrak{c}_2))$$        
\hf (iii) $S\cap T\neq\varnothing$, $S\cap T\subset S$: In this case, $|S|=2$ and $|S\cap T|=1$. Thus $\mathscr{F}_T(\Gamma,\mathfrak{c}_1\sqcup\mathfrak{c}_2)=\mathscr{F}_T(\Gamma_1\sqcup\Gamma_2,\mathfrak{c}_1\sqcup\mathfrak{c}_2)  $, we get 
$$f_{\sigma_\gamma}(\mathscr{F}_T(\Gamma,\mathfrak{c}_1\sqcup\mathfrak{c}_2))=f_{\sigma}(\mathscr{F}_T(\Gamma_1\sqcup\Gamma_2,\mathfrak{c}_1\sqcup\mathfrak{c}_2))$$
A summing and counting argument (as in Proposition \ref{prop:self-loop}, for instance) now gives the result.
\end{proof}
\hf If $e$, joining $v$ and $w$ in a connected graph $\Gamma$, is a bridge, then $\Gamma/e$ is the wedge sum of the two components of $\Gamma-e$, obtained by identifying $v$ and $w$ together. 
\begin{defn}[Wedge sum]
A pointed graph is a graph $\Gamma$ with a choice of a vertex. Given a family of pointed graphs $\{\Gamma_\alpha\}_{\alpha\in J}$, each with a chosen vertex $v_\alpha\in V(\Gamma_\alpha)$, we define the {\it wedge sum} to be graph $\vee_\alpha \Gamma_\alpha$ with vertex set obtained from $\sqcup_\alpha  V(\Gamma_\alpha)$ by identifying all the chosen vertices together and edges determined naturally from the graphs. This graph can be given the quotient topology if $\Gamma_\alpha$'s are equipped with a topology. 
\end{defn}
Quite often, graphs constructed out of iterated (binary) wedge sums are more commonplace than a single wedge sum of graphs (refer to figure \ref{fig:wedge}).
We may think of $\vee_\alpha \Gamma_\alpha$ as the graph obtained from creating bridges between all possible $\Gamma_\alpha$ and $\Gamma_\beta$ with $\alpha\neq \beta$ and then contracting these bridges altogether. For the wedge sum of two graphs, we use the notation $\Gamma_1\vee \Gamma_2$, analogous to what is standard in the topological category. There is an isomorphism 
\bgd
(\Gamma_1\!\!\multimapdotboth\!\Gamma_2)/e\cong \Gamma_1\vee \Gamma_2
\edd
Combining this isomorphism with Proposition \ref{prop:disjoint-union} applied to a trivially colored $\Gamma_1\sqcup \Gamma_2$, Theorem \ref{thm:collapse} and Proposition \ref{prop:bridge-del} applied to $\Gamma_1\!\multimapdotboth\Gamma_2$ imply the following equality of island polynomials
\begin{equation}\label{eqn:wedge}
\beta_{(\Gamma_1\vee \Gamma_2,\sigma_1\vee\sigma_2)}(x)=(1+x)^{n_2-1}\beta_{(\Gamma_1,\sigma_1)}(x)+(1+x)^{n_1-1}\beta_{(\Gamma_2,\sigma_2)}(x)-(1+x)^{n_1+n_2-2}
\end{equation}
where $\Gamma_i$ has $n_i$ vertices. This generalizes the case of graphs with an appendix (refer to Proposition \ref{prop:Gamma-app}, \eqref{eqn:pen-new-color}), for the island polynomial. An iterated application of this gives the formula for wedge sums.
\begin{cor}\label{cor:wedge}
Let $(\Gamma_1\vee\cdots\vee \Gamma_k,\sigma\vee\cdots\vee\sigma_k)$ be an embedded trivially colored wedge sum graph, where $\Gamma_j$ has $n_j$ vertices. With $N=n_1+\cdots+n_k$ we have
\bgd
\beta_{(\Gamma_1\vee\cdots\vee \Gamma_k,\sigma\vee\cdots\vee\sigma_k)}(x)=(1+x)^{N-n_i-(k-1)}\beta_{(\Gamma_1,\sigma_1)}(x)+\cdots+(1+x)^{N-n_i-(k-1)}\beta_{(\Gamma_1,\sigma_1)}(x)-(k-1)(1+x)^{N-k}
\edd
\end{cor}

\subsection{Basic Transformations: subdividing edges}
\label{subsec:split_edge}

\hf Given a graph $\Gamma$, the subdivision of an edge $e$, joining $v$ and $w$, refers to the resulting graph $\Gamma_\textup{sp}$ obtained by inserting a vertex $v_0$ in the interior of $e$. For this process (as opposed to contracting an edge), we allow edges to be loops, i.e., $v=w$ is allowed. If $\sigma:\Gamma\hookrightarrow \Sigma$ is an embedding, then this extends to an embedding $\tilde{\sigma}:\Gamma_\textup{sp}\hookrightarrow \Sigma$. Note that the isotopy class of $(\tilde{\sigma},\Gamma_\textup{sp})$ is independent of the position of the new vertex. The graph $\Gamma_\textup{sp}$ is called an {\it edge subdivision} of $\Gamma$.\\
\begin{figure}[!h]
\centering
\includegraphics[scale=0.50]{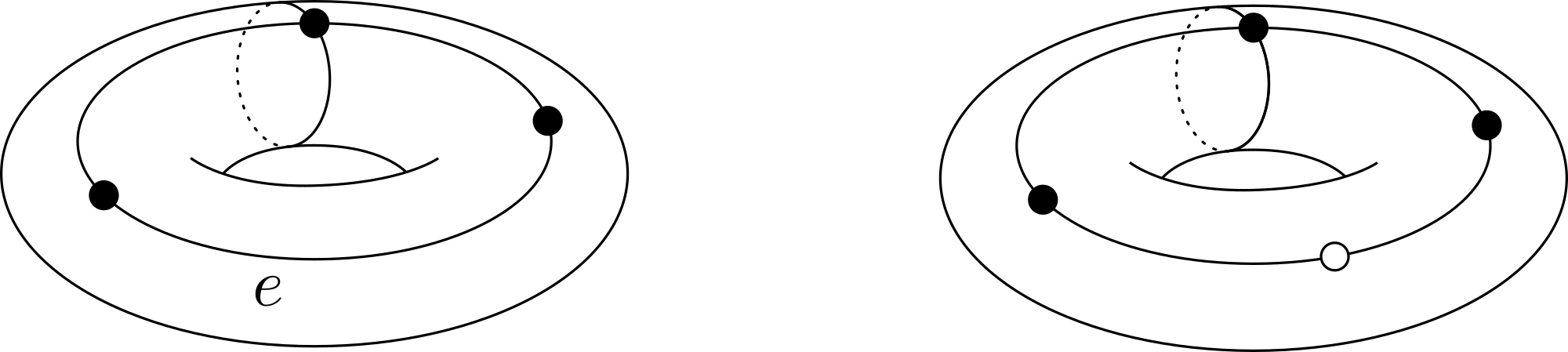}
\caption{Subdividing an edge}
\label{fig:split}
\end{figure}
\hf Given a vertex-coloring $\mathfrak{c}:V(\Gamma)\to\mathcal{C}$, there are two possibilities for a vertex-coloring of $\Gamma_\textup{sp}$:\\
\hf (i) color $v_0$ using $S=\{\mathfrak{c}(v),\mathfrak{c}(w)\}$;\\
\hf (ii) color $v_0$ using a new color $c_0$.\\
Note that in case (i), either $S$ is of size $2$ which means $e$ is necessarily not a loop and $v_0$ may get assigned one of two possible colors, or $S$ is of size $1$ and $v_0$ gets assigned that unique color. We shall denote the coloring in case (i) by $\mathfrak{c}_\textup{ext}$. In case (ii), the coloring is denoted by $\hat{\mathfrak{c}}$. 
\begin{prop}\label{prop:subdivide}
Let $(\Gamma,\sigma,\mathfrak{c})$ be a graph, vertex-colored by $N$ colors. Let $(\Gamma_\textup{sp},\sigma)$ be the graph obtained from $\Gamma$ by subdividing an edge $e$, joining $v$ and $w$. With $S=\{\mathfrak{c}(v),\mathfrak{c}(w)\}$, we have the identities
\begin{eqnarray}
{\beta}_{(\Gamma_\textup{sp},\sigma,\mathfrak{c}_\textup{ext})}(x) & = & \beta_{(\Gamma,\sigma,\mathfrak{c})}(x) \label{eqn:split1}\\
{\beta}_{(\Gamma_\textup{sp},\sigma,\hat{\mathfrak{c}})}(x) & = & x{\beta}_{(\Gamma,\sigma,\mathfrak{c})}(x)+{\beta}_{(\Gamma-e,\sigma,\mathfrak{c})}(x)+(1+x)^{N-|S|}\label{eqn:split2}
\end{eqnarray}
for case (i) and case(ii) respectively.
\end{prop}
\begin{proof}
The first case is clear as $\mathscr{F}_T(\Gamma_\textup{sp},\sigma,\mathfrak{c}_\textup{ext})=\mathscr{F}_T(\Gamma,\sigma,\mathfrak{c})$ if $S\cap T=\varnothing$ and $S=\{\mathfrak{c}(v),\mathfrak{c}(w)\}$. If $C\cap T\neq \varnothing$, then $\mathscr{F}_T(\Gamma_\textup{sp},\sigma,\mathfrak{c}_\textup{ext})$ is either isomorphic to $\mathscr{F}_T(\Gamma,\sigma,\mathfrak{c})$ (when $S\subseteq T$) or $\mathscr{F}_T(\Gamma_\textup{sp},\sigma,\mathfrak{c}_\textup{ext})$ is obtained from $\mathscr{F}_T(\Gamma,\sigma,\mathfrak{c})$ by adding a pendant edge (when $|S\cap T|=1$ and $|S|=2$). In all cases, $\mathscr{D}_T(\Gamma_\textup{sp},\sigma,\mathfrak{c}_\textup{ext})=\mathscr{D}_T(\Gamma,\sigma,\mathfrak{c})$, establishing \eqref{eqn:split1}. Note that if the new vertex is colored by $\mathfrak{c}(v_0)$ and $e_1$ is the new edge joining the new vertex and $v_0$, then $\Gamma_\textup{sp}/e_1$ is isomorphic to $\Gamma$. An application of Corollary \ref{cor:collapse-tree} (refer to \eqref{eqn:collapse-tree}) then implies \eqref{eqn:split1}.\\
\hf We need to analyze $\mathscr{D}_T(\Gamma_\textup{sp},\sigma,\hat{\mathfrak{c}})$ for the following cases:\\
\hf (a) $c_0\not\in T$: the face count equals $\mathscr{D}_T(\Gamma-e,\sigma,\mathfrak{c})$;\\
\hf (b) $c_0\in T$, $S\cap T=\varnothing$: as $\mathscr{F}_T(\Gamma_\textup{sp},\sigma,\hat{\mathfrak{c}})=\mathscr{F}_{T-c_0}(\Gamma,\sigma,\mathfrak{c})\sqcup \{v_0\}$, the face count equals $\mathscr{D}_{T-c_0}(\Gamma,\sigma,\mathfrak{c})+1$;\\
\hf (c) $c_0\in T$, $S\cap T\neq \varnothing$: if $S\subseteq T$, then $\mathscr{F}_T(\Gamma_\textup{sp},\sigma,\hat{\mathfrak{c}})\cong \mathscr{F}_{T-c_0}(\Gamma,\sigma,\mathfrak{c})$ as the underlying spaces are isomorphic. If $S\not\subseteq T$, then it follows that $|S\cap T|=1$ and $|S|=2$. In this case, $\mathscr{F}_T(\Gamma_\textup{sp},\sigma,\hat{\mathfrak{c}})$ is $\mathscr{F}_{T-c_0}(\Gamma,\sigma,\mathfrak{c})$ with a pendant edge (the appropriate half of $e$). Thus, the face count equals $\mathscr{D}_{T-c_0}(\Gamma,\sigma,\mathfrak{c})$ in both subcases.\\
\hf We now sum up the above counts to obtain
\begin{align*}
\mathscr{D}_i(\Gamma_\textup{sp},\sigma,\hat{\mathfrak{c}}) & = \sum_{|T|=i, c_0\not\in T} \mathscr{D}_T(\Gamma_\textup{sp},\sigma,\hat{\mathfrak{c}})+\sum_{|T|=i,c_0\in T} \mathscr{D}_T(\Gamma_\textup{sp},\sigma,\hat{\mathfrak{c}})\\
& = \sum_{|T|=i, c_0\not\in T} \!\!\!\!\mathscr{D}_T(\Gamma - e, \sigma, \mathfrak{c}) + \sum_{|T|=i,c_0\in T, S\cap T=\varnothing} \!\!\!\!\!\!\!\!\!\big(\mathscr{D}_{T-c_0}(\Gamma,\sigma,\mathfrak{c})+1\big)+\sum_{|T|=i,c_0\in T, S\cap T\neq \varnothing} \!\!\!\!\!\!\!\!\!\mathscr{D}_{T-c_0}(\Gamma,\sigma,\mathfrak{c})\\
& = \mathscr{D}_i(\Gamma - e, \sigma, \mathfrak{c})+{N-|S|\choose i-1}+\sum_{|T'|=i-1} \!\mathscr{D}_{T'}(\Gamma,\sigma,\mathfrak{c})
\end{align*}
Equation \eqref{eqn:split2} now follows by multiplying by $x^{i-1}$ and summing over $i$.      
\end{proof}


\subsection{Other Transformations: adding parallel edges}
\label{subsec:arrange_nearestneighbor}


\hf Given a graph $(\Gamma,\sigma,\mathfrak{c})$ in $\Sigma$, we consider adjacent vertices $v_1$ and $v_2$. We assume that $v_1\neq v_2$ as the case of adding a loop was discussed in \S \ref{subsec:edge-loop}.  
\begin{defn}\label{defn:sim_adj}
Given $(\Gamma,\sigma,\mathfrak{c})$ in $\Sigma$, creating a graph $(\tilde{\Gamma},\tilde{\sigma},\mathfrak{c})$ by adding a parallel edge between adjacent vertices $v_1$ and $v_2$ is choosing an extension of the embedding $\sigma$ to the graph $\tilde{\Gamma}$ obtained by adjoining a parallel edge $e'$, joining $v_1$ to $v_2$, to $\Gamma$. 
\end{defn}
\begin{figure}[!ht]
\centering
\includegraphics[scale=0.52]{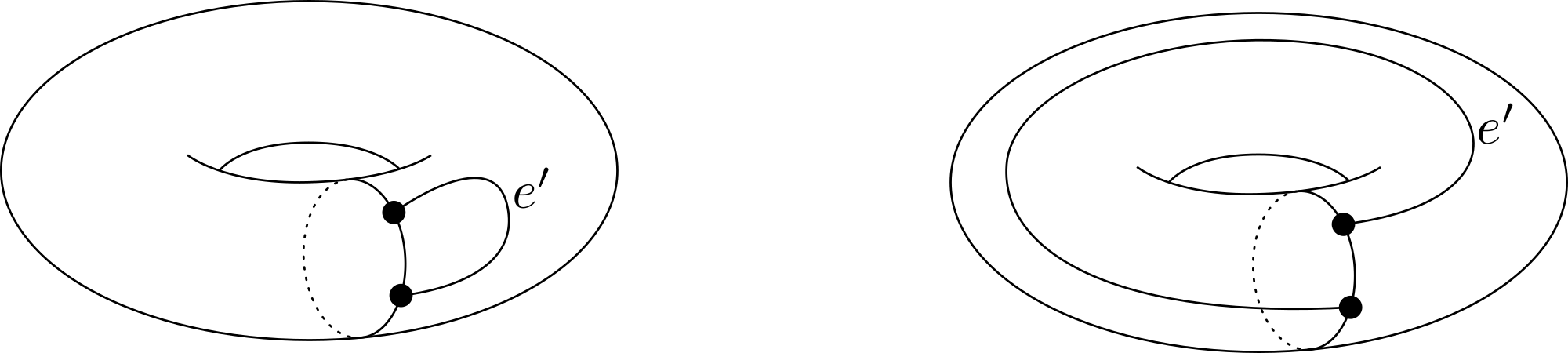}
\caption{Two non-isotopic extensions of $C_2$}
\label{fig:nn-adj}
\end{figure}
\begin{prop}\label{prop:nn-adj}
Let $(\tilde{\Gamma},\tilde{\sigma},\mathfrak{c})$ be a graph on $n>2$ vertices, obtained from $(\Gamma,\sigma,\mathfrak{c})$ by adding a parallel edge $e'$, joining vertices $v_1$ and $v_2$. Let $N$ be the number of colours and $S=\{\mathfrak{c}(v_1),\mathfrak{c}(v_2)\}$. \\
(i) If $f_{\tilde{\sigma}}(\mathscr{F}_S(\tilde{\Gamma},\mathfrak{c}))=f_{\sigma}(\mathscr{F}_S(\Gamma,\mathfrak{c}))+1$ then
\begin{equation}\label{eqn:parallel-edge}
\beta_{(\tilde{\Gamma},\tilde{\sigma},\mathfrak{c})}(x)=\beta_{(\Gamma,\sigma,\mathfrak{c})}(x)+x^{|S|-1}(1+x)^{N-|S|}.
\end{equation}
(ii) If $\tilde{\sigma}(e')$ does not disconnect $(\Sigma-\sigma(\Gamma))\cup\{v_1,v_2\}$, then ${\beta}_{(\tilde{\Gamma},\tilde{\sigma},\mathfrak{c})}(x)={\beta}_{(\Gamma,\sigma,\mathfrak{c})}(x)$.
\end{prop}
\begin{proof}
We argue case by case.\\
\textbf{Case i}: $f_{\tilde{\sigma}}(\mathscr{F}_S(\tilde{\Gamma},\mathfrak{c}))=f_{\sigma}(\mathscr{F}_S(\Gamma,\mathfrak{c}))+1$\\
If $T\subset \mathfrak{c}(V(\Gamma))$, then either $S\subseteq T$ or $T$ does not contain $S$. In the second case, $\mathscr{D}_T(\tilde{\Gamma},\tilde{\sigma},\mathfrak{c})=\mathscr{D}_T(\Gamma,\sigma,\mathfrak{c})$. In the first case, the hypothesis (of case i) implies that 
$\mathscr{D}_T(\tilde{\Gamma},\tilde{\sigma},\mathfrak{c})=\mathscr{D}_T(\Gamma,\sigma,\mathfrak{c})+1$. Thus,
$$\mathscr{D}_j(\tilde{\Gamma},\tilde{\sigma},\mathfrak{c})=\sum_{|T|=j}\mathscr{D}_T(\tilde{\Gamma},\tilde{\sigma},\mathfrak{c})=\sum_{|T|=j, S\not\subseteq T}\mathscr{D}_T(\Gamma,\sigma,\mathfrak{c}) + \sum_{|T|=j, S\subseteq T}\big(\mathscr{D}_T(\Gamma,\sigma,\mathfrak{c})+1\big)  $$
where the right hand side simplifes to $\mathscr{D}_j(\Gamma,\sigma,\mathfrak{c})+{N-|S|\choose j-|S|}$ as we add $1$ for all subgraphs on $j$ colours containing $S$. We obtain \eqref{eqn:parallel-edge}
\begin{equation*}
\beta_{(\tilde{\Gamma},\tilde{\sigma},\mathfrak{c})}(x)=\beta_{(\Gamma,\sigma,\mathfrak{c})}(x)+x^{|S|-1}(1+x)^{N-|S|}.
\end{equation*}
\textbf{Case ii}: $\tilde{\sigma}(e')$ does not disconnect $(\Sigma-\sigma(\Gamma))\cup\{v_1,v_2\}$\\
The condition means that no new cycles are formed in $\tilde{\Gamma}$ that increase the face count. In this case the counts $\mathcal{D}_j$ for $\Gamma$ and $\tilde{\Gamma}$ agree, implying an equality of polynomials ${\beta}_{(\tilde{\Gamma},\tilde{\sigma},\mathfrak{c})}={\beta}_{(\Gamma,\sigma,\mathfrak{c})}$.
\end{proof}
We have considered the two cases in Proposition \ref{fig:nn-adj} as these represent the two extreme situations. Adding a parallel edge in am embedded graph may not create new regions (case (ii)) or create a new cycle. Choose a cycle (containing $e'$) with the minimum number of vertices. Look at the set $S'$ of colors of the vertices of this cycle. As $S'\supseteq S$, case (i) is asserting that a cycle is created on $|S|$ colors.
Note that the island polynomial typically changes after adding a parallel edge.
\begin{figure}[!ht]
\centering
\includegraphics[scale=0.62]{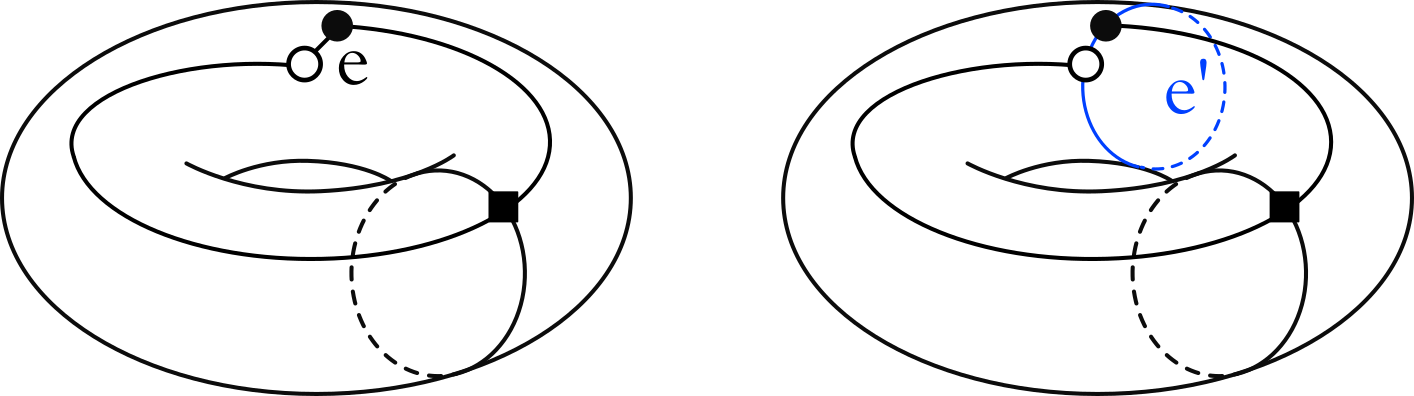}
\caption{Adding a parallel edge in the torus}
\label{fig:sim-adj}
\end{figure}
In figure \ref{fig:sim-adj}, the island polynomials for the two graphs are $3+3x+x^2$ and $3+3x+2x^2$ respectively.


\subsection{Other Transformations: adding shortcuts}\label{subsec:arrange_shortcircuit}

\begin{figure}[!h]
\centering
\includegraphics[scale=0.25]{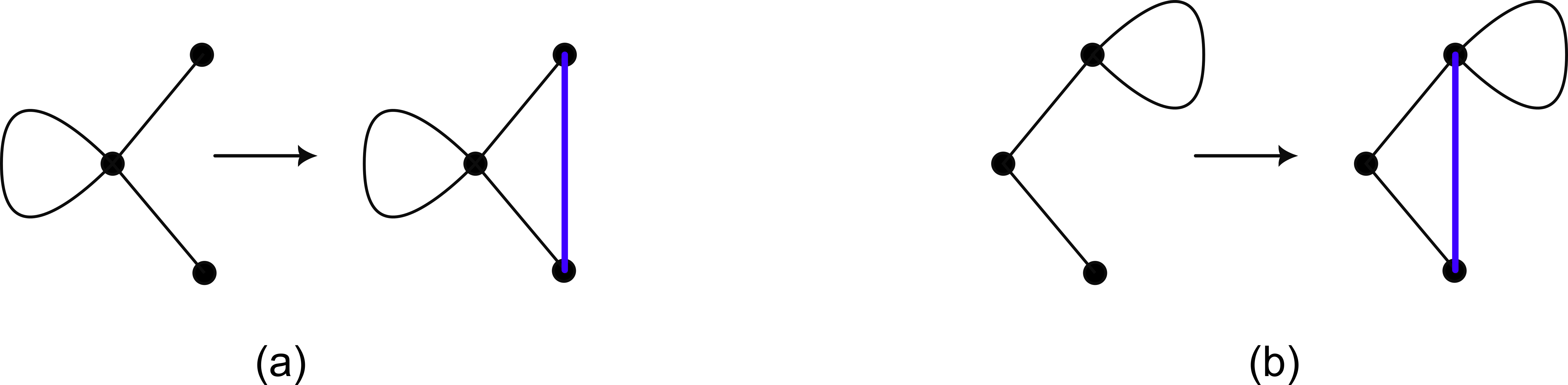}
\caption{(a) a shortcut which is not clean; (b) a clean shortcut}
\label{fig:shortcut}
\end{figure}

\begin{figure}[!ht]
\centering
\includegraphics[scale=0.55]{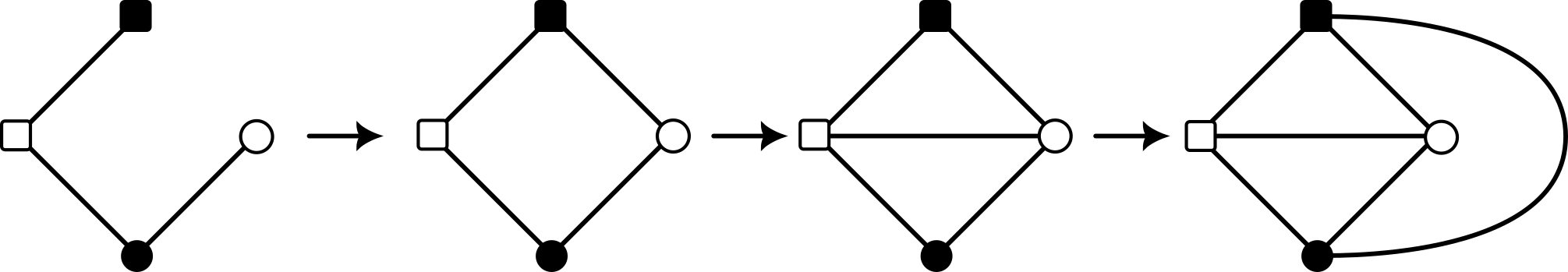}
\caption{Iterated shortcuts creating a planar $K_4$}
\label{fig:sc-K4}
\end{figure}

\hf Given a connected graph $\Gamma$, let $v$ and $w$ be vertices which are non-adjacent. We may add a new edge $e$ joining $v$ and $w$. Then $\Gamma\cup e$ will be termed as a graph obtained from $\Gamma$ by a {\it shortcut}, with $e$ serving as the shortcut. This decreases the distance between the vertices $v$ and $w$ to $1$. Creating a shortcut has interesting impact on embedded graphs and (colored) island polynomial. Although we did not encounter this operation in the literature, this is a very natural operation, when interpreted in the context of topological entanglement entropy in physics \cite{SiddharthaTEE2021}.  
\begin{defn}[Clean shortcut]\label{defn:csc}
Let $\Gamma\cup e$ be obtained from a connected graph $\Gamma$ via a shortcut $e$. We will call $e$, joining vertices $v\neq w$, a {\it clean shortcut} if there exists a path $P_k, k\geq 3$ in $\Gamma$ joining $v$ and $w$ such that every non-pendant vertex of $P_k$ has valency two in $\Gamma$. We shall call such a $P_k$ a {\it clean path} with respect to the shortcut $e$.
\end{defn}
Note that $k\geq 3$ is enforced to distinguish clean shortcuts from parallel edges. A clean path need not be an induced subgraph (see figure \ref{fig:shortcut} (b)). Moreover, a clean shortcut means that $\Gamma\cup e$ may be obtained by taking $(\Gamma-P_k)\cup e$ and then adding a parallel edge $e'$, replicating $e$. Subdividing $e'$ $k-2$ times results in $\Gamma\cup e$.\\
\hf Figure \ref{fig:sc-K4} depicts a planar $K_4$ obtained by iterated shortcuts, two of which are clean while the last one is not.
The planar island polynomials for the four graphs in figure \ref{fig:sc-K4} are $4+9x+6x^2+x^3$, $4+8x+4x^2+2x^3$, $4+7x+6x^2+3x^3$ and $4+6x+8x^2+4x^3$ respectively. 
\begin{figure}[!ht]
\centering
\includegraphics[scale=0.5]{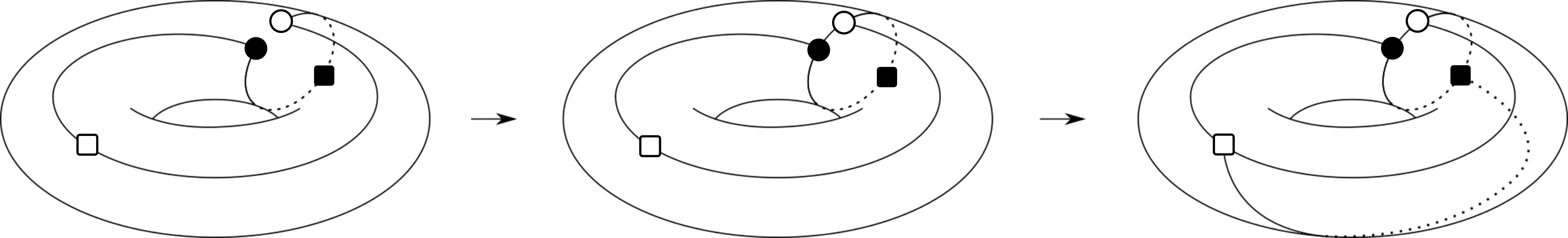}
\caption{Iterated shortcuts creating a toric $K_4$}
\label{fig:torus-sc-K4}
\end{figure}
The island polynomials for the graphs in figure \ref{fig:torus-sc-K4}, obtained by a clean shortcut followed by a shortcut, are $4x^2+x^3, 4+7x+4x^2+x^3$ and $4+6x+5x^2+2x^3$ respectively. These examples suggest that the island polynomials of graphs $\Gamma$ obtained as clean shortcuts of other graphs satisfy $\beta_\Gamma(-1)=0$. This is true in a much broader context and follows from the subsequent result.
\begin{prop}\label{prop:clean-sc}
Let $(\Gamma,\sigma)$ be a connected graph with a clean shortcut $e$ joining $v_1$ to $v_2$. Let $P_k$ denote a clean path joining $v_1$ and $v_2$. Let $\Gamma_1:=(\Gamma-P_k)\cup e$ and assume that $|V(\Gamma_1)|\geq 3$. Let $\mathfrak{c}$ be a vertex-coloring of $\Gamma$ such that all the non-pendant vertices of $P_k$ receive distinct colors different from the $N$ colors in $\mathfrak{c}(\Gamma_1)$. Let $\beta(x):=\beta_{(\Gamma_1,\sigma,\mathfrak{c})}(x)$ and $S=\{\mathfrak{c}(v_1),\mathfrak{c}(v_2)\}$.\\
(i) Let $\Gamma'$ be the graph induced on the vertex set $V(P_k)$. If $f_{\sigma}(\Gamma')=f_{\sigma}(\Gamma'-P_k)+1$, then 
\begin{equation}\label{eqn:clean-sc1}
    \beta_{(\Gamma,\sigma,\mathfrak{c})}(x) = \left\{\begin{array}{rl}
    (1+x)\beta(x)+(1+x^{|S|})(1+x)^{N-|S|} & \textup{if $k=3$}\\
    (1+x)^{k-2}\beta(x)+x^{k-3}(1+x^{|S|})(1+x)^{N-|S|} +\sum_{j\geq 2} jx^{k-2-j}(1+x)^{N-3+j} & \textup{if $k\geq 4$}
    \end{array}\right.
\end{equation}
(ii) If $\sigma(e)$ does not disconnect $(\Sigma-\sigma(\Gamma-P_k))\cup\{v_1,v_2\}$, then 
\begin{equation}\label{eqn:clean-sc2}
    \beta_{(\Gamma,\sigma,\mathfrak{c})}(x) = \left\{\begin{array}{rl}
    (1+x)\beta(x)+(1+x)^{N-|S|} & \textup{if $k=3$}\\
    (1+x)^{k-2}\beta(x)+x^{k-3}(1+x)^{N-|S|} +\sum_{j\geq 2} jx^{k-2-j}(1+x)^{N-3+j} & \textup{if $k\geq 4$}
    \end{array}\right.
\end{equation}
In particular, if $N\geq 3$, then $\beta_{(\Gamma,\sigma,\mathfrak{c})}(-1)=0$ in both cases (i) and (ii).
\end{prop}
\begin{proof}
The path graph $P_k$, along with $e$, causes a clean shortcut in $\Gamma$. Then we can define a sequence of embedded, vertex-colored graphs, starting with $\Gamma_1$, as follows:\\
\hf (1) define $\Gamma_2$ to be $\Gamma_1\cup e'$, where $e'$ is an edge parallel to $e$ (see Definition \ref{defn:sim_adj}) and embeds via the image of $P_k$;\\
\hf (2) subdivide $e'$ into two edges, identified with the image of $P_{k-1}$ and $P_1$, to obtain $\Gamma_3$;\\
\hf (3) subdivide the edge identified with $P_{k-1}$ (if $k>3$) into two edges, identified with the image of $P_{k-2}$ and $P_1$, to obtain $\Gamma_4$;\\
\hf (4) iterate this process to obtain $\Gamma_j$'s and note that $\Gamma_k=\Gamma$.\\
In case (i), i.e., when $f_{\sigma}(\Gamma')=f_{\sigma}(\Gamma'-P_k)+1$, we infer that 
$f_{\tilde{\sigma}}(\mathscr{F}_S(\Gamma_2,\mathfrak{c}))=f_{\sigma}(\mathscr{F}_S(\Gamma_1,\mathfrak{c}))+1$, where $S=\{\mathfrak{c}(v_1),\mathfrak{c}(v_2)\}$. Applying Proposition \ref{prop:nn-adj} (i), we obtain 
\begin{equation}\label{eq:G2-G1}
\beta_{(\Gamma_2,\sigma,\mathfrak{c})}(x)=\beta_{(\Gamma_1,\sigma,\mathfrak{c})}(x)+x^{|S|-1}(1+x)^{N-|S|}
\end{equation}
where $N$ is the number of colors required for $\Gamma_1$. For $\Gamma_3$, Proposition \ref{prop:subdivide} implies 
\begin{eqnarray}
\beta_{(\Gamma_3,\sigma,\mathfrak{c})}(x) & = & x\beta_{(\Gamma_2,\sigma,\mathfrak{c})}(x)+\beta_{(\Gamma_1,\sigma,\mathfrak{c})}(x)+(1+x)^{N-|S|}\nonumber \\
& = & (1+x)\beta_{(\Gamma_1,\sigma,\mathfrak{c})}(x)+(1+x^{|S|})(1+x)^{N-|S|}.\label{eqn:G3-G1}
\end{eqnarray}
where the simplification is due \eqref{eq:G2-G1}. If $\Gamma=\Gamma_3$, then we are done. Otherwise, we obtain $\Gamma_4$ by splitting the edge $[w_1v_2]$ in $\Gamma_3$ joining $v_2$ to $w_1$ at a vertex $w_2$. Note that $\Gamma_3-[w_1v_2]$ is $\Gamma_1$ with a pendant edge attached. Applying \eqref{eqn:split2} with our coloring choice, we obtain
\begin{eqnarray*}
\beta_{(\Gamma_4,\sigma,\mathfrak{c})}(x) & = & x\beta_{(\Gamma_3,\sigma,\mathfrak{c})}(x)+\beta_{(\Gamma_3-[w_1v_2],\tilde{\sigma},\hat{\mathfrak{c}})}(x)+(1+x)^{N-1}\\
& = & x\big((1+x)\beta_{(\Gamma_1,\sigma,\mathfrak{c})}(x)+(1+x^{|S|})(1+x)^{N-|S|}\big)\\
&  & +(1+x)\beta_{(\Gamma_1,\sigma,\mathfrak{c})}(x)+(1+x)^{N-1}+(1+x)^{N-1}\\
& = & (1+x)^2\beta_{(\Gamma_1,\sigma,\mathfrak{c})}(x)+x(1+x^{|S|})(1+x)^{N-|S|}+2(1+x)^{N-1}.
\end{eqnarray*}
We may iterate this process (and use \eqref{eqn:it-pendant}) to obtain, for $j\geq 4$, 
\begin{eqnarray*}
\beta_{(\Gamma_j,\sigma,\mathfrak{c})}(x) & = & (1+x)^{j-2}\beta_{(\Gamma_1,\sigma,\mathfrak{c})}(x)+x^{j-3}(1+x^{|S|})(1+x)^{N-|S|}\\
&  & +(j-2)(1+x)^{N+j-5}+(j-3)x(1+x)^{N+j-6}+\cdots+2x^{j-4}(1+x)^{N-1}.
\end{eqnarray*}
Along with \eqref{eqn:G3-G1}, the above (with $j=k$) implies \eqref{eqn:clean-sc1}. \\
\hf In case (ii), Proposition \ref{prop:nn-adj} (ii) implies that 
$$\beta_{(\Gamma_2,\sigma,\mathfrak{c})}(x)=\beta_{(\Gamma_1,\sigma,\mathfrak{c})}(x).$$
For $\Gamma_3$, Proposition \ref{prop:subdivide} implies 
\begin{eqnarray}
\beta_{(\Gamma_3,\sigma,\mathfrak{c})}(x) & = & x\beta_{(\Gamma_2,\sigma,\mathfrak{c})}(x)+\beta_{(\Gamma_1,\sigma,\mathfrak{c})}(x)+(1+x)^{N-|S|}\nonumber \\
& = & (1+x)\beta_{(\Gamma_1,\sigma,\mathfrak{c})}(x)+(1+x)^{N-|S|}.\label{eqn:G3-G1-2}
\end{eqnarray}
Applying \eqref{eqn:split2} and \eqref{eqn:pen-new-color} we obtain
\begin{equation*}
\beta_{(\Gamma_4,\sigma,\mathfrak{c})}(x)=(1+x)^2\beta_{(\Gamma_1,\sigma,\mathfrak{c})}(x)+2(1+x)^{N-1}+x(1+x)^{N-|S|}.
\end{equation*}
Iterating this process, we obtain for $j\geq 4$,
\begin{eqnarray*}
\beta_{(\Gamma_j,\sigma,\mathfrak{c})}(x) & = & (1+x)^{j-2}\beta_{(\Gamma_1,\sigma,\mathfrak{c})}(x)+x^{j-3}(1+x)^{N-|S|}\\
&  & + (j-2)(1+x)^{N+j-5}+(j-3)x(1+x)^{N+j-6}+\cdots+2x^{j-4}(1+x)^{N-1}.  
\end{eqnarray*}
Together with \eqref{eqn:G3-G1-2} and the above for $j=k$, we get \eqref{eqn:clean-sc2}.
\end{proof}
\begin{rem}
(1) Some topological hypothesis (case (i) and (ii)) in Proposition \ref{prop:clean-sc} is necessary to ensure a tractable and explicit formula. Figure \ref{fig:torus-sc2} illustrates a clean shortcut which has been subdivided with non-zero $\beta_{(\Gamma,\sigma)}(-1)$. In fact, the island polynomial for figure \ref{fig:torus-sc2} (c) is $4+7x+4x^2+2x^3$.\\
(2) We have chosen the coloring hypothesis on $P_k$ to ensure explicit formula. In principle, one can compute the colored island polynomial for a clean shortcut by iterated applications of Propositions \ref{prop:Gamma-app} (specifically Corollary \ref{cor:it-pendant}) and \ref{prop:subdivide}, once Proposition \ref{prop:nn-adj} has been applied to compute $\beta_{(\Gamma_2,\sigma,\mathfrak{c})}$. 
\end{rem}
\begin{figure}[!ht]
\centering
    \subfloat[]{{\includegraphics[scale=0.52]{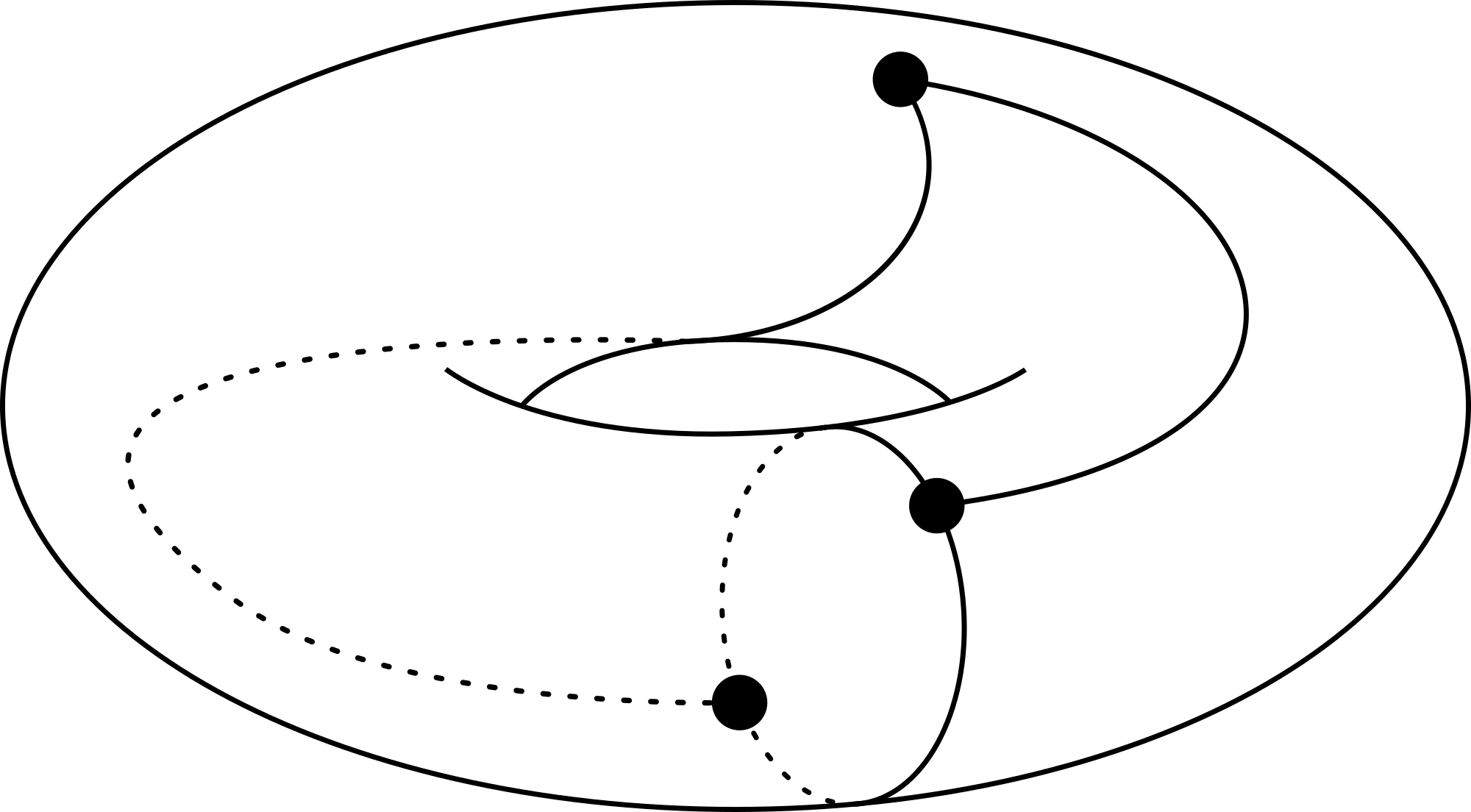} }}
    \qquad
    \subfloat[]{{\includegraphics[scale=0.52]{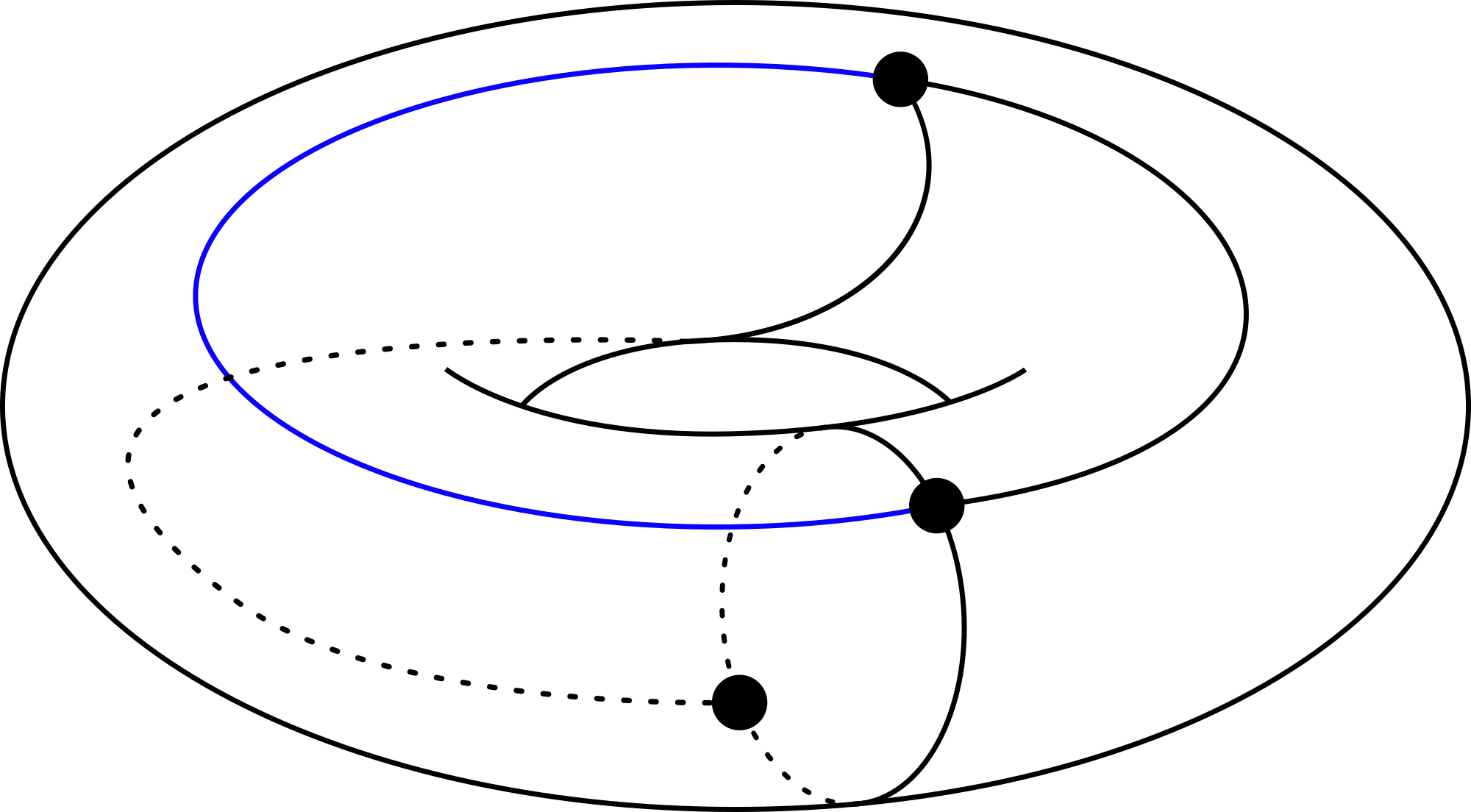} }}
    \qquad
    \subfloat[]{{\includegraphics[scale=0.52]{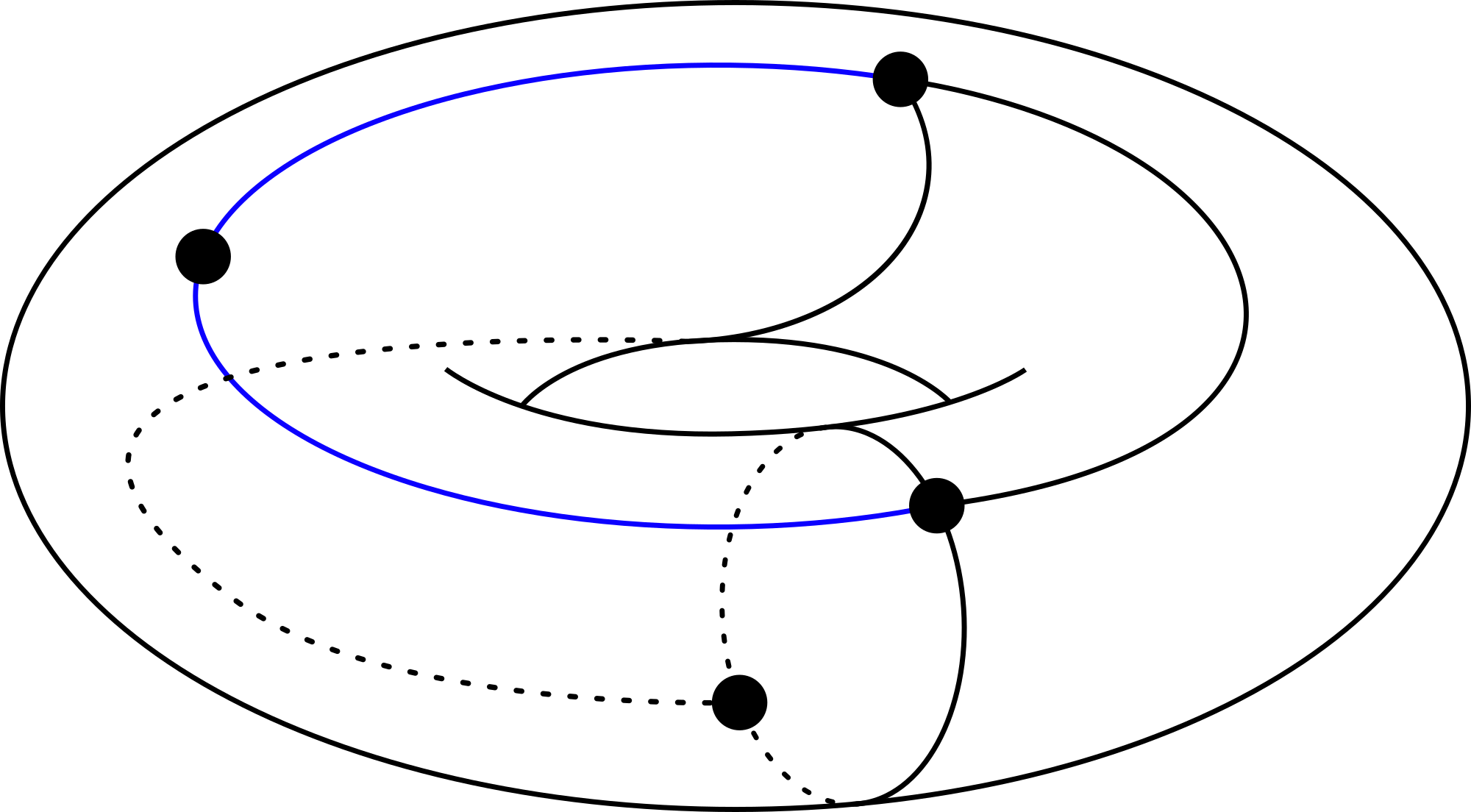} }}
\caption{(a) A toric graph (b) After a clean shortcut (c) Subdividing the shortcut}
\label{fig:torus-sc2}
\end{figure}


\subsection{Properties of the structure set}\label{subsec:str-set}

If $\Gamma$ has at least $2$ vertices, then Proposition \ref{prop:Gamma-app} implies that
$$\mathscr{S}_{d+1}(\Gamma_\textup{app})=\mathscr{S}_{d+1}(\Gamma)\cup \big(\mathscr{S}_{d+1}(\Gamma)+(1+x)^{d}\big)\cup\big(\mathscr{S}_{d}(\Gamma)(1+x)+(1+x)^{d}\big)\label{eqn:S-append}$$
This implies the following.
\begin{cor}\label{cor:str-tree}
If $\Gamma, \Gamma'$ are graphs such that $\mathscr{S}(\Gamma)=\mathscr{S}(\Gamma')$, then $\mathscr{S}(\Gamma_\textup{app})=\mathscr{S}(\Gamma'_\textup{app})$. In particular, the sets $\mathscr{S}(T)$ and $\mathscr{S}(T')$ are equal if $T$ and $T'$ are trees on the same number of vertices.
\end{cor}
\begin{eg}\label{eg:bigon}
It can be seen that 
$$\mathscr{S}(C_2)=\{1,2,x+2, 2x+2\}.$$
Conversely, if $\mathscr{S}(\Gamma)$ equals the above set, then $\Gamma$ has $2$ vertices. As $1\in \mathscr{S}(\Gamma)$, it follows that $\Gamma$ is connected and its complement (with respect to some embedding) is also connected. Thus, $\Gamma$ is built out of $P_2$, adding loops and parallel edges, say $k_1$ loops at $v_1$, $k_2$ loops at $v_2$ and $k$ parallel edges. Any such graph admits a planar embedding. A trivial coloring of $\Gamma$, embedded in the plane, would result in island polynomial $(k+k_1+k_2+1)x+(k_1+k_2+2)$. This forces $k_1=k_2=0$ and $k\in \{0,1\}$. But $k=0$ implies that $\Gamma=P_2$ and $\mathscr{S}(P_2)=\{1, x+2\}$ gives rise to a contradiction. Thus, $k=1$ and $\Gamma=C_2$.
\end{eg}
\begin{eg}
We now consider an example, quite in contrast to example \ref{eg:bigon}. Start with $P_2$ and build two graphs $\Gamma_1$ and $\Gamma_2$ by attaching one loop at each vertex and attaching two loops at one vertex respectively. These have identical structure set of polynomials, i.e.,
$$\mathscr{S}(\Gamma_j)=\{1,2,3,x+2, 2x+2, 2x+3, 3x+4\},\,\,j=1,2.$$
\begin{figure}[H]
\centering
\includegraphics[trim={0 45 0 0}, clip, scale=0.28]{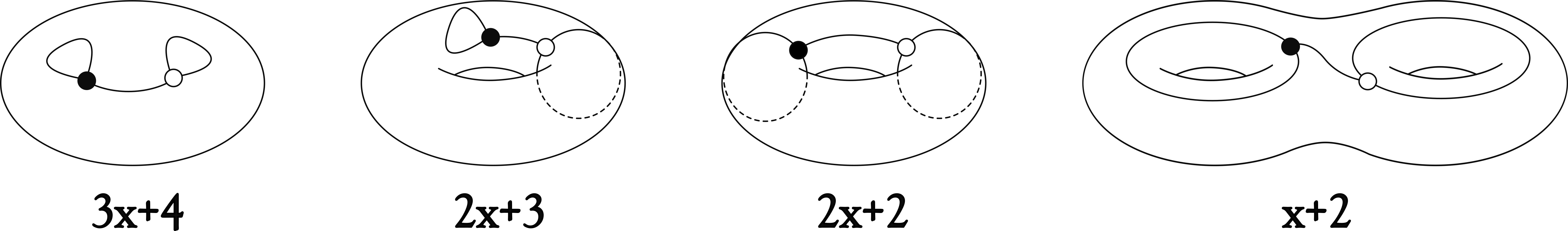}
\caption{Island polynomials of $\Gamma_1$: $3x+4$ (planar), $3x+2$, $2x+2$, $x+2$}
\label{fig:glasses-island}
\end{figure}
The structure set for the {\it theta graph} $\Gamma_3$ is
$$\mathscr{S}(\Gamma_3)=\{1,2,3,x+2, 2x+2, 3x+2\}.$$
\begin{figure}[H]
\centering
\includegraphics[trim={0 45 0 0}, clip, scale=0.35]{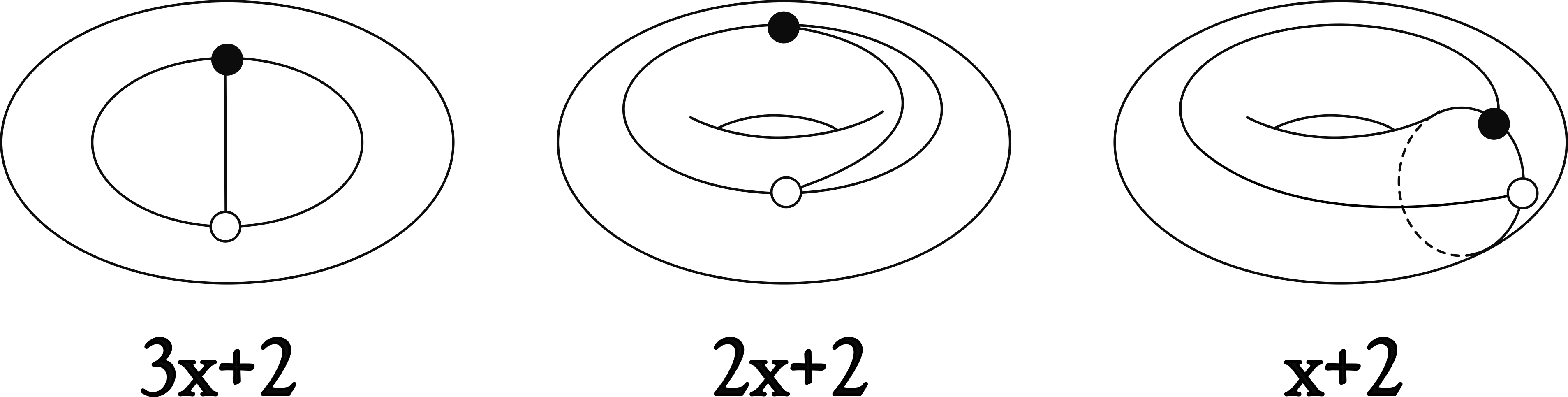}
\caption{Island polynomials of $\Gamma_3$: $3x+2$ (planar), $2x+2$, $x+2$}
\label{fig:theta-island}
\end{figure}
Note that $\Gamma_j$'s are mutually non-isomorphic and collapsing the maximal tree in $\Gamma_j$ results in $B_2$. 
\end{eg}
We have the following properties:\\
\hf (i) Lemma \ref{lmm:collapse} implies $\mathscr{S}(\Gamma/e)\subset \mathscr{S}(\Gamma)$. \\
\hf (ii) Let $\tilde{\Gamma}$ denote $\Gamma$ with a loop attached. Then Proposition \ref{prop:self-loop} implies that
$$\mathscr{S}_d(\tilde{\Gamma})\supseteq \mathscr{S}_d(\Gamma)\cup (\mathscr{S}_d(\Gamma)+(1+x)^d).$$
It is difficult (and expectedly so) to get explicit relationship between the structure set of a graph and one obtained from it (by splitting an edge or adding a parallel edge or adding a clean shortcut) due to topological constraints imposed (cf. Propositions \ref{prop:clean-sc}, \ref{prop:nn-adj}, \ref{prop:subdivide}) on the edge in order to obtain a tractable and explicit formula (see \eqref{eqn:split1}, \eqref{eqn:split2}, \eqref{eqn:parallel-edge}, \eqref{eqn:clean-sc1}, \eqref{eqn:clean-sc2}).\\
\hf A vertex-colored graph $(\Gamma,\mathfrak{c})$ is called proper if no adjacent vertices receive the same color. Given a graph $(\Gamma,\sigma,\mathfrak{c})$, we may iteratively use Lemma \ref{lmm:collapse} to obtain a properly colored embedded graph $\Gamma'$ having the same colored island polynomial as that of $\Gamma$. 
\begin{defn}\label{defn:ps-set}
Let
$$\mathscr{PS}_d(\Gamma):=\{\beta_{(\Gamma,\sigma,\mathfrak{c})}(x)\,|\,\sigma:\Gamma\hookrightarrow \Sigma, \mathfrak{c}:V(\Gamma)\twoheadrightarrow \mathcal{C}, |\mathcal{C}|=d+1,\,\textup{$\mathfrak{c}$ is proper}\}$$
denote the colored island polynomials obtained from proper vertex-colorings of $\Gamma$ using $d+1$ colors, varying over embeddings. We define the {\it proper structure set of polynomials} of $\Gamma$ to be
$$\mathscr{PS}(\Gamma):=\sqcup_{d\geq 0} \,\mathscr{PS}_d(\Gamma).\label{eqn:proper-str-set}$$
\end{defn}
\begin{rem}\label{rem:V-PS}
If $\Gamma$ is a proper vertex-colored graph on $n$ vertices, then we may embed this into a surface of high enough genus such that the face count of induced subgraph on a vertex is $1$, i.e., no loop disconnects the surface. If $\sigma$ denotes one such embedding, then $\beta_{(\Gamma,\sigma,\mathfrak{c})}(0)=n$ is the number of vertices of $\Gamma$. On the other hand, for a trivial coloring of $\Gamma$ and any embedding $\sigma$, $\beta_{(\Gamma,\sigma)}(0)\geq n$ and strict inequality holds if $\Gamma$ has any self-loop which disconnects $\Sigma$. If we choose the embedding in a high enough genus, then we can ensure that $\beta_{(\Gamma,\sigma)}(0)=n$. Thus, the number of vertices may be considered as the minimum of $\beta(0)$, where $\beta\in \mathscr{PS}(\Gamma)$. It is also the minimum of $\beta(0)$, where $\beta$ ranges over island polynomials.
\end{rem}
\begin{eg}\label{eg:cycle-tree-PS}
For any tree $T_n$ on $n\geq 2$ vertices, any proper vertex-coloring requires at least $2$ and at most $n$ colors. For the cycle graph $C_n$ on $n\geq 2$ vertices, any proper coloring requires the following:\\
\hf (i) at least $2$ and at most $n$ colors, if $n$ is even;\\
\hf (ii) at least $3$ and at most $n$ colors, if $n$ is odd.
The parity plays a significant role as there is no proper coloring of $C_{2k+1}$ with two colors. We can now use Theorem \ref{thm:tree-cycle-color-detect} (in particular \eqref{eqn:tree-color}, \eqref{eqn:cycle-color}) to conclude that
\begin{eqnarray*}
    \mathscr{PS}(T_n) & = & \left\{\begin{array}{rl}
    \{(1+x)^{n-1}+(n-1)(1+x)^{n-2},\ldots,(1+x)^2+(n-1)(1+x),(1+x)+n-1\} & \textup{$n$ even}\\
    \{(1+x)^{n-1}+(n-1)(1+x)^{n-2},\ldots,(1+x)^2+(n-1)(1+x)\} & \textup{$n$ odd}\\
    \end{array}\right.\\
    \mathscr{PS}(C_n) & = & \{\alpha x^{n-1}+n(1+x)^{n-2},\ldots,\alpha x^2+n(1+x),\alpha x +n\,|\,\alpha=1,2\}.
\end{eqnarray*}
\end{eg}


\subsection{Comparison with other invariants}\label{subsec:Krushkal}

\hf We shall list some of the properties of the Bollob\'{a}s-Riordan polynomial (see \cite[Theorem 27.5, chapter 27]{EmMof22}), following notations used therein. Let $R(G):=R(G;X,Y,Z,W)$ denote the Bollob\'{a}s-Riordan polynomial for a graph $G$. We note
that
$$R(G)=\left\{\begin{array}{rl}
R(G/e)+R(G-e) & \textup{if $e$ is an edge that is not a loop or a bridge}\\
X-R(G/e) & \textup{if $e$ is a bridge}\\
(1+Y)R(G-e) & \textup{if $e$ is an oriented trivial loop}\\
(1+WYZ)R(G-e) & \textup{if $e$ is non-oriented trivial loop}
\end{array}\right.$$
Consider the modern Krushkal polynomial $P_G$ due to Butler \cite{But18}, a polynomial in variables $X,Y,A,B$. Here we shall abbreviate $P_{G,\Sigma}$ by $P_G$. Edge contraction and deletion are related by the following formula (cf. \cite[Lemma 2.2]{Kru11}):
$$P_G=\left\{\begin{array}{rl}
P_{G/e}+P_{G-e} & \textup{if $e$ is an edge that is not a loop or a bridge}\\
(1+X)P_{G/e} & \textup{if $e$ is a bridge}\\
(1+Y)P_{G-e} & \textup{if $e$ is a loop such that $\Sigma-e$ has $2$ components}
\end{array}\right.$$
\hf We note that there is no deletion-contraction formula for colored island polynomial. In fact, for any edge $e$ (which is not a loop) in $(\Gamma,\sigma,\mathfrak{c})$, if the two vertices have the same color, then $\beta_{(\Gamma/e,\overline{\sigma},\mathfrak{c})}=\beta_{(\Gamma,\sigma,\mathfrak{c})}$ (see Lemma \ref{lmm:collapse}). Even if the two vertices have different colors, there is no simple formula relating the colored island polynomial of $\Gamma/e$, $\Gamma-e$ and $\Gamma$. If $e$ is a loop, then deleting it results in three possibilities, of which two admit explicit formula (cf. \eqref{eqn:loop-trivial}, \eqref{eqn:loop-non-trivial}) depending on whether 
\begin{equation*}
  {\beta}_{(\Gamma-e,\sigma,\mathfrak{c})}(x) = \left\{\begin{array}{rl}
  {\beta}_{(\Gamma,\sigma,\mathfrak{c})}(x)+(1+x)^{c-1} & \textup{if $\sigma(e)$ is homologically trivial}\\
  {\beta}_{(\Gamma,\sigma,\mathfrak{c})}(x) & \textup{if $\sigma(e)$ does not disconnect $\Sigma-\sigma(\Gamma-e)$}
  \end{array}\right.
\end{equation*}
where $c$ is the number of colors used in $\Gamma$. As noted following Proposition \ref{prop:self-loop}, there is no complete formula for the third possibility where a loop $\sigma(e)$ could be homologically non-trivial and disconnect $\Sigma-\sigma(\Gamma-e)$. \\
\hf For a disjoint union of two embedded graphs $(\Gamma_1,\sigma_i,\mathfrak{c}_i)$ with disjoint vertex-coloring, we have (cf. \eqref{eqn:disjoint-color})
$${\beta}_{(\Gamma_1\sqcup\Gamma_2,\sigma_1\sqcup \sigma_2,\mathfrak{c}_1\sqcup\mathfrak{c}_2)}(x)=(1+x)^{c_2}{\beta}_{(\Gamma_1,\sigma_1,\mathfrak{c}_1)}(x)+(1+x)^{c_1}{\beta}_{(\Gamma_2,\sigma_2,\mathfrak{c}_2)}(x).$$
This is analogous to the multiplicative property of $P$ for disjoint union of two graphs.


\section{Applications}
\label{sec:appl}

\hf We discuss several consequences and applications of the theory developed in \S \ref{sec:invariant}.

\subsection{Characterizing trees and cycles}\label{subsec:cycle-tree}

\subsubsection{Island polynomials}

\hf We have seen in \eqref{eqn:tree-island} that 
$\beta_{T_n}(x)=(1+x)^{n-1}+(n-1)(1+x)^{n-2}$
for any tree on $n$ vertices. In particular, the island polynomial is independent of the embedding and the combinatorial structure of the tree. Conversely, let $(\Gamma,\sigma)$ be an embedded graph with island polynomial $\beta(x)=(1+x)^{n-1}+(n-1)(1+x)^{n-2}$, i.e., $\beta$ is independent of the embedding. This implies the following:\\
\hf (a) $\Gamma$ has $n$ vertices as the degree of the island polynomial is one less that $|V(\Gamma)|$;\\
\hf (b) $\Gamma$ is connected as the top coefficient of $\beta$ is $1$;\\
\hf (c) if $\Gamma$ admits a planar embedding $\sigma$, then $S^2-\sigma(\Gamma)$ is connected, implying that $\Gamma$ has no cycles; \\
\hf (d) if there are no planar embeddings and $\Gamma$ has a cycle, let $C_r$ be one such cycle which is also a connected induced subgraph of $\Gamma$. Then one may choose an embedding in some surface of sufficiently high genus such that $\sigma(C_r)$ bounds a disk $D$ and $\sigma(\Gamma)\cap \textup{int}(D)=\varnothing$. The island polynomial for this embedding has top coefficient greater than $1$, a contradiction to (b).
\begin{theorem}\label{thm:tree-detection}
If $\Gamma$ is a graph such that the island polynomial $\beta_{(\Gamma,\sigma)}(x)=(1+x)^{n-1}+(n-1)(1+x)^{n-2}$ for any embedding, then $\Gamma$ is a tree on $n$ vertices. Moreover, if $\Gamma$ is a connected planar graph such that the island polynomial is given by
\bgd
{\beta}_{\Gamma}(x)=a(1+x)^{n-1}+b(1+x)^{n-2}
\edd
for $a\in \mathbb{N}$, then $\Gamma$ is obtained from a tree on $n$ vertices by adding $(a+b-n)$ loops and $(n-b-1)$ parallel edges. 
\end{theorem}
\begin{proof}
The properties (a)-(d) establish the first result. In fact, if the island polynomial is as claimed for any planar embedding, then by (a)-(c), $\Gamma$ is a tree on $n$ vertices. For the second result, ${\beta}_{\Gamma}(x)=a(1+x)^{n-1}+b(1+x)^{n-2}$ with $a\neq 0$ implies that $\Gamma$ has $n$ vertices. As $\beta_\Gamma(0)=a+b$, it implies that collectively the $n$ vertices have $a+b-n$ loops. We may consider the graph $\Gamma'$ obtained by deleting these loops. Using \eqref{eqn:loop-trivial} we obtain
$$\beta_{\Gamma'}(x)= (n-b)(1+x)^{n-1}+b(1+x)^{n-2}.$$
By Euler's Theorem applied to the connected planar graph $\Gamma'$, we get
$$f(\Gamma')=2+e(\Gamma')-v(\Gamma')=2+e-n$$ 
on the one hand while $f=n-b$ due to the top coefficient being $n-b$ in $\beta_{\Gamma'}(x)$. Therefore,
\begin{equation}\label{eqn:n-1}
n-1=\tfrac{e+b}{2}
\end{equation}
is the number of edges in any maximal tree in $\Gamma'$; fix one such tree $T$. Let $\rho$ denote the number of parallel edges, i.e., edges in $\Gamma'$ which are not edges in $T$ that join two adjacent vertices of $T$. Let $\epsilon$ denote the number of edges win $\Gamma'$ that are neither in $T$ and nor parallel edges, i.e., these edges join two vertices which are non-adjacent in $T$. This implies that
\begin{equation}\label{eqn:e-as-sum}
e=n-1+\epsilon+\rho=\tfrac{e+b}{2}+\epsilon+\rho.
\end{equation}
where we have used \eqref{eqn:n-1} in the last equality.\\
\hf We now compute the coefficient of $x^{n-2}$ in two ways. From the polynomial it is $n^2-n-b$. We can count $\mathscr{D}_2(\Gamma')$ directly. If $\alpha$ denotes the number of (unordered) pairs of non-adjacent vertices $(v,w)$ in $\Gamma'$, then 
\begin{equation}\label{eqn:alpha-n}
\alpha={n-1 \choose 2}-\epsilon
\end{equation}
as there are ${n-1 \choose 2}$ pairs of non-adjacent vertices in $T$. Each of the non-adjacent pairs contribute $2$ to the island boundary count while each of ${n \choose 2}-\alpha$ adjacent pairs contribute $1$ along with $\rho$, the number of parallel edges. Thus, we obtain
\begin{equation*}
n^2-n-b=2\alpha+{n\choose 2}-\alpha+\rho.
\end{equation*}
Using \eqref{eqn:alpha-n} in the above, we obtain a simplification
\begin{equation}\label{eqn:rho}
    \rho=\tfrac{e-b}{2}+\epsilon.
\end{equation}
Plugging \eqref{eqn:rho} in \eqref{eqn:e-as-sum} we obtain that $\epsilon=0$, i.e., all the edges outside of $T$ are parallel edges.
\end{proof}
The assumption of connectivity on $\Gamma$ cannot be omitted as otherwise, $\Gamma'$ (as in the proof) may be a disjoint union of $n-b$ trees (refer to Corollary \ref{cor:tree-island}, \eqref{eqn:r-comp-tree}). Thus, $\Gamma$ with the prescribed polynomial may be realized as a disjoint union of $n-b$ connected graphs, each of which is obtained from a tree by attaching loops. \\
\begin{rem}
1. If a connected embedded graph $(\Gamma,\sigma)$ has $a(1+x)^{n-1}+b(1+x)^{n-2}$ as the island polynomial (with $a>0$), then $\Gamma$ need not be obtained from a tree by adding loops and parallel edges, unless the embedding is planar; look at figure \ref{fig:island-beta-qtree} below 
\begin{figure}[!h]
    \centering
    \includegraphics[trim={0 45 0 0}, clip, scale=0.25]{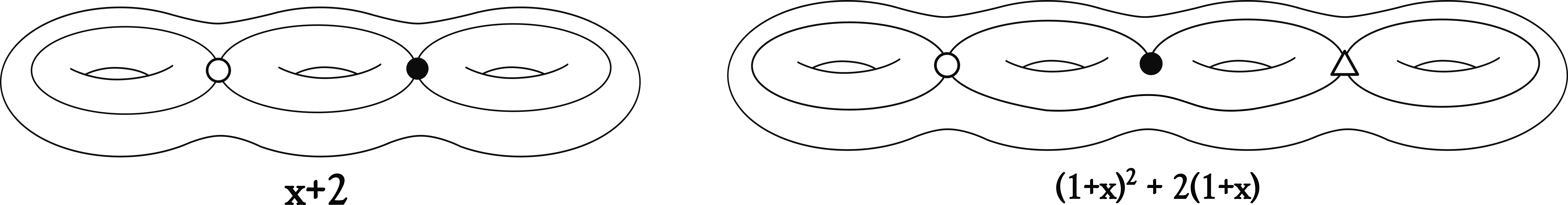}
    \caption{Graphs built out of $P_2$ and $P_3$ with island polynomials $x+2$ and $(1+x)^2+2(1+x)$ respectively }
    \label{fig:island-beta-qtree}
\end{figure}
\newline 2. There exists non-tree graphs embedded in surfaces of higher genus having island polynomials that correspond to planar trees.
\begin{figure}[!h]
    \centering
    \includegraphics[scale=0.5]{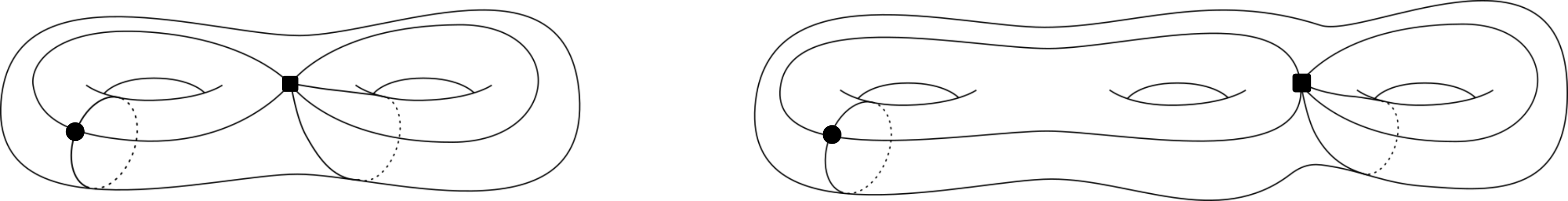}
    \caption{Two embeddings with island polynomial $x+2$}
    \label{fig:beta-same-poly}
\end{figure}
This is possible due to the extra handles present in higher genus surfaces. For instance, there is no embedding of the graph in figure \ref{fig:beta-same-poly} inside the torus with polynomial $x+2$. Note that the graph in figure \ref{fig:beta-same-poly} is built out of $P_2$ by adding three loops and one parallel edge.
\end{rem}

\hf Recall that the cycle graph $C_n$ consists of $n$ vertices labelled $1$ through $n$ such that the vertex labelled $j$ is adjacent to vertices labelled $j-1$ and $j+1$, where we are counting modulo $n$. We will assume that $n\geq 3$ whenever we discuss $C_n$.
\begin{rem}\label{rem:TEE}
We have found \cite{SiddharthaTEE2021} that the island count (which is the alternative sum of the count $\mathscr{D}_j(C_n)$) $\beta_{C_n}(-1)$ for any planar embedding, is proportional to the $n$-partite information and topological entanglement entropy. Cycle graphs are crucially related to the topological entanglement entropy measure of a collection of subsystems arranged in an annulus. 
\end{rem}
Note that the terms involved in computing $\beta_{(C_n,\sigma)}(x)$ deal with proper subgraphs of $C_n$. These subgraphs are disjoint union of trees and thus $\mathscr{D}_j(C_n,\sigma)$ is independent of the embedding if $j<n$. It is the last count $\mathscr{D}_n(C_n,\sigma)$ that detects non-triviality in homology as follows. An embedded cycle graph $\sigma(C_n)\subset \Sigma$ is an embedded loop in a connected surface. It is well-known that if $\Sigma$ is a closed\footnote{A closed manifold is a compact manifold without boundary.}, oriented surface then an embedded loop $\gamma$ disconnects $\Sigma$ if and only if the loop is homologically trivial, i.e., $\Sigma-\gamma$ has two components if and only if $[\gamma]=0\in H_1(\Sigma;\mathbb{Z})$. 
\begin{figure}[!ht]
	\centering
	\includegraphics[scale=0.6]{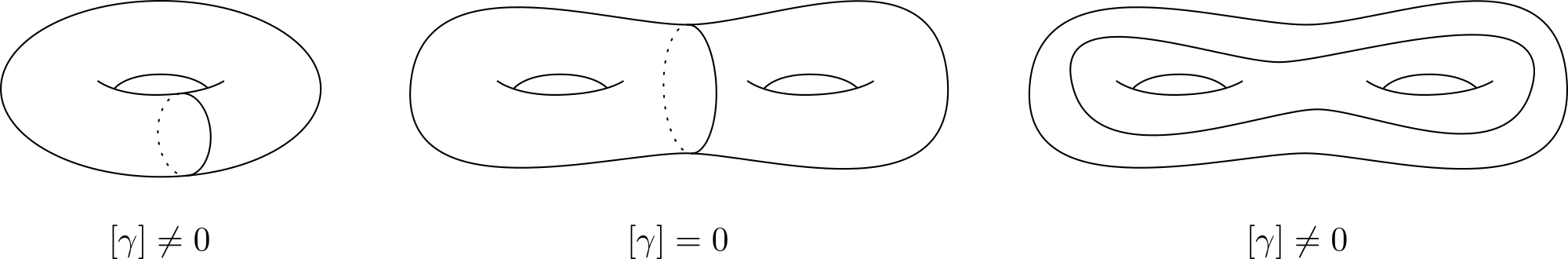}
	\caption{Various types of embedded cycles}
\end{figure}
When $\Sigma=S^2$, any embedded loop disconnects it (Jordan Curve Theorem), i.e., $\mathscr{D}_n(C_n,\sigma)=2$ for $\sigma:C_n\hookrightarrow S^2$. The coefficients $\mathscr{D}_i(C_n), i<n$ of $\beta_{C_n}$ can be computed explicitly for cycle graphs.  
\begin{theorem}\label{thm:Dnm}
The numbers $\mathscr{D}_m(C_n)$, for $1\leq m<n$, are given by
\begin{equation}
\mathscr{D}_m(C_n)=\sum_{j=0}^{m-1} (-1)^{j}~ (m-j) {n \choose m-j}=n  {n-2 \choose m-1}. \label{eq:Dcycle}
\end{equation}
\end{theorem}
The first equality in \eqref{eq:Dcycle} is the original formula that was predicted in \cite{SiddharthaTEE2021}, based on an intuitive counting method similar to the inclusion-exclusion principle. However, a simple proof using inclusion-exclusion seems elusive. The proof of Theorem \ref{thm:Dnm} uses elementary methods but is slightly lengthy; we have moved it to appendix \ref{cycle-app}.
\begin{theorem}\label{thm:cycle-detect}
The island polynomial for $(C_n,\sigma)$ is given by
\begin{equation}\label{eqn:cycle-pol}
{\beta}_{(C_n,\sigma)}(x)=\left\{\begin{array}{rl}
n(1+x)^{n-2}+2x^{n-1} & \textup{if $\Sigma-\sigma(C_n)$ has $2$ components}\\
n(1+x)^{n-2}+x^{n-1} & \textup{otherwise}.
\end{array}\right.
\end{equation}
If $\Gamma$ is a graph such that the island polynomials satisfy
$$\left\{\beta_{(\Gamma,\sigma)}(x)\,|\,\sigma:\Gamma\hookrightarrow \Sigma\right\}\subseteq\{x^{n-1}+n(1+x)^{n-2},2x^{n-1}+n(1+x)^{n-2}\}=\mathscr{S}(C_n)$$
then $\Gamma$ is a cycle on $n$ vertices.
\end{theorem}
\begin{proof}
The island polynomial for $C_n$ \eqref{eqn:cycle-pol} follows from Theorem \ref{thm:Dnm}. For the second claim, if $n=2$, then the result follows from Example \ref{eg:bigon}. Let $n\geq 3$ for the rest of the proof. The graph $\Gamma$ has $n$ vertices as either polynomial is of degree $n-1$. As the coefficient of $x^{n-2}$ in a possible $\beta_{(\Gamma, \sigma)}$ is $n$, for each vertex $v\in V(\Gamma)$, the graph $\Gamma-v$ is connected and $\Sigma-\sigma(\Gamma-v)$ is also connected. As $n\geq 3$, choose two distinct vertices $v_1, v_2$. For any $w\in V(\Gamma)-\{v_1,v_2\}$, $w$ is connected to $v_1$ in $\Gamma-v_2$ and $w$ is connected to $v_2$ in $\Gamma-v_1$. Thus, $\Gamma$ is connected. For either of the possible polynomials, it takes the value $n$ at $x=0$. Thus, $\Gamma$ has no loops at any vertex as if it did, then one can construct an embedding $\sigma$ with $\beta_{(\Gamma,\sigma)}(0)>n$, a contradiction. \\
\hf The coefficient for $x$ is $n(n-2)$; we are crucially using $n\geq 3$ here. As the number of faces $f_\sigma(\Gamma)\leq 2$, there can be at most one parallel edge. If there are no parallel edges and $\Gamma$ has $r$ edges, then 
\bgd
\mathcal{D}_1(\Gamma)=2\bigg({n \choose 2}-r\bigg)+r=n(n-2).
\edd
This implies that $r=n$. Now, either $\Gamma$ is a cycle or $\Gamma$ is obtained from a maximal tree $T$ by adding an edge joining two distinct vertices. As $n\geq 3$, in the latter case, $\Gamma$ will have a pendant vertex. By \eqref{eqn:pen-new-color}, $\beta_{(\Gamma,\sigma,\mathfrak{c})}(-1)=0$ which contradicts the possible values $\{(-1)^{n-1}, 2(-1)^{n-1}\}$. Thus, $\Gamma$ is a cycle.\\
\hf If there was a parallel edge $e'$, then consider $\Gamma':=\Gamma-e'$. We may embed $\Gamma'$ and then extend it to an embedding of $\Gamma$ such that $\sigma(e')$ disconnects $\Sigma-\sigma(\Gamma')$. Then \eqref{eqn:parallel-edge} implies that
$$\beta_{(\Gamma',\sigma,\mathfrak{c})}\in \{x^{n-1}+n(1+x)^{n-2}-x(1+x)^{n-2},x^{n-1}+n(1+x)^{n-2}-x(1+x)^{n-2}\}.$$
The coefficient of $x$ is $n(n-2)-1$ and we get that $\Gamma$ has $n+1$ edges. Thus, $\Gamma$ is obtained from a tree (on $n$ vertices) by adding a parallel edge and another edge, i.e., $\Gamma$ is obtained from a unicyclic graph by adding a parallel edge. As unicyclis graphs are planar, $\Gamma$ admits planar embeddings. But any such embedding of $\Gamma$ would have $f_\sigma(\Gamma)=3$, a contradiction. 
\end{proof}
\begin{rem}
Induction allows us to prove
$$\sum_{j=1}^m jx^{m-j}(1+x)^{j-1}=(m-x)(1+x)^m+x^{m+1}.$$
This may be used in \eqref{eqn:clean-sc1} (see Proposition \ref{prop:clean-sc}) to simplify the colored island polynomial for a planar graph $\Gamma$ built out of adding a pendant edge $e$ to a trivially colored $C_k$, where the new vertex is colored the same as its neighbour. On the one hand, this gives us
$$\beta_{(\Gamma,\mathfrak{c})}(x)=k(1+x)^{k-2}+2x^{k-1}$$
while on the other hand, by Corollary \ref{cor:collapse-tree}   
$$\beta_{C_k}(x)=\beta_{\Gamma/e}(x)=\beta_{(\Gamma,\mathfrak{c})}(x).$$
This gives a completely different proof of \eqref{eqn:cycle-pol}.
\end{rem}
\begin{rem}
It is known that the Tutte polynomial for any tree $T_n$ on $n$ vertices is given by
$$T_{T_n}(x,y)=x^{n-1}$$
and if $T_\Gamma(x,y)=x^n$, then $\Gamma$ is a tree on $n$ vertices. The Tutte polynomial for the cycle graph $C_n$ is
$$T_{C_n}(x,y)=y+x+x^2+\cdots+x^{n-1}$$
and it detects cycle graphs.    \\
\hf The Krushkal polynomial for $T_n$ embedded in a surface of genus $g$ is given by
$$P_{T_n}(X,Y,A,B)=(1+X)^{n-1}B^g.$$
The Krushkal polynomial for the cycle graph $C_n$, embedded in a genus $g$ surface $\Sigma$, is
$$P_{C_n}(X,Y,A,B)=(X^{n-1}+nX^{n-2}+\cdots+n)B^g+Y^{\alpha-1}B^{g+\alpha-2}$$
where $\alpha$ is the number of path components of $\Sigma-C_n$.
\end{rem}

\subsubsection{Colored island polynomials}

\hf We shall deal with proper vertex-colorings. Recall that given a graph $(\Gamma,\sigma,\mathfrak{c})$, we may iteratively use Lemma \ref{lmm:collapse} to obtain a properly colored embedded graph $\Gamma'$ having the same colored island polynomial as that of $\Gamma$. We want to generalize Theorems \ref{thm:tree-detection} and \ref{thm:cycle-detect} to the colored island polynomial setting.
\begin{theorem}\label{thm:tree-cycle-color-detect}
Let $(T_n,\mathfrak{c})$ be proper vertex-colored tree on $n$ vertices. Then the colored island polynomial is given by
\begin{equation}\label{eqn:tree-color}
{\beta}_{(T_{n},\sigma,\mathfrak{c})}(x)= (1+x)^{N-1}+(n-1)(1+x)^{N-2}.
\end{equation}
where $N=|\mathfrak{c}(V(T_n))|$.\\
\hf Let $(C_n,\mathfrak{c})$ be proper vertex-colored cycle on $n$ vertices. Then the colored island polynomial is given by
\begin{equation}\label{eqn:cycle-color}
{\beta}_{(C_{n},\sigma,\mathfrak{c})}(x)=\alpha x^{N-1}+n(1+x)^{N-2}.
\end{equation}
where $N=|\mathfrak{c}(V(C_n))|$ and $\alpha$ is the number of components of $\Sigma-\sigma(C_n)$.
\end{theorem}
\begin{proof}
We observe that a proper vertex-colored tree is built out of $P_1$ by adding pendant edges in one of the following ways:\\
\hf (I) the vertex being added is colored by a new color (case (i) of Proposition \ref{prop:Gamma-app});\\
\hf (II) the vertex being added is colored by a color already used but not the color of its (only) adjacent vertex (case (iii) of Proposition \ref{prop:Gamma-app}). \\
In fact, the given tree is built out of a specific sequence of $c$ type I operations and $n-c$ type II operations. Let $\Gamma_{\textup{I}\circ \textup{II}}$ (respectively $\Gamma_{\textup{II}\circ \textup{I}}$) denote the colored graph obtained from $\Gamma$ by applying II followed by I (respectively I followed by II). It follows from Proposition \ref{prop:Gamma-app} that the colored island polynomials of $\Gamma_{\textup{I}\circ \textup{II}}$ and $\Gamma_{\textup{II}\circ \textup{I}}$ are identical. Thus, to compute ${\beta}_{(T_{n},\sigma,\mathfrak{c})}(x)$, we may compute it for a tree built out of $T_c$ (obtained by $c-1$ type I operations) with trivial coloring (with colored island polynomial $\beta(x)=(1+x)^{c-1}+(c-1)(1+x)^{c-2}$) and then applying II $n-c$ times. Proposition \ref{prop:Gamma-app}, \eqref{eqn:pen-same-non-adj-color} applied $(n-c)$ times to $\beta(x)$ implies \eqref{eqn:tree-color}.\\
\hf Consider $(C_n,\mathfrak{c})$, a cycle graph on $n$ vertices equipped with a proper vertex-coloring requiring $N$ colors. Let $e$ denote the edge joining vertex $1$ and $n$. Then $C_n-\{e\}$ is a path graph $P_n$ on $n$ vertices with an induced proper vertex-coloring $\mathfrak{c}$ requiring $N$ colors. We compare $\mathscr{D}_i(C_n,\sigma,\mathfrak{c})$ with $\mathscr{D}_i(P_n,\sigma,\mathfrak{c})$.\\
\hf (i) $\mathscr{D}_1(C_n,\sigma,\mathfrak{c})=n=\mathscr{D}_1(P_n,\sigma,\mathfrak{c})$ due to the coloring being proper;\\
\hf (ii) for $N>i\geq 2$, we have
\begin{eqnarray*}
    \mathscr{D}_i(C_n,\sigma,\mathfrak{c}) & = & \sum_{\mathscr{F}_i(C_n,\mathfrak{c})}f_\sigma(\Gamma')\\
    & = & \sum_{\Gamma'\in\mathscr{F}_i(C_n,\mathfrak{c}),\,\{1,n\}\subseteq V(\Gamma')}f_\sigma(\Gamma')+\sum_{\Gamma'\in\mathscr{F}_i(C_n,\mathfrak{c}),\,\{1,n\}\not\subseteq V(\Gamma')}f_\sigma(\Gamma')\\
    & = & \sum_{\Gamma'-e\in\mathscr{F}_i(P_n,\mathfrak{c}),\,\{1,n\}\subseteq V(\Gamma'-e)}\big(f_\sigma(\Gamma'-e)-1\big)+\sum_{\Gamma'\in\mathscr{F}_i(P_n,\mathfrak{c}),\,\{1,n\}\not\subseteq V(\Gamma')}f_\sigma(\Gamma')\\
    & = & \mathscr{D}_i(P_n,\sigma,\mathfrak{c})-{N-2 \choose i-2}
\end{eqnarray*}
Thus, we can compute colored island polynomial for $(C_n,\sigma,\mathfrak{c})$ as follows:
\begin{eqnarray*}
\beta_{(C_n,\sigma,\mathfrak{c})}(x) & = & \alpha x^{N-1}+\sum_{i=1}^{N-1} \mathscr{D}_i(C_n,\sigma,\mathfrak{c})x^{i-1}\\
& = & \alpha x^{N-1}+\sum_{i=1}^{N-1} \mathscr{D}_i(P_n,\sigma,\mathfrak{c})x^{i-1}-\sum_{i=1}^{N-1}{N-2 \choose i-2}x^{i-1}\\
& = & \alpha x^{N-1}+\sum_{i=1}^{N-1} \mathscr{D}_i(P_n,\sigma,\mathfrak{c})x^{i-1}-\sum_{j=1}^{N-3}{N-2 \choose j}x^{j+1}\\
& = & \alpha x^{N-1}+\sum_{i=1}^{N} \mathscr{D}_i(P_n,\sigma,\mathfrak{c})x^{i-1}-\sum_{j=1}^{N-2}{N-2 \choose j}x^{j+1}\\
& = & \alpha x^{N-1}+(1+x)^{N-1}+(n-1)(1+x)^{N-2}-x(1+x)^{N-2}\\
& = & \alpha x^{N-1}+n(1+x)^{N-2}
\end{eqnarray*}
where we have used \eqref{eqn:tree-color} in the second-last equality. 
\end{proof}
\begin{cor}
The structure set for any tree $T_n$ on $n$ vertices is given by
\begin{equation}\label{eqn:str-tree}
\mathscr{S}(T_n)=\{(1+x)^j+k(1+x)^{j-1}\,|\,0\leq j \leq n-1,\,j\leq k\leq n-1\}.
\end{equation}
Conversely, if $\mathscr{S}(\Gamma)\subseteq \mathscr{S}(T_n)$ then $\Gamma=T_k$ for some $k\leq n$, with $\Gamma=T_n$ if and only if $\mathscr{S}(\Gamma)= \mathscr{S}(T_n)$.\\
\hf The structure set for the cycle $C_n$ on $n$ vertices is given by
\begin{equation}\label{eqn:str-cycle}
\mathscr{S}(C_n)=\{\alpha x^j+k(1+x)^{j-1}\,|\,\alpha\in\{1,2\},\,0\leq j \leq n-1,\,j+1\leq k\leq n\}.
\end{equation}
Conversely, if $\mathscr{S}(\Gamma)\subseteq \mathscr{S}(C_n)$ then $\Gamma=C_k$ for some $k\leq n$, with $\Gamma=C_n$ if and only if $\mathscr{S}(\Gamma)= \mathscr{S}(C_n)$.
\end{cor}
\begin{proof}
The claims about the structure sets follow from Theorem \ref{thm:tree-cycle-color-detect} and any vertex-colored graph can be transformed by iterative edge collapses to a proper vertex-colored graph without changing the colored island polynomial (see Lemma \ref{lmm:collapse}). In the case of $T_n$ (respectively $C_n$), such a reduction is another tree (respectively a cycle) on smaller number of vertices. \\
\hf For the converse part, we first prove it for $T_n$. Suppose that $\Gamma$ has $k$ vertices; the island polynomials of $\Gamma$ are of degree $k-1$. As $\mathscr{S}(\Gamma)\subseteq \mathscr{S}(T_n)$, by \eqref{eqn:str-tree}, the island polynomials of $\Gamma$ must be of the form 
$$\beta_{(\Gamma,\sigma)}(x)\in \{(1+x)^{k-1}+l(1+x)^{k-2}\,|\,n-1\geq l\geq k-1\}.$$
The only way to obtain $k=\beta_{(\Gamma,\sigma)}(0)$ is by $\beta_{(\Gamma,\sigma)}(x)=(1+x)^{k-1}+(k-1)(1+x)^{k-2}$. By remark \ref{rem:V-PS}, $k$ is realized as the minimum among $\beta(0)$ among island polynomials. This forces $\Gamma$ to have $(1+x)^{k-1}+(k-1)(1+x)^{k-2}$ as one of its island polynomial. If there is another possible island polynomial $\gamma(x)=(1+x)^{k-1}+j(1+x)^{k-2}$, with $j\geq k$, then $\gamma(0)=j+1>k$, implying that there are loops in $\Gamma$. We now look at $f_\sigma(\Gamma)$ across embeddings. By the existence of loops $f_\sigma(\Gamma)\geq 2$ for a suitable embedding. These $f_\sigma(\Gamma)$'s correspond to the constant polynomials obtained by uni-colorings of $\Gamma$, i.e., elements of $\mathscr{S}_1(\Gamma)\subseteq \mathscr{S}_1(T_n)$. As $\mathscr{S}(T_n)=\{1\}$, we arrive at a contradiction. Thus, $\Gamma$ has $(1+x)^{k-1}+(k-1)(1+x)^{k-2}$ as its island polynomial; it is independent of the embedding. Thus, by Theorem \ref{thm:tree-detection}, $\Gamma$ is a tree on $k$ vertices.\\
\hf Suppose $\mathscr{S}(\Gamma)\subseteq \mathscr{S}(C_n)$ and $\Gamma$ has $k$ vertices; the island polynomials of $\Gamma$ are of degree $k-1$. As $\mathscr{S}(\Gamma)\subseteq \mathscr{S}(C_n)$, by \eqref{eqn:str-cycle}, the island polynomials of $\Gamma$ must be of the form 
$$\beta_{(\Gamma,\sigma)}(x)\in \{\alpha x^{k-1}+l(1+x)^{k-2}\,|\,\alpha\in\{1,2\},\,n\geq l\geq k\}.$$
If the set of island polynomials of $\Gamma$ is contained in $\{\alpha x^{k-1}+k(1+x)^{k-2}\,|\,\alpha=1,2\}$, then by Theorem \ref{thm:cycle-detect} we conclude that $\Gamma$ is the cycle graph on $k$ vertices. This is certainly the case when $k=1,n$. Otherwise, assume that $1<k<n$ and $\beta_{(\Gamma,\sigma)}(x)=\alpha x^{k-1}+l(1+x)^{k-2}$ for some $n\geq l>k$ and form some embedding $\sigma$. As $\beta_{(\Gamma,\sigma)}(0)=l>k$, this implies that $\Gamma$ has $l-k$ loops. Thus, $f_\sigma(\Gamma) \geq  l-k+1$ for a suitable embedding. These $f_\sigma(\Gamma)$'s correspond to the constant polynomials obtained by uni-colorings of $\Gamma$, i.e., elements of $\mathscr{S}_1(\Gamma) \subseteq \mathscr{S}_1(C_n)=\{1,2\}$. This forces $l=k+1$, i.e., $\Gamma$ is connected and has exactly one loop $e$. If $\Gamma'=\Gamma-e$, then it follows that $\Sigma-\sigma(\Gamma')$ is connected for any $\sigma$. This implies that $\Gamma'$ has no parallel edges or cycles (and hence so does $\Gamma$). Thus, $\Gamma'$ is a tree on $k$ vertices and $\Gamma$ is a graph obtained by attaching a loop to $\Gamma'$. Iteratively using \eqref{eqn:pen-new-color}, we obtain the island polynomial for $\Gamma$ (thought of as built out of the loop $e$ and then adding pendant edges in succession)
$$\beta_{(\Gamma,\sigma)}(x)=\alpha(1+x)^{k-1}+(k-1)(1+x)^{k-2}$$
where $\alpha$ is the number of components of $\Sigma-\sigma(e)$. We now have an equality of sets
$$\{(1+x)^{k-1}+(k-1)(1+x)^{k-2}, 2(1+x)^{k-1}+(k-1)(1+x)^{k-2}\}= \{x^{k-1}+l(1+x)^{k-2}, 2 x^{k-1}+l(1+x)^{k-2}\}.$$
This is a contradiction; a way to conclude this is to evaluate both sides at $x=-1$. Thus, the set of island polynomials of $\Gamma$ is contained in $\{\alpha x^{k-1}+k(1+x)^{k-2}\,|\,\alpha=1,2\}$ and we are done.
\end{proof}

\subsection{Tree-cycle graphs}\label{subsec:tree-cycle-graphs}

\hf We introduce and discuss a special class of planar graphs for which ${\beta}_\Gamma(-1)$ vanishes. On the one hand, graphs of the form $\Gamma_1\vee\Gamma_2$, with $v(\Gamma_i)\geq 2$, have this property; this is a consequence of \eqref{eqn:wedge}. We shall denote by $\mathcal{W}$ the collection of all wedge sum of graphs having two or more vertices. On the other hand, consider a graph $(\Gamma,\sigma)$ with an edge $e$, which is not a loop, between $v$ and $w$. We may add another edge $e'$ between $v$ and $w$ (parallel edge) to create $\Gamma_0(e')$ and then subdivide it $k\geq 1$ times, resulting in a graph $\Gamma_k(e')$. We call such a graph as obtained from $\Gamma$ via {\it iterated split parallel edge}. The edge $e'$ now looks like $P_{k+2}$. If we have chosen an extension of this $P_{k+2}$ into $S^2$, then we may compare ${\beta}_{\Gamma_k(e')}$ with ${\beta}_{\Gamma}$. In the language of \S \ref{subsec:arrange_shortcircuit}, $P_{k+2}$ is a clean path between $v$ and $w$ while $e$ is called a clean short circuit. It follows from Proposition \ref{prop:clean-sc} that 
\begin{equation}\label{eqn:ssa}
{\beta}_{\Gamma_k(e')}(-1)=0
\end{equation}
if $k\geq 1$. Let $\mathcal{S}$ denote the collection of graphs obtained by repeated applications of iterated split parallel edge. 
\begin{figure}[!ht]
\centering
\includegraphics[scale=0.5]{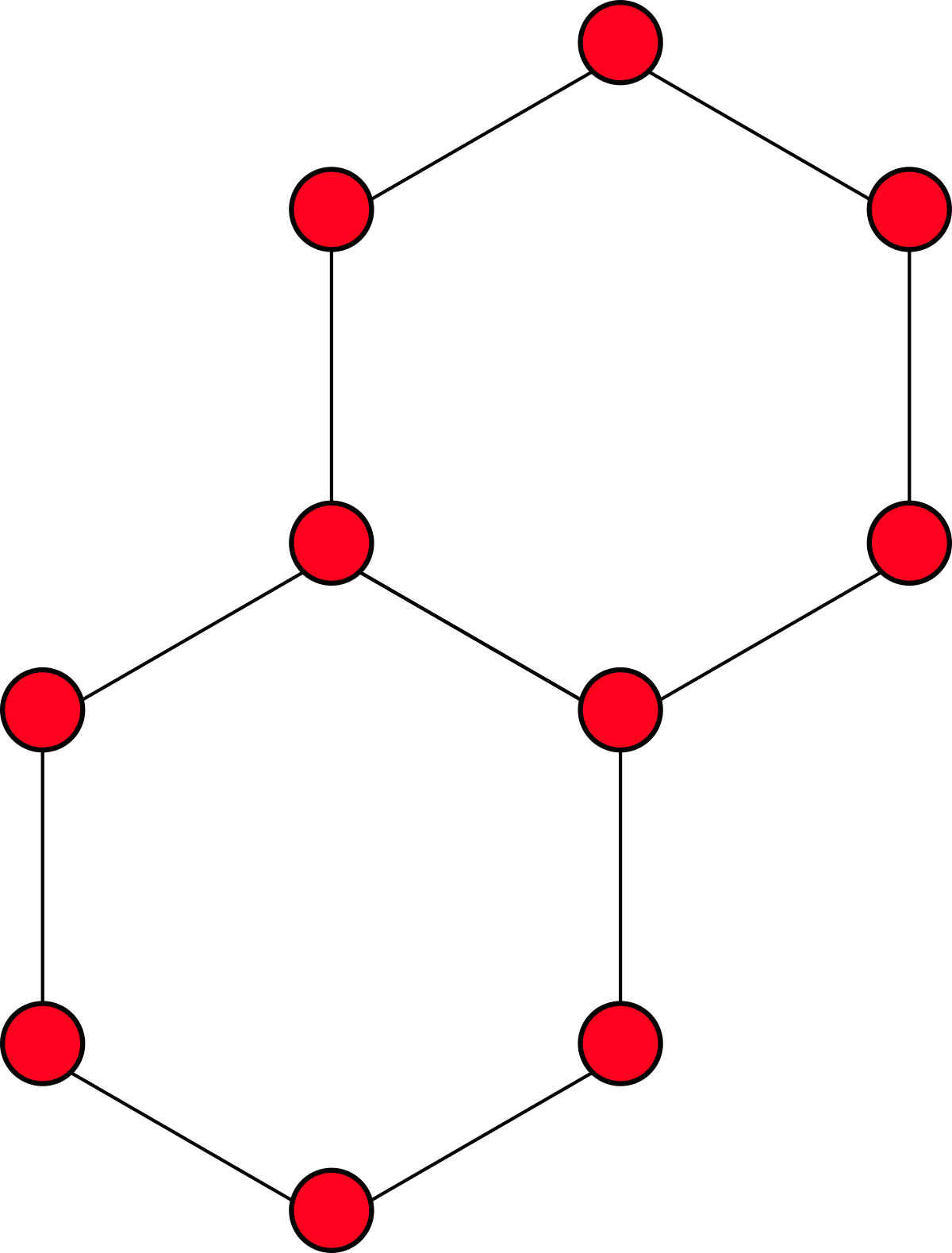}
\caption{This graph is in $\mathcal{S}$ but not in $\mathcal{W}$}
\label{fig:split-sim-adj}
\end{figure}
By \eqref{eqn:ssa}, any $\Gamma\in \mathcal{S}$ satisfies ${\beta}_\Gamma(-1)=0$. Note that (figure \ref{fig:split-sim-adj}) neither $\mathcal{S}$ or $\mathcal{W}$ is contained in the other. As the example presented in figure~\ref{fig:C_nonzero} shows,
\begin{figure}[!ht]
\centering
\includegraphics[scale=0.5]{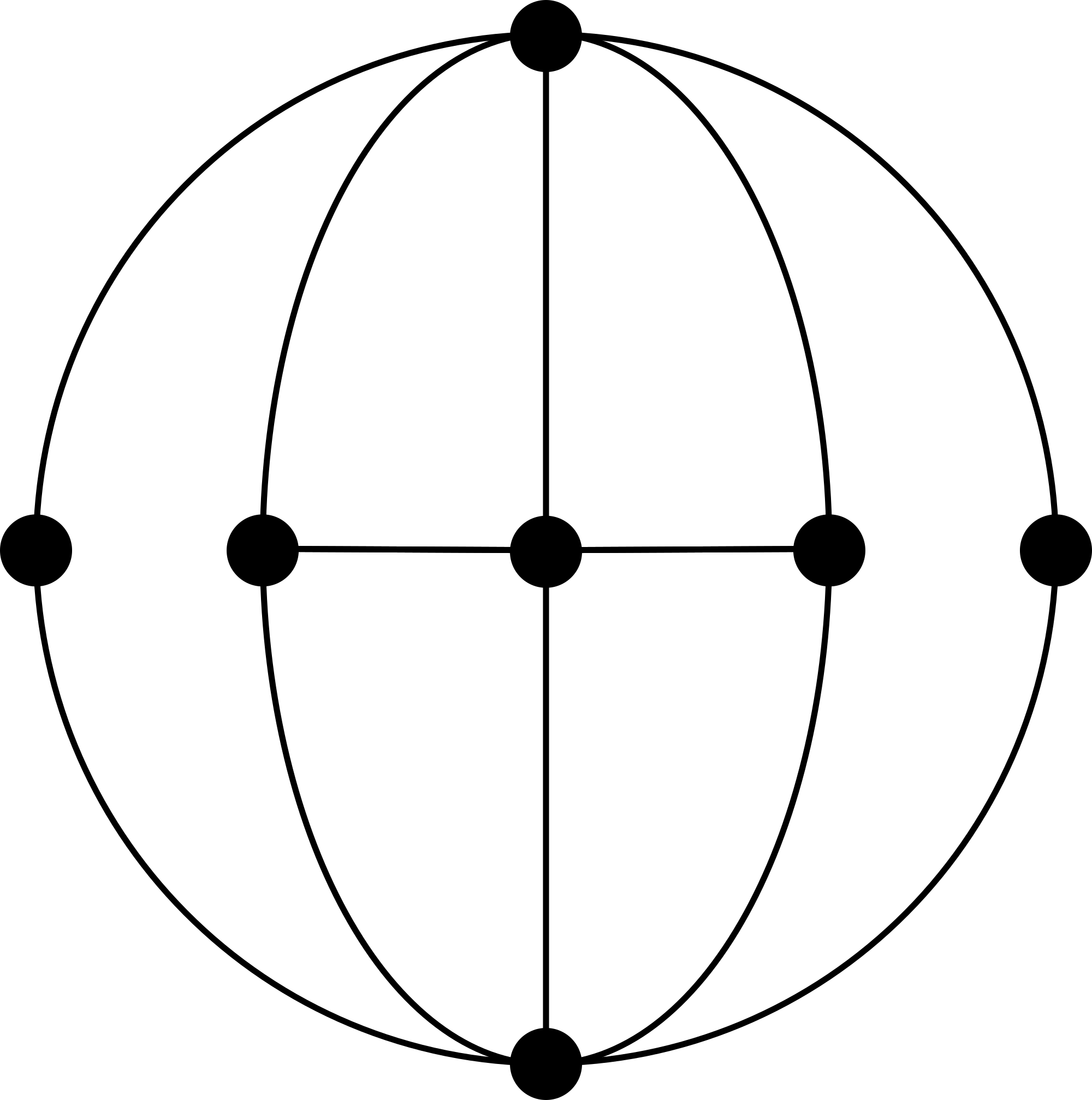}
\caption{The graph $\Gamma$ has non-zero ${\beta}_\Gamma(-1)$.}
\label{fig:C_nonzero}
\end{figure}
there are graphs which are not in $\mathcal{W}\cup\mathcal{S}$ with non-zero total connected induced subgraph count. 

\hf The intersection $\mathcal{W}\cap\mathcal{S}$ is non-empty and is an interesting collection. We shall, however, analyze a subset of $\mathcal{W}$ which is constructive in nature.
\begin{figure}[!h]
    \centering
    \includegraphics[scale=0.6]{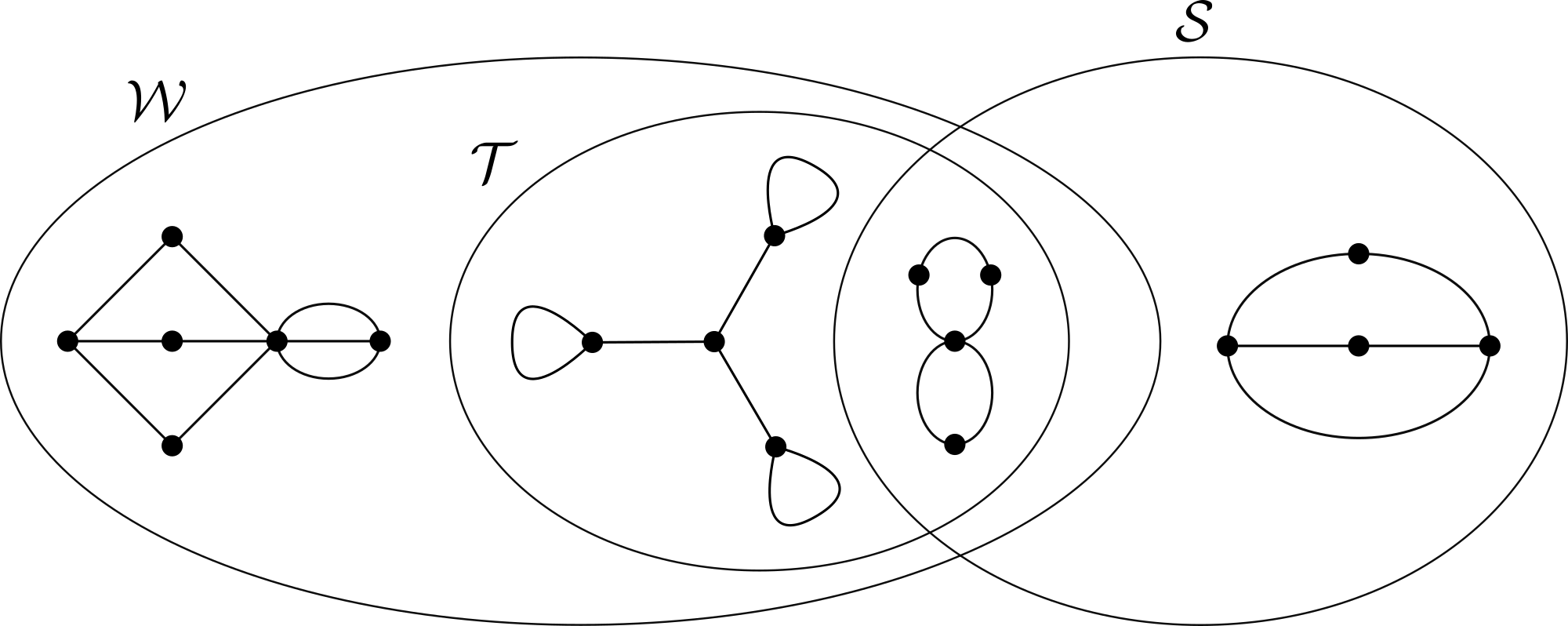}
    \caption{A Venn diagram indicating the relationship between  $\mathcal{W}, \mathcal{S}$ and $\mathcal{T}$}
    \label{fig:beta-zero}
\end{figure}
We first consider the following two operations on a graph~$\Gamma$.\\
\hf \textbf{Operation I}: Add a loop $e$ and subdivide this loop $k\geq 0$ many times, i.e., $e$ after this iterated subdivision looks like $C_{k+1}$ and the resulting graph is $\Gamma\vee C_{k+1}$.\\
\hf \textbf{Operation II}: Add a parallel edge $e$ (see Definition \ref{defn:sim_adj}) and subdivide this new edge $k\geq 0$ many times, i.e., $e$ after this iterated subdivision looks like $P_{k+2}$. \\
An operation of type I or II is called admissible. Note that in operation I (resp. II), we may add a loop (resp. parallel edge) and not subdivide it. This corresponds to $k=0$.
\begin{defn}\label{tree-cycle}
A {\it tree-cycle graph} is a planar graph obtained from a (finite) tree $T$, with $v(T)\geq 3$, by applying a finite sequence of admissible operations on it. We shall denote by $\mathcal{T}$ the collection of all tree-cycle graphs.
\end{defn}
By definition, all trees are tree-cycle graphs as no operations of type I or II are used. Note that we have already encountered tree-cycle graphs in the proof of Theorem \ref{thm:tree-detection}, where the island polynomial detects tree-cycle graphs formed out of operations I and II, but no edge subdivision is used in either operations. In fact, operations I and II without edge subdivisions were also used, in reverse, in the proof of Theorem \textcolor{red}{\ref{thm:cycle-detect}}.
\begin{figure}[H]
    \centering
    \includegraphics[scale=0.5]{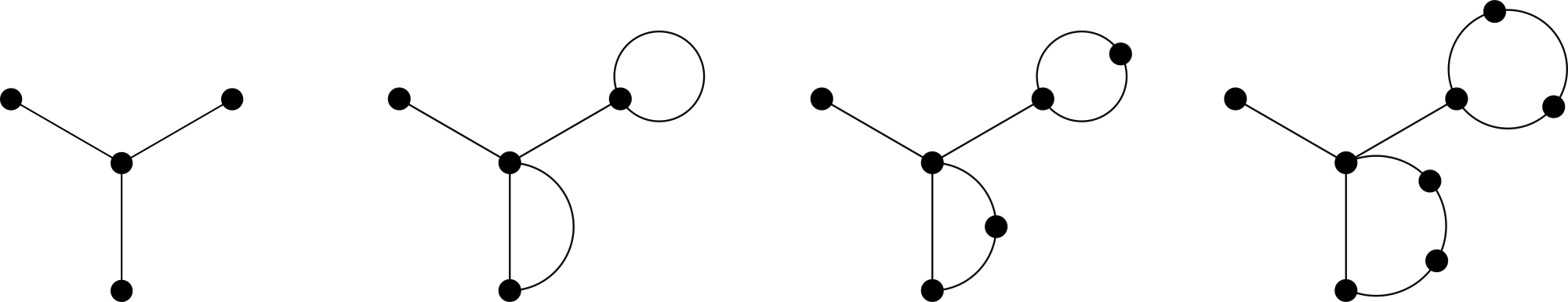}
    \caption{Creation of a tree-cycle graph}
    \label{fig:tree_cycle-cr}
\end{figure}
The graphs in figure \ref{fig:bigon} are
\begin{figure}[!h]
    \centering
    \includegraphics[scale=0.5]{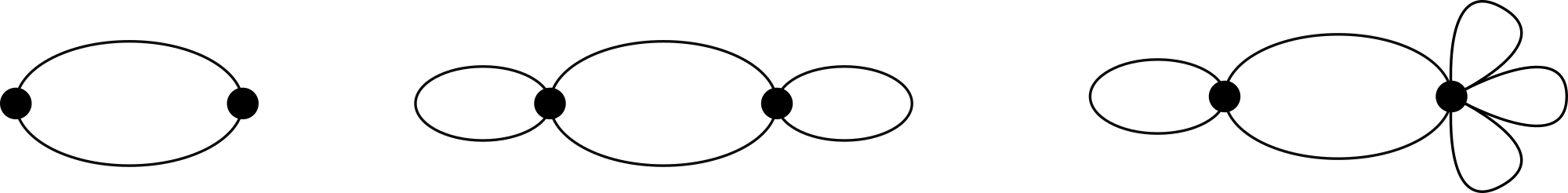}
    \caption{Graphs built out of operations I and II from $P_2$}
    \label{fig:bigon}
\end{figure}
not considered tree-cycle graphs although all of these satisfy ${\beta}_\Gamma(-1)=0$.
\begin{prop}\label{prop:tree-cycle-gr}
The island count for any tree-cycle graph is zero.
\end{prop}
There are two ways of proving the claim. Any tree-cycle graph is in $\mathcal{W}$ and ${\beta}_\Gamma(-1)=0$ for $\Gamma\in \mathcal{W}$. The direct intuitive explanation for Proposition \ref{prop:tree-cycle-gr} is the following. Subdividing a loop is equivalent to adding an appendix (which kills the total connected induced subgraph count) and then replicating this appendix. Adding a parallel edge and then subdividing it is equivalent to creating a short circuit. Both these operations, for planar graphs, kill the island count. Note that Proposition \ref{prop:tree-cycle-gr} can be generalized to embedded tree-cycle graphs with the appropriate assumptions similar to those appearing in Propositions \ref{prop:self-loop}, \ref{prop:nn-adj} and \ref{prop:clean-sc}.

\subsection{Topological entanglement entropy}\label{subsec:TEE}

\hf The island polynomial $\beta_{(\Gamma,\sigma)}$ has applications in theoretical physics and applied branches of science. The island count $\beta_{(\Gamma,\sigma)}(-1)$, for a planar embedding, was first used in a quantum condensed matter physics problem and called the {\it total connected induced subgraph count} in that context. The application of quantum information-theoretic tools in understanding strongly correlated many-body systems has increased drastically in recent decades. Von Neumann entanglement entropy is one such measure that quantifies the entanglement of a region of a system with the rest. Due to the presence of strong correlations, standard perturbation theory cannot be used to study such strongly correlated systems. A special class of these strongly correlated systems are topologically ordered phases. Long range entanglement entropy and dependence on the topology of the underlying manifold on which this is embedded are some of the important features of these phases. In this subsection, we supress the embedding $\sigma$ as we deal with planar graphs and denote the island polynomial by $\beta_{\Gamma}$. \\
\hf Von Neumann entanglement entropy $S_A$ of a subsystem $A$, in two spatial dimension, measures the entanglement between the subsystem $A$ and the rest. In a general system, this measure $S_A$ is dependent on the geometry of the subsystem $A$. Topologically order systems are special calls of systems. This entanglement measure $S_A$ has a piece that depends on the topology of the subsystem $A$ and the topology of the underlying manifold, called Topological Entanglement Entropy (TEE) \cite{Kitaev_Preskill_2006,Levin_Wen_2006,SiddharthaTEE2021}, along with a geometry dependent piece. Symbolically
\begin{eqnarray}
S_A=\alpha L_A -\gamma_A~,~~\gamma_A=\mathcal{J}_A \Omega
\label{eq:S_A}
\end{eqnarray}
where $L_A$ is the perimeter of the 2-dimensional subsystem $A$ and $\gamma_A$ is the geometry independent piece, $\mathcal{J}_A$ is the number of boundary components of the subsystem $A$ and $\Omega$ is a characteristic of the topologically ordered phase. To capture the topological entanglement entropy, multipartite information was defined such that all the geometry-dependent pieces cancel with each other. 
\begin{eqnarray*}
I^N_\mathcal{A} &=& \sum_{i=1}^{N} (-1)^{i-1} \sum_{Q\in P_i(\mathcal{A}) } S_{Q}
\end{eqnarray*}
where $\mathcal{A}=\{A_1,A_2,\cdots,A_N\}$ represents the collection of $N$ subsystems, and $P_i(\mathcal{A})$ is the set of all possible subsets with $i$ number of subsystems in it. It was shown that for this particular entanglement measure the net geometry dependent piece (first term in the right hand side of the first equation in \eqref{eq:S_A}) is zero. Thus the problem of entanglement entropy can be solved by properly computing the number of boundary components.
\begin{eqnarray*}
I^N_{\mathcal{A}} &=& -\Omega \sum_{i=1}^{N} (-1)^{i-1} \sum_{Q\in P_i(\mathcal{A}) } \mathcal{J}_{Q}
\end{eqnarray*}

\begin{figure}[!ht]
\centering
\includegraphics[scale=0.4]{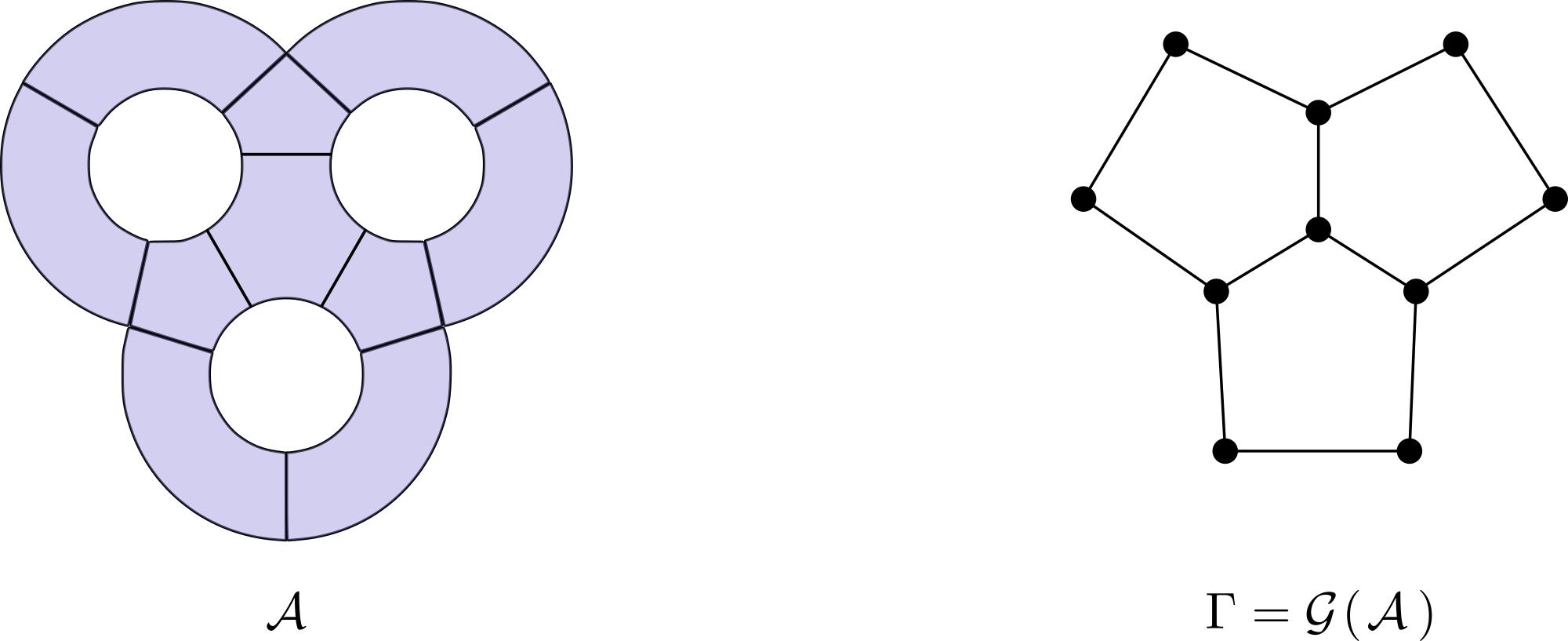}
\caption{Converting a CSS to a graph}
\label{fig:css-to-graph}
\end{figure}
\hf This multipartite entanglement problem can be expressed in the language of graph theory, where each subsystem is represented by a vertex and the connection between two subsystems is represented by an edge (figure \ref{fig:css-to-graph}). We focus on planar graphs. If $n_h$ is the number of holes in the graph, then $n_h+1=f$, where $f$ is the number of faces (connected components) in the complement of the graph. In figure  \ref{fig:css-to-graph} the number of holes $n_h=3$. For planar graphs, $I^N_{\mathcal{A}} = -\Omega {\beta}_{\Gamma}(-1)$, whence one can describe the vanishing of multipartite information measure $I^N_{\mathcal{A}}$ by showing that ${\beta}_{\Gamma}(-1)=0$. 

\begin{figure}[!h]
\centering
\includegraphics[scale=0.44]{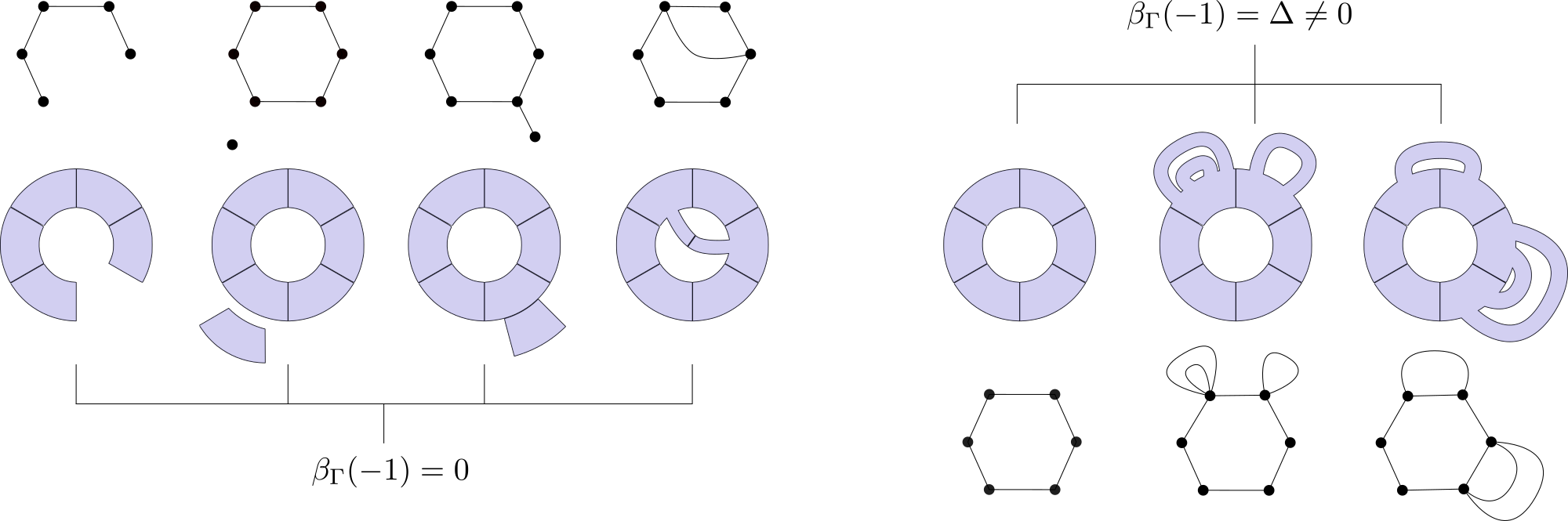}
\caption{Summary of the effect of certain transformations}
\label{fig:summary}
\end{figure}
We summarize the results for various types of CSS (collection of subsystems) and their corresponding graphs in figure \ref{fig:summary}. We know that ${\beta}_{\Gamma}(-1)$ and the corresponding multipartite information measure are both non-zero for cycle graphs with at least $3$ vertices. 
This is robust against certain transformations like adding loops or adding parallel edges between adjacent vertices (rightmost figures in figure \ref{fig:summary}). On the other hand, the island count is zero (leftmost figures in figure \ref{fig:summary}) for path graph, a union of disjoint connected graphs, graphs with a pendant edge or having a clean shortcut.\\
\hf Apart from applications to topological entanglement entropy, our graph-theoretic formalism has potential applications in machine learning \cite{ML_srinivasan2018quantifying,ML_zhang2017network} and artificial intelligence \cite{AI_ibeas2006artificial,AI_kannaiyan2021review}, both of which are data-driven methods with the potential to solve many complex problems by pattern recognition. Storing data and their interconnections are often represented in terms of a graph. Storing data in subgraphs can help drastically simplify the problem \cite{ko2017graph}. We believe that the colored island polynomial will be computationally useful and relevant in the context of these problems.

\subsection{Emergence of Euler characteristic}\label{subsec:Euler} 

\begin{figure}[!ht]
\centering
\includegraphics[scale=0.5]{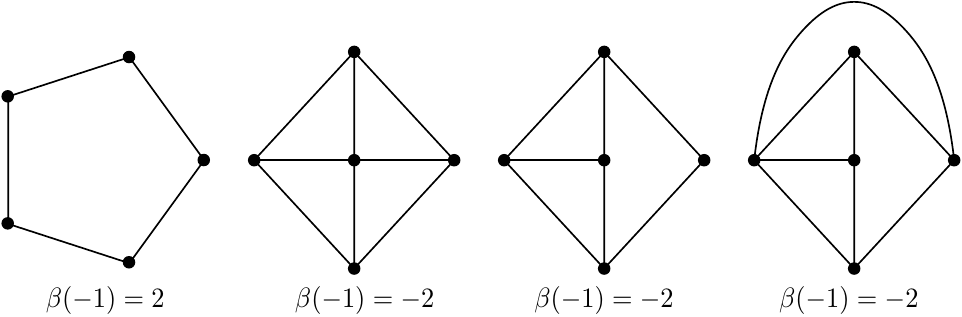}
\caption{Graphs on $5$ vertices with non-zero $\beta(-1)$.} 
\label{fig:beta-5}
\end{figure}

\hf Recall (cf. Definition \ref{defn:beta} and \eqref{eqn:col-isl-pol}) that the island polynomial for an embedded graph $(\Gamma,\sigma)$ on $n$ vertices is given by
\begin{equation}\label{eqn:beta-rec}
\beta_{(\Gamma,\sigma)}(x)=\sum_{i=1}^n \Big(\sum_{\Gamma'\in\mathscr{F}_i(\Gamma)} f_\sigma(\Gamma')\Big)x^{i-1}=f_\sigma(\Gamma)x^{n-1}+\sum_{i=1}^{n-1} \Big(\sum_{\Gamma'\in\mathscr{F}_i(\Gamma)} f_\sigma(\Gamma')\Big)x^{i-1}.
\end{equation}
This can be rewritten, recursively using \eqref{eqn:beta-rec}, as 
\begin{eqnarray*}
\beta_{(\Gamma,\sigma)}(x) & = &  f_\sigma(\Gamma)x^{n-1}+\sum_{\Gamma'\in\mathscr{F}_{n-1}(\Gamma)} f_\sigma(\Gamma')x^{i-1}+\sum_{i=1}^{n-2} \Big(\sum_{\Gamma'\in\mathscr{F}_i(\Gamma)} f_\sigma(\Gamma')\Big)x^{i-1}\\
& = & f_\sigma(\Gamma)x^{n-1}+\sum_{\Gamma'\in\mathscr{F}_{n-1}(\Gamma)} \Big(\beta_{(\Gamma',\sigma)}(x)-
\sum_{j=1}^{n-2} \Big(\sum_{\Gamma''\in\mathscr{F}_j(\Gamma')} f_\sigma(\Gamma'')\Big)x^{j-1}\Big)+
\sum_{i=1}^{n-2} \Big(\sum_{\Gamma'\in\mathscr{F}_i(\Gamma)} f_\sigma(\Gamma')\Big)x^{i-1}\\
& = & f_\sigma(\Gamma)x^{n-1}+\sum_{\Gamma'\in\mathscr{F}_{n-1}(\Gamma)} \beta_{(\Gamma',\sigma)}(x)-
\sum_{j=1}^{n-2} \Big(\sum_{\Gamma'\in\mathscr{F}_i(\Gamma)} (n-j-1)f_\sigma(\Gamma')\Big)x^{i-1}\\
\end{eqnarray*}
where the factor of $(n-j-1)$ appears because a subgraph $\Gamma''$ on $j$ vertices lies inside $n-j$ subgraphs with $n-1$ vertices. Iterating this we obtain
\begin{equation}\label{eqn:beta-rec-full}
\beta_{(\Gamma,\sigma)}(x) = f_\sigma(\Gamma)x^{n-1}+ \Big(\sum_{\Gamma'\in\mathscr{F}_{n-1}(\Gamma)} \beta_{(\Gamma',\sigma)}(x)\Big)- \Big(\sum_{\Gamma'\in\mathscr{F}_{n-2}(\Gamma)} \beta_{(\Gamma',\sigma)}(x)\Big)+\cdots+(-1)^{n-2}\Big( \sum_{\Gamma'\in\mathscr{F}_{1}(\Gamma)} \beta_{(\Gamma',\sigma)}(x)\Big).
\end{equation}
In particular, the island count for $\Gamma$ is, up to an additive factor of $\pm f_\sigma(\Gamma)$, the alternating sum of the island counts of its subgraphs. Recall that graphs have island count zero if one of the following holds:\\
\hf (a) it is not connected \eqref{eqn:disjoint-pol}\\
\hf (b) has at least three vertices and a pendant vertex (Proposition \ref{prop:Gamma-app}, \eqref{eqn:pen-new-color})\\
\hf (c) is a wedge sum of graphs (Corollary \ref{cor:wedge})\\
\hf (d) is a bridge graph (Proposition \ref{prop:bridge-del}, \eqref{eqn:bridge-del})\\
We assume that $\Gamma$ is connected and cellularly embedded in a connected surface $\Sigma$. Let $v, e$ and $f$ denote the number of vertices, edges and faces of $\Gamma\hookrightarrow \Sigma$ respectively. The only contributions to the island count $\beta_{(\Gamma,\sigma)}(-1)$ come from the following: \\
\hf (i) $(-1)^n v$, the last term in the right hand side (RHS) of \eqref{eqn:beta-rec-full}\\
\hf (ii) $(-1)^{n-1}f$, the first term in the RHS of \eqref{eqn:beta-rec-full}\\
\hf (iii) $(-1)^{n-1}e$, the second last term in the RHS of \eqref{eqn:beta-rec-full}\\
\hf (iv) $\beta_{\Gamma'}(-1)$ for subgraphs $\Gamma'$ not satisfying (a)-(d) and $\Gamma'\in\mathscr{F}_j(\Gamma)$, where $3\leq j\leq n-1$. \\
The contribution from (i) to (iii) is
\bgd
(-1)^n(v-e-f)=(-1)^n(v-e+f)-(-1)^n 2f=(-1)^n(\chi -2f),
\edd
where $\chi$ is the Euler characteristic of $\Sigma$. \\
\hf \textbf{Planar graphs}: The contributions from (iv) for planar graphs include signed counts of $k$-cycles in $\Gamma$. In particular, if $c_k$ denotes the number of $k$-cycles, then 
\begin{eqnarray}
(-1)^n\beta_\Gamma(-1) & = & \chi-2f+\sum_{j=3}^{n-1}\sum_{\Gamma'\in \mathscr{F}_j}(-1)^{j-1}\beta_{\Gamma'}(-1)\nonumber \\
& = & \chi-2f+ \beta_{C_3}(-1)c_3+\sum_{j=4}^{n-1}\sum_{\Gamma'\in \mathscr{F}_j}(-1)^{j-1}\beta_{\Gamma'}(-1). \nonumber \\
& = & \chi-2f+2c_3-((-2)c_4+2c_{K_4})+\sum_{j=5}^{n-1}\sum_{\Gamma'\in \mathscr{F}_j}(-1)^{j-1}\beta_{\Gamma'}(-1) \label{eqn:Euler}
\end{eqnarray}
where $c_{K_4}$ is the number of $K_4$'s in $\Gamma$ and $\beta_{K_4}(-1)=2$. One can expand this further with a lot of care. For instance, the only planar graphs on $5$ vertices with non-zero island count (in particular, these do not satisfy (a)-(d)) are depicted in figure \ref{fig:beta-5}. Thus, the contribution from $\mathscr{F}_5(\Gamma)$ would be counting induced subgraphs in $\Gamma$ of one of the types shown in figure \ref{fig:beta-5} with weights $2, -2, -2$ and $-2$ respectively.\\
\hf \textbf{Cellularly embedded graphs $(\Gamma,\sigma)$}: If $\Gamma$ has no loops or parallel edges, then we have a version of \eqref{eqn:Euler}. Consider evaluating $(-1)^n\beta_{(\Gamma,\sigma)}(-1)$ via \eqref{eqn:beta-rec-full}:
\begin{eqnarray*}\label{eqn:Euler-em}
(-1)^n\beta_{(\Gamma,\sigma)}(-1) & = & -f+ \sum_{\Gamma'\in \mathscr{F}_{1}(\Gamma)}\beta_{(\Gamma',\sigma)}(-1)-\sum_{\Gamma'\in \mathscr{F}_{2}(\Gamma)}\beta_{(\Gamma',\sigma)}(-1)+\sum_{j=3}^{n-1}\sum_{\Gamma'\in \mathscr{F}_{j}(\Gamma)}(-1)^{j-1}\beta_{(\Gamma',\sigma)}(-1)\\
& = & -f+ v-\sum_{\Gamma'\in \mathscr{F}_{2}(\Gamma)}\beta_{(\Gamma',\sigma)}(-1)+\sum_{j=3}^{n-1}\sum_{\Gamma'\in \mathscr{F}_{j}(\Gamma)}(-1)^{j-1}\beta_{(\Gamma',\sigma)}(-1)\\
& = & -f+ v-e+\sum_{j=3}^{n-1}\sum_{\Gamma'\in \mathscr{F}_{j}(\Gamma)}(-1)^{j-1}\beta_{(\Gamma',\sigma)}(-1)
\end{eqnarray*}
where $\beta_{(\Gamma',\sigma)}(x)$ is either $2(1+x)$ or $2+x$, depending on whether $\Gamma'\in\mathscr{F}_2(\Gamma)$ is $E_2$ or an edge of $\Gamma$ respectively. This implies 
$$\chi-2f=(-1)^n\beta_{(\Gamma,\sigma)}(-1) + \sum_{\Gamma'\in \mathscr{F}_{n-1}(\Gamma)}(-1)^{n-1}\beta_{(\Gamma',\sigma)}(-1)+\cdots+\sum_{\Gamma'\in \mathscr{F}_{3}(\Gamma)}(-1)^{3}\beta_{(\Gamma',\sigma)}(-1)
$$
From a computational perspective (and certainly from a topological entanglement entropy point of view), we can recover $\chi-2f$ as the alternating sum of island counts of all induced subgraphs of $\Gamma$ on at least $3$ vertices. In physics parlance, $\chi-2f$ is an emergent feature arising from computing topological entanglement entropy on all subsystems.

\subsection{An important transformation for detecting changes in topology} \label{subsec:pair-of-pants}

\begin{figure}[H]
\centering
\includegraphics[scale=0.45]{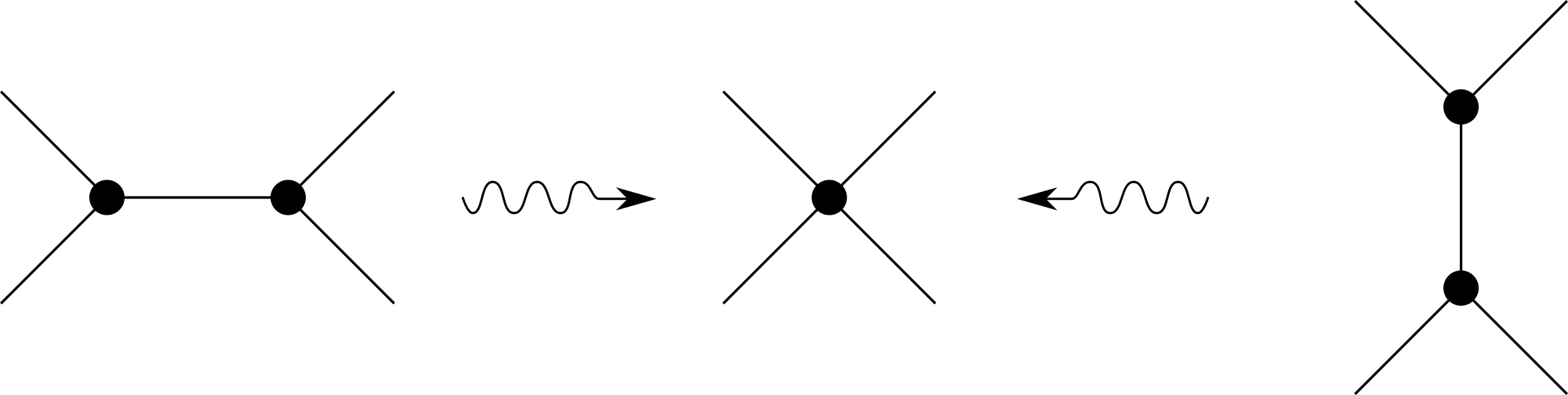}
\caption{A common degeneration of two structures}
\label{fig:>-<I}
\end{figure}

\hf An important structural transformation that appears in many contexts, including homotopy associativity of algebraic structures~\cite{HAH1, HAH2}, pair of pants decomposition of Riemann surfaces~\cite{Rat06}, Feynman diagrams in physics (e.g., the four-point correlation in the $s$ and $t$ channels~\cite{peskinschroeder}), Morse theory and singularity theory~\cite{bhattacharjee2017topology,strogatz2018nonlinear}, reconnection of vortices in quantum fluids~\cite{fonda2019}, Lifshitz transitions of the Fermi volume~\cite{lifshitz1960,volovik2007,volovik2017} and quantum transport of electrons in wire networks~\cite{lal2007,lal2008}, is given in figure \ref{fig:>-<I}.

\begin{figure}[!h]
\centering
\includegraphics[scale=0.45]{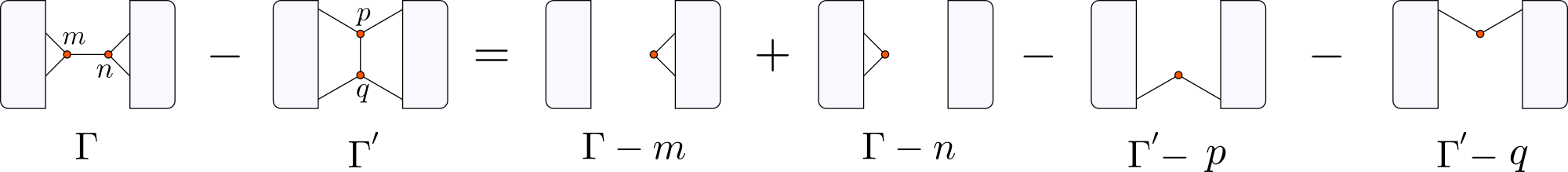}
\caption{Difference in the island polynomials}
\label{fig:>-<}
\end{figure}

\hf We shall analyze what happens under this transformation in the context of island polynomials associated to graphs. Let $\Gamma$ and $\Gamma'$ be graphs obtained from $\Gamma_1\sqcup\Gamma_2$ by attaching bridges of type I and type II respectively (see figure \ref{fig:>-<}). As all our embeddings are planar, we shall use $\beta_\Gamma$ to denote $\beta_{(\Gamma,\sigma)}$. Using Theorem \ref{thm:color} we get, 
\begin{eqnarray}
{\beta}_{\Gamma}(x)&=& x{\beta}_{\Gamma_{m=n}}(x)+{\beta}_{\Gamma-m}(x)+{\beta}_{\Gamma-n}(x)-(1+x){\beta}_{\Gamma-\{ m,n\}}(x) \nonumber\\
{\beta}_{\Gamma'}(x)&=& x{\beta}_{\Gamma'_{p=q}}(x)+{\beta}_{\Gamma'-p}(x)+{\beta}_{\Gamma'-q}(x)-(1+x){\beta}_{\Gamma'-\{p,q\}}(x). \nonumber
\end{eqnarray}
As $\Gamma-\{m,n\}=\Gamma-\{p,q\}$ and $\Gamma_{m=n}=\Gamma_{p=q}$, taking the difference of the two equations above we get
\begin{equation}
\label{eqn:>-<}
{\beta}_{\Gamma}(x)-{\beta}_{\Gamma'}(x) =  \Big[ {\beta}_{\Gamma-m}(x)+{\beta}_{\Gamma-n}(x) \Big] - \Big[{\beta}_{\Gamma'-p}(x)+{\beta}_{\Gamma'-q}(x) \Big]
\end{equation}
The difference in the island polynomials of $\Gamma$ and $\Gamma'$ is more succinctly depicted in figure \ref{fig:>-<}. 
We shall record the following observations:\\
\hf (a) The difference is a polynomial of lower degree with zero as a root.\\
\hf (b) As $\Gamma'-p$ and $\Gamma'-q$ are both obtained by subdividing the bridge for $\Gamma_1\sqcup\Gamma_2$, they have identical polynomials. In fact, using \eqref{eqn:bridge-del} and \eqref{eqn:disjoint-color}, followed by \eqref{eqn:split2}, we can infer that
\bgd
{\beta}_{\Gamma'-p}(x)=(1+x)^{n_1+1}{\beta}_{\Gamma_2}(x)+(1+x)^{n_2+1}{\beta}_{\Gamma_1}(x)+(1-x)(1+x)^{n_1+n_2-1}.
\edd
\hf (c) Note that $\Gamma-m$ and $\Gamma-n$ need not have the same polynomial. However, as $\Gamma -m$ is a disjoint union of at least two graphs, by \eqref{eqn:disjoint-pol}, the island polynomial has $(1+x)$ as a factor.\\
Combining the last two observations, we conclude that ${\beta}_\Gamma(-1)={\beta}_{\Gamma'}(-1)$. 
\begin{eg}
Consider a special case depicted below. 
\begin{figure}[!h]
    \centering
    \includegraphics[scale=0.5]{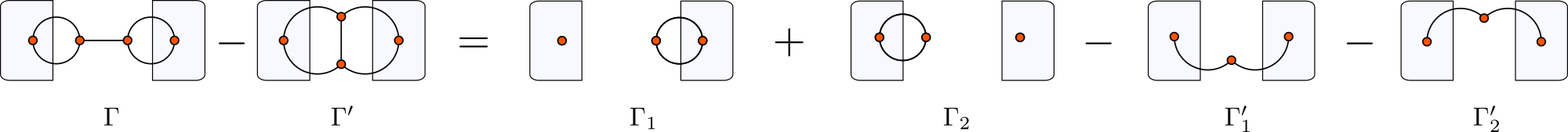}
    \caption{Difference between ``pair of glasses" graph and ``theta" graph}
    \label{fig:coloring_simplecase_1}
\end{figure}
For a planar embedding of the graph in figure~\ref{fig:coloring_simplecase_1}, we calculate the polynomials on the right side of the identity in figure \ref{fig:coloring_simplecase_1} to be
\begin{eqnarray*}
{\beta}_{\Gamma_1}(x)={\beta}_{\Gamma_2}(x) &=& 3+6x+3x^2 \nonumber\\
{\beta}_{\Gamma'_1}(x) = {\beta}_{\Gamma'_2}(x) &=& 3+4x+x^2
\end{eqnarray*}
This implies that
\bgd
{\beta}_{\Gamma}(x)-{\beta}_{\Gamma'}(x) = 4x(x+1)
\edd
and the polynomials for $\Gamma$ and $\Gamma'$ take the same value at $x=0,-1$.
\begin{figure}[!h]
    \centering
    \includegraphics[scale=0.7]{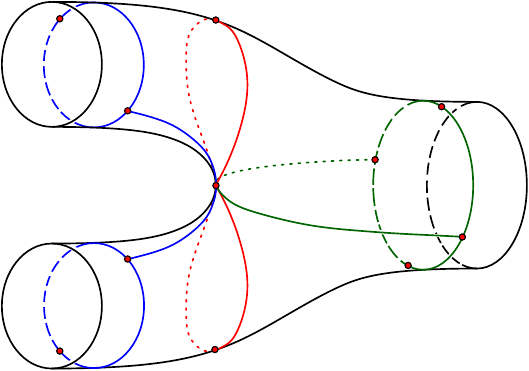}
    \caption{The transformation viewed via a pair of pants ala Morse theory}
    \label{fig:pair-of-pants}
\end{figure}
The red graph (on three vertices) in figure \ref{fig:pair-of-pants} is a degeneration of both the blue graph as well as the green graph, both on four vertices. Note that in Morse theory, the change in topology going from the left pair of circles to the right circle is detected by the index of the saddle point. Such saddle points are known to characterize transitions in various contexts in physics. For instance, in dynamical systems with a few degrees of freedom, a diagram such as figure \ref{fig:pair-of-pants} can represent the change in the topology of the phase portrait of the system due to a saddle point in the energy of the system (a Morse function). This change represents a dynamical transition (or 
 bifurcation)~\cite{strogatz2018nonlinear,bhattacharjee2017topology}. An example involves the stabilisation of the inverted pendulum~\cite{kapitza1951}. Analogously, in the Ginzburg-Landau-Wilson approach to critical phenomena~\cite{kardartext}, saddle points in the free energy of a system in which a macroscopic number of degrees of freedom are interacting with one another can signify a phase transition (i.e., a transformation of the phase of the system). Here, the free energy is written in terms of a coarse-grained order parameter determined purely by the symmetries and spatial dimensionality of the system. Another important example involves the case of Lifshitz transitions, i.e., phase transitions that involves changes in the Euler characteristic of the zero-temperature Fermi volume of a system of electrons~\cite{lifshitz1960}. Using the fact that the electronic dispersion relation $E_{\vec{k}}$ (i.e., the energy-momentum relation) for a system of electrons is a Morse function, the Euler characteristic of the Fermi volume (i.e., the set of all electronic states that are occupied at zero temperature) can be obtained from the critical points of $E_{\vec{k}}$ by the application of Morse theory~\cite{vanHove1953}. This Euler characteristic can then be shown to change when the application of external parameters (e.g., pressure, chemical potential etc.) tunes extrema of $E_{\vec{k}}$ (e.g., saddles points, maxima and minima) through the Fermi energy $E_{F}$ (i.e., the greatest occupied energy state, and denoting the Fermi surface in the $D+1$ dimensional energy-momentum space of a $D$-dimensional system of electrons). These extrema are the singularities of the Morse function $E_{\vec{k}}$, and reflect the existence of van Hove singularities in the electronic density of states~\cite{vanHove1953}. While such Lifshitz transitions have typically been studied in systems of non-interacting electrons, a lot of interest in quantum condensed matter physics presently involves searching for similar Lifshitz transitions of the Fermi volume in strongly correlated electronic systems~\cite{volovik2007,volovik2017}.\\
\begin{figure}[!h]
    \centering
    \includegraphics[scale=0.6]{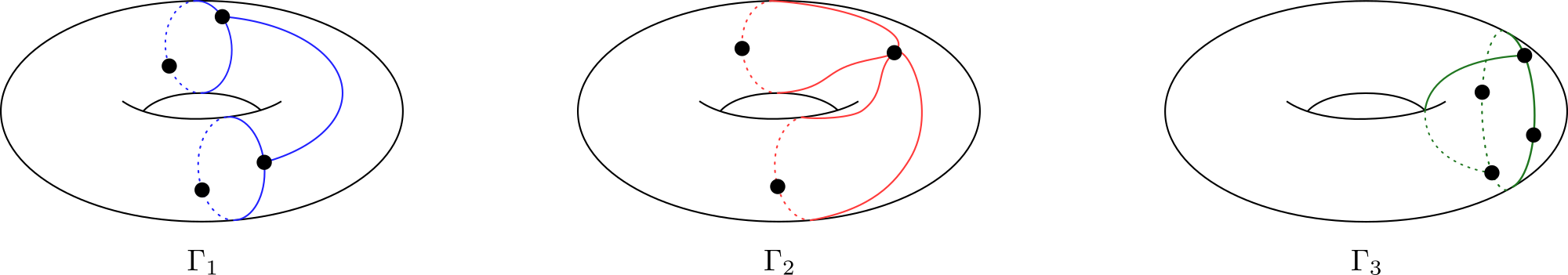}
    \caption{Pair of pants inside a torus}
    \label{fig:transformation_on_torus}
\end{figure}
\hf One can embed the object of figure \ref{fig:pair-of-pants} on a torus as shown in figure \ref{fig:transformation_on_torus}. In that case the island polynomials are given by
\begin{eqnarray*}
\beta_{(\Gamma_1,\sigma)}(x) &=& 4+9x+6x^2+2x^3 \\
\beta_{(\Gamma_2,\sigma)}(x) &=& 3+4x+2x^2\\
\beta_{(\Gamma_3,\sigma)}(x) &=& 4+7x+4x^2+2x^3
\end{eqnarray*}
We obtain the relation
\begin{eqnarray*}
\beta_{(\Gamma_1,\sigma)}(x)-\beta_{(\Gamma_3,\sigma)}(x)=\beta_{(\Gamma_2,\sigma)}(x)-(2x+3)=2x(x+1)
\end{eqnarray*}
The difference in the topology between the embedded graphs $\Gamma_1$ and $\Gamma_3$ is captured by the difference polynomial. For instance, when $x=-3/2$ the difference $\beta_{(\Gamma_1,\sigma)}(-3/2)-\beta_{(\Gamma_3,\sigma)}(-3/2)$ exactly matches with the value $\beta_{(\Gamma_2,\sigma)}(-3/2)$ of the transition configuration. However, for $x=0,-1$ the difference $\beta_{(\Gamma_1,\sigma)}(x)-\beta_{(\Gamma_3,\sigma)}(x)$ vanishes.
\end{eg}
\hf For the general planar case, assume that a graph $\tilde{\Gamma}$ has a local bridge of type I (figure \ref{fig:gen>-<}). 
\begin{figure}[!h]
    \centering
    \includegraphics[scale=0.78]{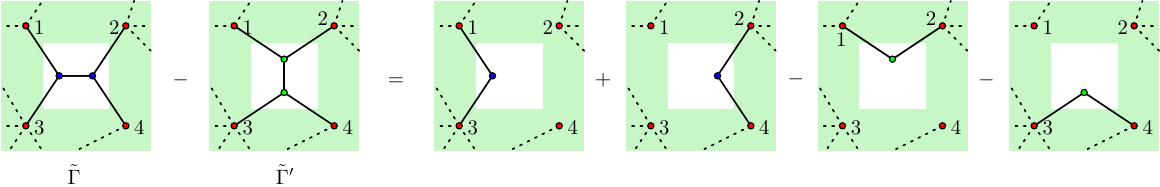}
    \caption{A graphical illustration of equation \eqref{eqn:>-<}}
    \label{fig:gen>-<}
\end{figure}
We make no assumptions regarding the connectivity of $\tilde{\Gamma}$ or of the graph $\Gamma$ obtained by removing the bridge. Let $\tilde{\Gamma}'$ be obtained from $\tilde{\Gamma}$ by removing the type I bridge and inserting a type II bridge. In the generic case, we assume that the four ends of a bridge of type I are four distinct vertices; these are labeled as $1$ through $4$. All the four graphs on the right hand side in figure \ref{fig:gen>-<} are obtained by subdividing an edge that is attached to $\Gamma$ in four possible configurations (figure \ref{fig:spl-att}). 
\begin{figure}[!h]
    \centering
    \includegraphics[scale=0.8]{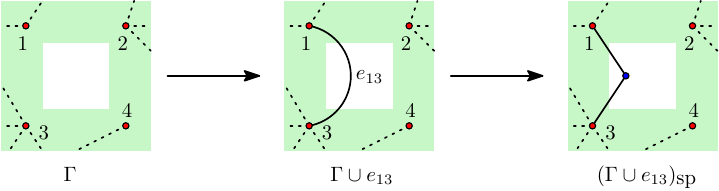}
    \caption{Subdividing after attaching an edge}
    \label{fig:spl-att}
\end{figure}
Let ${\beta}_{ij}$ denote the polynomials of these graphs. From \eqref{eqn:split2}
\begin{equation}\label{eqn:Feyn1}
{\beta}_{ij}(x)=x{\beta}_{\Gamma\cup e_{ij}}(x)+{\beta}_\Gamma(x)+(1+x)^{n-2}.
\end{equation}
As we are working with planar graphs, it follows from comparing coefficients of $\mathcal{D}_k(\Gamma)$ and $\mathcal{D}_k(\Gamma\cup e_{ij})$ we infer that
\begin{equation*}
\mathcal{D}_k(\Gamma\cup e_{ij})=\mathcal{D}_k(\Gamma)+2s^{ij}_k-{n-2 \choose j-2},
\end{equation*}
where $s_k^{ij}$ is the number of subgraphs $\Gamma'$ of $\Gamma$ on $k$ vertices such that vertices $i$ and $j$ belong the same component of $\Gamma'$. In particular, 
\begin{equation}\label{eqn:Feyn2}
{\beta}_{\Gamma\cup e_{ij}}(x)={\beta}_\Gamma(x)+x(1+x)^{n-2}+\sum_k 2s_k^{ij}x^{k-1}.
\end{equation}
Using \eqref{eqn:Feyn2} in \eqref{eqn:Feyn1}, we get
\begin{equation}\label{eqn:Feyn3}
{\beta}_{ij}(x)=2x^2\Big(\sum_k s_k^{ij} x^{k-2}\Big)+(1+x){\beta}_\Gamma(x)+(1+x^2)(1+x)^{n-2}.
\end{equation}
We may now use \eqref{eqn:>-<} and \eqref{eqn:Feyn3} to compute the difference polynomial. It is given by
\begin{equation*}\label{eqn:diff-poly}
{\beta}_{\tilde{\Gamma}}(x)-{\beta}_{\tilde{\Gamma}'}(x)=2x^2\Big(\sum_{k\geq 2} \big(s_k^{13}+ s_k^{24}-s_k^{12}-s_k^{34}\big)x^{k-2}\Big).
\end{equation*}
\begin{rem}
The coefficient $s_k^{13}+ s_k^{24}-s_k^{12}-s_k^{34}$ is reminiscent of an additive version of cross-ratio (in geometry) as well as the sum of the three four-point correlation functions in quantum field theory written in terms of the Mandelstam variables~\cite{peskinschroeder}. 
\end{rem}
\hf The examples and discussions presented here indicate that the graph polynomial is capable of tracking the changes in topology of the embedded graph upon carrying out certain transformations on it. This is of likely relevance to the study of dynamical transitions in non-linear dynamical systems~\cite{strogatz2018nonlinear,bhattacharjee2017topology} and phase transitions in statistical mechanics within the Ginzburg-Landau-Wilson paradigm~\cite{kardartext}. Our formalism is likely to be relevant to the paradigm of fermionic criticality, i.e., the Lifshitz phase transitions of systems of interacting fermions that involve changes in the topology of the Fermi volume~\cite{lifshitz1960,volovik2007,volovik2017,anirbanurg1,anirbanurg2,anirbanmotti,anirbanmott2,MukherjeeTLL,siddharthacpi,MukherjeeKondo2022,mukherjeeMERG2022,pal2019}.


\appendix
\section{Proof of Theorem \ref{thm:Dnm}}\label{cycle-app}
\begin{defn}
	Let $A\subseteq[n]=\{1, 2, \dots, n\}$. Let us call a subset $B\subseteq A$ an \textit{island} if there exist integers $i$ and $j$ such that $B=\{i, i+1, \dots ,j\}$ and $i-1, j+1\notin A$. Basically, $B$ is a connected induced subgraph of $A$. Let $\mathcal I(A)$ denote the set of all islands of $A$ and define $\B^n_m:=\sum_{|A|=m}|\mathcal I(A)|$.
\end{defn}
\begin{defn}
	Let $A\subseteq[n]=\{1, 2, \dots, n\}$. Let us call a subset $B\subseteq A$ an \textit{island on the circle} if there exist integers $i$ and $j$ such that $B=\{i, i+1, \dots ,j\}$ and $i-1, j+1\notin A$. Here the addition is carried out mod $n$. Basically, $B$ is a connected induced subgraph of the cycle graph $C_n$. Let $\mathcal J(A)$ denote the set of all islands of $A$ on the circle and define $\D^n_m=\sum_{|A|=m}|\mathcal J(A)|$.
\end{defn}
Note that $\mathcal{D}^n_m=\mathscr{D}_m(C_n)$. We need some preliminary observations.
\begin{lemma}
	The numbers $\B_m^n$ satisfy the following recurrence relation
	\begin{equation}
	\B_m^n=\B_{m}^{n-1}+\left(\B_{m-1}^{n-2}+{n-2 \choose m-1}\right)+\dots +\left(\B^{n-m}_1+{n-m \choose 1}\right)+1. 
	\label{eq1}
	\end{equation}
\end{lemma}
\begin{proof}
Let us consider a set $A\subset [n]$ of $m$ elements. Let $I$ be the connected induced subgraph in $A$ containing $n$. Then $|I|=r$ for some $0\leq r\leq m$ and the $r=0$ case takes into account that there may be no such connected induced subgraph $I$ in $A$. \\
\hf We now count $\mathcal{B}^n_m$ in terms of $r$. If we set $r=0$, then $I=\varnothing$ and $A\subseteq [n-1]$. This means that sum of all possible islands of all the $A$'s satisfying this case will be $\B^{n-1}_m$. If we put $r=1$, then the only possibility is $I=\{n\}$. Therefore, it follows that $n\in A$ and $n-1\notin A$, whence $A\setminus\{n\}\subseteq[n-2]$ 
As the total number of possible $A$'s is ${n-2 \choose m-1}$, the contribution of this case is given by	
	\bgd
	\sum_{|A|=m, n\in A, n-1\notin A}|\mathcal I(A)| =\B_{m-1}^{n-2}+{n-2 \choose m-1}.
	\edd
Let us now consider the general case, $|I|=r$. In this case the only possibility is $I=\{n-r+1, n-r+2, \dots ,n\}$. It follows that $A\setminus I\subseteq[n-r-1]$. As the total number of possible $A$'s is ${n-r-1 \choose m-r}$, the contribution of this case is given by	
	\bgd
	\sum_{|A|=m, I\subseteq A, n-r\notin A}|\mathcal I(A)| = \B_{m-r}^{n-r-1}+{n-r-1 \choose m-r}.
	\edd
Summing over all possible values of $r$ gives us the required result.
\end{proof}
A repeated use of Pascal's triangle formula helps us simplify \eqref{eq1}:
\begin{equation}
\B_m^n=\B_m^{n-1} +\B_{m-1}^{n-2} + \dots + \B_1^{n-m} + {n-1 \choose m-1}.\label{eq2}
\end{equation}
Using (\ref{eq2}), with $n$ replaced by $n-1$, we obtain
$$\B^{n-2}_{m-1} + B^{n-3}_{m-2} + \dots + \B^{n-m}_1 = \B^{n-1}_{m-1} - {n-2 \choose m-2}.$$
Putting this in (\ref{eq2}) and simplifying, we obtain
\begin{equation}
\B^n_m = \B^{n-1}_m + \B^{n-1}_{m-1} + {n-2 \choose m-1}.
\label{eq3}
\end{equation}
These $\mathcal{B}^n_m$'s obey a concrete formula, i.e., 
\begin{equation}
\B^n_m = (n-m+1){n-1 \choose m-1}.
\label{eq4}
\end{equation}
This follows from the fact that $\B^n_m$ and $(n-m+1){n-1 \choose m-1}$ both satisfy the recurrence relation \eqref{eq3} and agree in the base case, i.e., $\B^n_1=(n-1+1).1=n$.\\
\hf We shall now give a recurrence relation for $\D^n_m$ in terms of $\B^i_j$'s.
\begin{lemma}
The numbers $\D^n_m$ satisfy the following recurrence relation:
\begin{equation}
	\D^n_m=\B^n_m + m+\sum_{j=1}^{m-1} j \left(\B^{n-2-j}_{m-j} + {n-2-j \choose m-j}\right).
\label{eq5}
\end{equation}
\end{lemma}
\begin{proof}
Let us fix a set $A\subseteq[n]$ with $|A|=m$. Let $I$ be the connected induced subgraph in $A$ containing $n$. The possibilities are $|I|=r$ for $0\leq r\leq m$. If $|I|=0$, then $n\notin A$ and $A$ is a subset of the line with $n-2$ elements.	
Hence, the total contribution of this case is $\B^{n-1}_m$. \\
\hf If $|I|=1$, then $I=\{n\}$ and $n-1, 1\notin A$. Hence $A\setminus\{n\}\subseteq[n-2]\setminus\{1\}$ 
Note that the total number of possible $A$'s is ${n-3 \choose m-1}$. Therefore the contribution of this case is given by
\bgd
\sum_{|A|=m, n\in A, n-1, 1\notin A}|\mathcal J(A)| =\B_{m-1}^{n-3}+{n-3 \choose m-1}.
\edd
\hf Let us now consider $|I|=r$, where there are $r$ many distinct possibilities of $I$, all of which are all symmetric, in our context of counting. Thus, we may consider one of them and multiply its contribution by $r$. Let us put $I=\{n-r+1, n-r+2, \dots, n\}$. In particular, as $I\subseteq A$ is an connected induced subgraph we must have $n-r, 1\notin A$. Hence $A\setminus I\subseteq[n-r-1]\setminus\{1\}$. The total number of possible $A$'s is ${n-r-2 \choose m-r}$ and we have
\bgd
\sum_{|A|=m, I\subseteq A, n-r, 1\notin A}|\mathcal J(A)| =\B_{m-r}^{n-r-2}+{n-r-2 \choose m-r}.
\edd
Thus, the total contribution of this case is $r\Big(\B_{m-r}^{n-r-2}+{n-r-2 \choose m-r}\Big)$ and we sum over $r$. 
\end{proof}

\begin{proof}[Proof of Theorem \ref{thm:Dnm}]
Let us prove the second equality. Note that if 
\begin{equation}
\sum_{j=0}^{m-1} (-1)^{j}~ (m-j) {N \choose m-j}=N  {N-2 \choose m-1}\label{id2}
\end{equation}
holds for $N$ and $m<N-1$, then consider the sum (and its simplification)
\begin{eqnarray*}
\sum_{j=0}^{m} (-1)^{j}~ (m+1-j) {N \choose m+1-j} & = & (m+1){N \choose m+1}-\bigg(\sum_{j=0}^{m-1} (-1)^{j}~ (m-j) {N \choose m-j}\bigg) =N {N-2 \choose m}.
\end{eqnarray*}
As \eqref{id2} holds for $m=1$ and any $N$, it holds for all $m$ and $N$ such that $1\leq m<N$. \\
\hf Instead of proving the first equality, we will show that $\mathcal{D}_m(C_N)=N {N-2 \choose m-1}$. By rearranging the terms of right hand side of \eqref{eq5} and repeatedly using \eqref{eq1} we obtain that
\bgd
\D^n_m = \B^n_m - {n-2 \choose m-2}.
\edd
Finally, using \eqref{eq4}, we may simplify the above formula to get our desired formula.
\end{proof}

\printbibliography

\end{document}